\documentclass
[nosumlimits,twoside]{amsart}
\usepackage{amsfonts, amsmath, amssymb}
\usepackage[bookmarksnumbered,plainpages,hypertex]{hyperref}
\usepackage{graphicx}
\usepackage{float}
\usepackage{srcltx}
\usepackage[all]{xy}
\usepackage[usenames]{color}

\newtheorem{theorem}{\sc Theorem}[section]
\newtheorem{proposition}[theorem]{\sc Proposition}

\newtheorem{lemma}[theorem]{\sc Lemma}
\newtheorem{corollary}[theorem]{\sc Corollary}
\theoremstyle{definition}
\newtheorem{definition}[theorem]{\sc Definition}

\theoremstyle{remark}
\newtheorem{remark}[theorem]{\sc Remark}

\newcounter{maint}

\newtheorem{mainthm}[maint]{Theorem}

\def\x1{x_1}
\def\x2{x_2}
\def\a1{a_1}
\def\a2{a_2}

\allowdisplaybreaks

\begin{document}

\title[Gauge deformations]{Gauge deformations for Hopf algebras with\\
the dual Chevalley property}

\author{Alessandro Ardizzoni}
\address{Department of Mathematics, University of Ferrara, Via Machiavelli
35, Ferrara I-44121, Italy}
\email{alessandro.ardizzoni@unife.it}
\urladdr{http://www.unife.it/utenti/alessandro.ardizzoni}

\author{Margaret Beattie}
\address{Department of Mathematics and Computer Science, Mount Allison
University, Sackville, NB, Canada, \indent E4L 1E6}
\email{mbeattie@mta.ca}
\urladdr{http://www.mta.ca/~mbeattie/}

\author{Claudia Menini}
\address{Department of Mathematics, University of Ferrara, Via Machiavelli
35, Ferrara I-44121, Italy}
\email{men@dns.unife.it}
\urladdr{http://www.unife.it/utenti/claudia.menini}

\thanks{M. Beattie's research was supported by an NSERC Discovery Grant.
Thanks to the University of Ferrara for the warm hospitality extended to M.
Beattie during her visit of 2010. This paper was written while A. Ardizzoni
and C. Menini were members of G.N.S.A.G.A. with partial financial support
from M.I.U.R. (PRIN 2007).}

\begin{abstract}
Let $A$ be a Hopf algebra over a field $K$ of characteristic zero such that
its coradical $H$ is a finite dimensional sub-Hopf algebra. Our main theorem
shows that there is a gauge transformation $\zeta $ on $A$ such that $%
A^{\zeta }\cong Q\#H$ where $A^\zeta$ is the dual quasi-bialgebra obtained
from $A$ by twisting its multiplication by $\zeta$, $Q$ is a connected dual
quasi-bialgebra in $^H_H\mathcal{YD}$ and $Q \#H $ is a dual quasi-bialgebra
called the bosonization of $Q$ by $H $.
\end{abstract}

\maketitle
\tableofcontents

\section{Introduction}

Let $A$ be a Hopf algebra over a field $K$ with characteristic $0$ and
suppose that $A$ has the dual Chevalley property, in other words, $H:=A_{0}$%
, the coradical of $A$, is a sub-Hopf algebra of $A$. For example pointed
Hopf algebras have the dual Chevalley property. Moreover we assume that $H$
is finite dimensional. Then if $\lambda \in H^{\ast }$ is the total integral
for $H$, $\lambda $ is left $H$-linear where the left $H$-action of $H$ on
itself is the adjoint action.

Under these conditions there is an $H$-bilinear coalgebra projection $\pi $
from $A$ to $H$ which splits the inclusion so that $A\cong R\#_{\xi }H$
where $(R,\xi )$ is a pre-bialgebra with cocycle in the category $_{H}^{H}%
\mathcal{YD}$ in the sense of \cite{A.M.St.-Small}. If $\xi
=u_{H}\varepsilon _{R\otimes R}$, then $R$ is a bialgebra in $_{H}^{H}%
\mathcal{YD}$ and $A$ is isomorphic to a Radford biproduct. In \cite{A.B.M.}%
, the question of finding a cocycle from $A \otimes A$ to $K$
twisting a Hopf algebra $A$ of the form $R\#_{\xi }H$ to a Radford
biproduct was studied. A bijection between the sets of cocycles on
$A$ and on $R$ shows that this question is equivalent to that of
finding a cocycle ${v}:R\otimes R\rightarrow K$ twisting $(R,\xi )$
to $(R^{{v}},\xi _{{v}}=\varepsilon _{R\otimes R}  )$. The map
${v}=(\lambda \xi )^{-1}$ does give $\xi
_{{v}}=\varepsilon_{R\otimes R}  $ but has not been shown to be a
cocycle; see \cite[Section 5.2]{A.B.M.}.

This paper shows that the correct setting for this problem is that of dual
quasi-bialgebras. Twisting $A$ by ${v}_{A}:A\otimes A\rightarrow K$, where ${%
\ v}_{A}$ is the map $v:=(\lambda \xi)^{-1}$ extended to $A\otimes A$,
yields a dual quasi-bialgebra which is the biproduct of a dual
quasi-bialgebra in $_{H}^{H}\mathcal{YD}$ and $H$.

The following is the main theorem of this paper:

\begin{mainthm}
\label{Theo1}  Let  $A$ be a Hopf algebra over a field of characteristic zero such that the coradical of $A$ is a sub-Hopf algebra.
 Assume the coradical $H$ of $A$ is
finite dimensional so that $H$ is
semisimple. Then there is a gauge transformation $\zeta $ on $A$ such that $%
A^{\zeta }\cong Q\#H$ where

\begin{itemize}
\item $A^\zeta$ is the dual quasi-bialgebra obtained from $A$ by twisting
its multiplication by $\zeta$;

\item $Q$ is a connected dual quasi-bialgebra in $^H_H\mathcal{YD}$;

\item $Q \#H $ is a dual quasi-bialgebra called the bosonization of $Q$ by $%
H $.
\end{itemize}
\end{mainthm}

%As is common with proofs involving dual quasi-bialgebras or
%quasi-bialgebras, the proof of this theorem requires some lengthy
%computations; we have attempted to make these as transparent as possible.
For
$H$  a cosemisimple Hopf algebra over a field $K$, let ${_{H}^{H}%
\mathcal{YD}}$ denote the braided monoidal category of Yetter-Drinfeld
modules over $H$. Throughout this paper, all maps will be $K$-linear.
Essentially the computations in this paper will    construct
functors between different categories and prove that the following diagram  commutes. \smallbreak

 \begin{equation}\label{diagram main}
\xymatrix{
 \mathcal{R} \ar[d]_{T_1} \ar[r]^{B_1} & \mathcal{A} \ar[d]^{T_2}  \\
\mathcal{Q} \ar[r]^{B_2} & \mathcal{B}             }%
\end{equation}%
where
the categories $\mathcal{A},\mathcal{B},%
\mathcal{Q},\mathcal{R}$ are defined as follows.

\begin{itemize}
\item $\mathcal{R}$ has objects the pairs $\left( R,\xi \right) $ where $R$
is a connected pre-bialgebra with cocycle $\xi :R\otimes R\rightarrow H$ in $%
{_{H}^{H}\mathcal{YD}}$ as defined in Subsection \ref{sect: prebialgs}. A morphism $%
f:\left( R,\xi \right) \rightarrow \left( R^{\prime },\xi ^{\prime }\right) $
in $\mathcal{R}$ is a morphism of pre-bialgebras with cocycle in ${%
_{H}^{H}\mathcal{YD}}$, i.e., a coalgebra homomorphism $f:R\rightarrow
R^{\prime }$ in ${_{H}^{H}\mathcal{YD}}$ which is multiplicative and unitary
and such that $\xi ^{\prime }\circ \left( f\otimes f\right) =\xi $ (see \cite%
[Definition 1.10]{SmallApp}).

\item $\mathcal{A}$ has objects the pairs $\left( A,\gamma \right) $ where $%
A $ is a bialgebra with coradical $H$ and $\gamma :A\otimes A\rightarrow K$
is a gauge transformation as defined in Subsection \ref{sec: dualquasiboso}.
A morphism $f:\left( A,\gamma \right) \rightarrow \left( A^{\prime },\gamma
^{\prime }\right) $ in $\mathcal{A}$ is a morphism of bialgebras $%
f:A\rightarrow A^{\prime }$ such that $f_{\mid H} =H$ and $\gamma
^{\prime }\left( f\otimes f\right) =\gamma .$

\item $\mathcal{Q}$ has objects the pairs $\left( Q,\alpha \right) $ where $%
Q $ is a connected dual quasi-bialgebra in ${_{H}^{H}\mathcal{YD}}$ and $%
\alpha \in {_{H}^{H}\mathcal{YD}}\left( Q^{\otimes 3},K\right) $ is the
corresponding reassociator, see Definition \ref{def: dual quasi braided}.
Morphisms in $\mathcal{Q}$ are morphisms of dual quasi-bialgebras   in ${_{H}^{H}\mathcal{YD}}$.

\item $\mathcal{B}$ has objects the pairs $\left( B,\beta \right) $ where $B$
is a dual quasi-bialgebra with coradical $H$ and $\beta :B^{\otimes
3}\rightarrow K$ is
%the associated $3$-cocycle.
\ the reassociator.
 A morphism $f:\left( B,\beta
\right) \rightarrow \left( B^{\prime },\beta ^{\prime }\right) $ in $%
\mathcal{A}$ is a morphism of dual quasi-bialgebras $f:B\rightarrow
B^{\prime }$ such that  $f_{\mid H} =H$, see Subsection \ref{sec:
dualquasiboso}.

\end{itemize}

Assume that $H$ has an ad-invariant integral $\lambda $.     The functors $T_i$ and $B_i$ are twisting and bosonization
functors where

\begin{itemize}
\item $B_{1}\left( R,\xi \right) :=\left( R\#_{\xi }H,v_{R\#_{\xi }H}\right)
$ as defined in Subsection \ref{subs: boso}. %Moreover $B_{1}\left( f\right)
%=f\otimes H$ for every morphism $f$ in $\mathcal{R}$.

\item $B_{2}\left( Q,\alpha \right) :=\left( Q\#H,\alpha _{Q\#H}\right) $ as
defined in Proposition \ref{teo: smash quasi YD}.% Moreover $B_{2}\left(
%f\right) =f\otimes H$ for every morphism $f$ in $\mathcal{Q}$.

\item $T_{1}\left( R,\xi \right) :=\left( R^{v},\partial _{R}^{2}v\right) $
as defined in Proposition \ref{pro: R dual quasi}.% Moreover $T_{1}\left(
%f\right) =f$ for every morphism $f$ in $\mathcal{R}$.

\item $T_{2}\left( A,\gamma \right) :=\left( A^{\gamma },\partial
_{A}^{2}\gamma \right) $ as defined in Subsection \ref{sec: dualquasiboso}.
%Moreover $T_{2}\left( f\right) =f$ for every morphism $f$ in $\mathcal{A}$.
\end{itemize}

Furthermore, for the bosonization functors $B_i$, we have  $B_i(f) = f \otimes H$ for every morphism $f$ in $\mathcal{R}$, $\mathcal{Q}$, and for  the twisting functors $T_i$, we have
$T_i(f) = f$ for every morphism $f$ in $ \mathcal{R}, \mathcal{A}$.

  Proposition \ref{teo: smash bilin} says that $B_{2}\circ
T_{1}=T_{2}\circ B_{1}.$ This result is then applied to prove Theorem \ref%
{Theo1}.

\smallbreak
As is common with proofs involving dual quasi-bialgebras or
quasi-bialgebras,
 some lengthy
computations are involved; we have attempted to make these as transparent as possible. We outline the proof here.
\smallbreak
Sections \ref{sec:preliminaries} and \ref{sect: cohomology} review definitions and set the stage for later computations.
Pre-bialgebras with cocycle, the objects in the category $\mathcal{R}$,
are familiar from \cite{A.M.S.}, \cite{A.M.St.-Small}, \cite{A.B.M.};
 the definitions are reviewed in Subsection \ref{sect: prebialgs}. As well, in Section \ref{sec:preliminaries},
 we review the idea of the   splitting datum $(A, H, \pi, \sigma)$ corresponding
to a pre-bialgebra with cocycle $(R,\xi)$ where $A= R \#_\xi H$ is a kind of  bosonization of $R$ and $H$,
 and remind the reader of some preliminary results from the literature.
In Section \ref{sect: cohomology}, we review the
cohomology of a bialgebra and define cohomology for a pre-bialgebra with cocycle $(R, \xi)$ in terms of the cohomology for the $K$-bialgebra
$A = R \#_\xi H$.  The computations in this section   simplify future calculations
 considerably.
\smallbreak
 For $(R,\xi)$ in $_H^H\mathcal{YD}$, the definition of the multiplication map $m_R: R \otimes R \rightarrow R$
involves the cocycle $\xi: R \otimes R \rightarrow H$. In order to twist $R$ we want to use a map $v: R \otimes R \rightarrow K$ which is related to $\xi$.
The bijection
between the set $\Xi$ of maps $\xi: R \otimes R \rightarrow H$ so that $(R,\xi)$ is a pre-bialgebra with cocycle and
a set  $V \subset \mathrm{Hom}(R \otimes R, K)$ is key to our computations.
This correspondence is established in Section \ref{sect: new interpretation}.  Essentially what happens in this section is that we show
that there is a isomorphism of categories from $\mathcal{R}$ to a category $\mathcal{R}^\prime$ whose objects are pairs $(R,v)$ where $R$ is a connected
pre-bialgebra in $_H^H\mathcal{YD}$ and $v \in V \subset \mathrm{Hom}(R \otimes R, K)$, with $V$ corresponding to $\Xi$.
 Morphisms in this category are multiplicative unitary
 coalgebra homomorphisms $f: R \rightarrow S$ in $_H^H\mathcal{YD}$ such that $v_S \circ (f \otimes f) = v_R$. Then a twisting functor
takes $\mathcal{R}^\prime$ to $\mathcal{Q}$ and $T_1$ is the composite of these.

\smallbreak
In Section \ref{sect: dual quasi}, dual quasi-bialgebras are reviewed, dual quasi-bialgebras
 in a Yetter-Drinfeld category are introduced and
we show that there is a process of  bosonization taking dual quasi-bialgebras in $_H^H\mathcal{YD}$ to dual quasi-bialgebras
in vector spaces. The language of cohomology from Section \ref{sec: coh prebi} is used throughout these computations.
In Subsection \ref{subs: boso} we show that Diagram \eqref{diagram main} commutes. The proof of the main theorem in
Section \ref{sect: main} follows.

\section{Preliminaries}\label{sec:preliminaries}

Throughout we work over a field $K$ of characteristic $0$, and all
maps are assumed to be $K$-linear. We will use Sweedler-Heyneman notation for the
comultiplication in a $K$-coalgebra $C$ but with the summation sign
omitted, namely $\Delta (x)=x_{(1)}\otimes x_{(2)}$ for $x\in C$.
For $C$ a coalgebra and $A$ an algebra the convolution
multiplication in $\mathrm{Hom}(C,A)$ will be denoted $\ast $.
Composition of functions will be denoted by $\circ $ or by
juxtaposition when the meaning is clear.
\smallbreak

 For a Hopf algebra $H$, and $M,N$ left $H$-modules,  we will denote by $\mathrm{Hom}_H(M,N)$ the left $H$-linear maps from $M$ to $N$.
If $M,N$ are $H$-bimodules, then $\mathrm{Hom}_{H,H}$ will denote the $H$-bimodule maps from $M$ to $N$. If $M,N$ are
left  $H$-comodules, then $\mathrm{Hom}^{H}(M,N) $ will denote the left
$H$-colinear maps from $M$ to $N$.  If $M,N$ are $H$-bicomodules, then the bicomodule morphisms from $M$ to $N$ will be denoted
by $\mathrm{Hom}^{H,H}(M,N)$.
\smallbreak If $C$ is a coalgebra and $A$ an algebra,
 $\mathrm{Reg}(C,A)$ denotes the convolution invertible maps
from $C$ to $A$.  If  $C$ is a  left $H$-module coalgebra and $A$ a left $H$-module algebra,
 then the convolution product of left $H$-linear maps from $C$ to $A$ is also
 left $H$-linear so that $\mathrm{Hom}_H(C,A)$ is a submonoid of $\mathrm{Hom}(C,A)$.
  The notation $\mathrm{Reg}_H(C,A)$ will denote the convolution invertible  left $H$-linear maps
 from $C$ to $A$.
Similarly we define $\mathrm{Reg}_{H,H}(C,A)$, $\mathrm{Reg}^{H}(C,A)$, etc.
\medskip

A Hopf algebra $H$ is a left $H$-module
under the adjoint action $h\rightharpoonup m=h_{(1)}m S_H (h_{(2)})$
and has a similar right adjoint action. The symbol $\rightharpoonup
$ may be omitted when the context makes the meaning clear. Recall
\cite[Definition 2.7]{A.M.S.} that a left and right integral
$\lambda \in H^{\ast }$ for $H$ is called
ad-invariant if $\lambda (1)=1$ and $\lambda $ is a left and right $H$%
-module map with respect to the left and right adjoint actions i.e., for all
$h,x\in H$, one has
\begin{equation*}
\lambda \left[ h_{\left( 1\right) }xS_{H}\left( h_{\left( 2\right)
}\right) \right] =\varepsilon _{H}\left( h\right) \lambda \left( x\right)
,\qquad  \lambda \left[ S_{H}\left( h_{\left( 1\right) }\right)
xh_{\left( 2\right) }\right] =\varepsilon _{H}\left( h\right) \lambda \left(
x\right) .
\end{equation*}%
If $H$ is semisimple and cosemisimple (for example,  since  characteristic zero is assumed, if $H$ is finite dimensional
cosemisimple), then the total integral
for $H$ is ad-invariant; see either \cite[Proposition 1.12, b)]{SvO} or \cite%
[Theorem 2.27]{A.M.S.}. If $H$ has an ad-invariant integral, then $H$ is
cosemisimple.

We point out that the group algebra $KG$, which is in general not
semisimple when $K$ has positive characteristic, always admits an $ad$%
-invariant integral, namely the map $\lambda $ defined on $G$ by setting $%
\lambda \left( g\right) =\delta _{g,e}$ (Kronecker delta) and extended by
linearity.\medskip

\textbf{\ Throughout this paper we will assume that $H$ is a $K$-Hopf
algebra with an ad-invariant integral $\lambda $. }Note that in this case $H$
is cosemisimple whence it has bijective antipode, see \cite[Theorem 3.3]%
{Larson}.

We assume familiarity with the general theory of Hopf algebras; good
references are \cite{Sw}, \cite{Mo}.

\subsection{The category $_{H}^{H}\mathcal{YD}$}

Let $(\mathcal{M},\otimes ,\mathbf{1})$ be a monoidal category. A
coaugmented coalgebra $(C,u_{C})$ consists of a coalgebra $C$ in $\mathcal{M}
$ and a coalgebra homomorphism $u_{C}:\mathbf{1}\rightarrow C$ in $\mathcal{M%
},$ called the coaugmentation map. In
particular, a coaugmented coalgebra in the category of vector spaces over a
field $K$ is a $K$-coalgebra $C$ and a $K$-linear map $u_{C}:K\rightarrow C$%
, such that $\varepsilon (1_{C})=1_{K}$ and $\Delta (1_{C})=1_{C}\otimes
1_{C}$ where $1_{C}:=u_{C}(1_{K})$. A coaugmented coalgebra $C$ is called
connected if $C_{0}=K1_{C}$.

\bigskip

Coalgebras in $_{H}^{H}\mathcal{YD}$, the category of left-left
Yetter-Drinfeld modules over $H$, will play a central role in this paper.
For $(V,\cdot )$ a left $H$-module, we write $h{v}$ for $h\cdot {v}$, the
action of $h$ on ${v}$, if the meaning is clear. The left $H$-module $H$
with the left adjoint action is denoted $(H,\rightharpoonup )$; the left and
right actions of $H$ on $H$ induced by multiplication will be denoted by
juxtaposition. For $(V,\rho )$ a left $H$-comodule, we write $\rho ({v})={v}%
_{\langle -1\rangle }\otimes {v}_{\langle 0\rangle }$ for the coaction.
Recall that if $V$ is a left $H$-module and a left $H$-comodule, then $V$ is
an object in $_{H}^{H}\mathcal{YD}$ if
\begin{equation*}
\rho (h\cdot {v})=h_{(1)}{v}_{\langle -1\rangle } S_H
(h_{(3)})\otimes h_{(2)}\cdot {v}_{\langle 0\rangle }.
\end{equation*}%
For example, $(H,\rightharpoonup ,\Delta _{H})$ is a Yetter-Drinfeld module.
The category $_{H}^{H}\mathcal{YD}$ is braided monoidal. The tensor product
of two Yetter-Drinfeld modules $V$ and $W$ is a Yetter-Drinfeld module with
structures $h\left( v\otimes w\right) =h_{\left( 1\right) }v\otimes
h_{\left( 2\right) }w$ and $\rho \left( v\otimes w\right) =v_{\left\langle
-1\right\rangle }w_{\left\langle -1\right\rangle }\otimes v_{\left\langle
0\right\rangle }\otimes w_{\left\langle 0\right\rangle }.$ The unit object
is $K$ with trivial structures i.e. $hk=\varepsilon _{H}\left( h\right) k$
and $\rho \left( k\right) =1_{H}\otimes k.$ The braiding $c_{V,W}:V\otimes
W\rightarrow W\otimes V$ is given by $c_{V,W}({v}\otimes w)={v}_{\langle
-1\rangle }w\otimes {v}_{\langle 0\rangle }$ with inverse $c_{V,W}^{-1}({w}%
\otimes v)=v_{\left\langle 0\right\rangle }\otimes S_{H}^{-1}\left(
v_{\left\langle -1\right\rangle }\right) {w}$.

For $C$ a coalgebra in $_{H}^{H}\mathcal{YD}$, we use a modified version of
the Sweedler notation, writing superscripts instead of subscripts, so that
comultiplication is written
\begin{equation*}
\Delta _{C}(x)=\Delta (x)=x^{(1)}\otimes x^{(2)}\text{, for every }x\in C%
\text{.}
\end{equation*}%
For $C$ a coaugmented coalgebra in $_{H}^{H}\mathcal{YD}$, then
$u_{C}$ is also required to be a map in the Yetter-Drinfeld
category, i.e.,
\begin{equation*}
{h}\cdot 1_{C}=\varepsilon _{H}(h)1_{C}\qquad \text{and}\qquad \rho
_{C}(1_{C})=1_{H}\otimes 1_{C}.  \label{eq:YD0'}
\end{equation*}%
If $C$ and $D$ are coalgebras in $_{H}^{H}\mathcal{YD}$, so is $C\otimes D$
with Yetter-Drinfeld module structure given as above, counit $\varepsilon
_{C\otimes D}=\varepsilon _{C}\otimes \varepsilon _{D}$ and $\Delta
_{C\otimes D}=(C\otimes c_{C,D}\otimes D)\circ (\Delta _{C}\otimes \Delta
_{D})$ , so that
\begin{equation*}
\Delta _{C\otimes D}(x\otimes y)=x^{(1)}\otimes x_{\langle -1\rangle
}^{(2)}y^{(1)}\otimes x_{\langle 0\rangle }^{(2)}\otimes y^{(2)}.
\end{equation*}
If $C$,$D$ are coaugmented so is $C \otimes D$ with $1_{C \otimes D} = 1_C
\otimes 1_D$.

Given three coalgebras $C,D,E$ in $_{H}^{H}\mathcal{YD}$ the coalgebras $%
\left( C\otimes D\right) \otimes E$ and $C\otimes \left( D\otimes E\right) $
are isomorphic and can both be denoted by $C\otimes D\otimes E.$

%with
%comultiplication
%\begin{equation*}
%\Delta _{C\otimes D\otimes E}\left( x\otimes y\otimes z\right) =x^{\left(
%1\right) }\otimes x_{\left\langle -2\right\rangle }^{\left( 2\right)
%}y^{\left( 1\right) }\otimes x_{\left\langle -1\right\rangle }^{\left(
%2\right) }y_{\left\langle -1\right\rangle }^{\left( 2\right) }z^{\left(
%1\right) }\otimes x_{\left\langle 0\right\rangle }^{\left( 2\right) }\otimes
%y_{\left\langle 0\right\rangle }^{\left( 2\right) }\otimes z^{\left(
%2\right) }.
%\end{equation*}%
\medskip

For $(C,\Delta _{C},\varepsilon _{C})$ a left $H$-comodule coalgebra, $%
\Delta _{C}$ and $\varepsilon _{C}$ are morphisms of left $H$-comodules i.e.
for all $c\in C$,%
\begin{eqnarray}
c_{\langle -1\rangle }\otimes (c_{\langle 0\rangle })^{(1)}\otimes
(c_{\langle 0\rangle })^{(2)} &=&(c^{(1)})_{\langle -1\rangle
}(c^{(2)})_{\langle -1\rangle }\otimes (c^{(1)})_{\langle 0\rangle }\otimes
(c^{(2)})_{\langle 0\rangle }\text{ and}  \label{form: Delta colin} \\
c_{\langle -1\rangle }\varepsilon _{C}(c_{\langle 0\rangle })
&=&1_{H}\varepsilon _{C}(c).  \label{form: eps colin}
\end{eqnarray}%
Similarly, for $(C,\Delta _{C},\varepsilon _{C})$ a left $H$-module
coalgebra, $\Delta _{C}$ and $\varepsilon _{C}$ are morphisms of left $H$%
-modules i.e. for all $h\in H,c\in C$,%
\begin{eqnarray}
(hc)^{(1)}\otimes (hc)^{(2)} &=&h_{(1)}c^{(1)}\otimes h_{(2)}c^{(2)}\text{
and}  \label{form: Delta lin} \\
\varepsilon _{C}(hc) &=&\varepsilon _{H}(h)\varepsilon _{C}(c).
\label{form: eps lin}
\end{eqnarray}%
For a left-left Yetter-Drinfeld module $C$, $(C,\Delta _{C},\varepsilon
_{C}) $ is a coalgebra in ${_{H}^{H}\mathcal{YD}}$ if $(C,\Delta
_{C},\varepsilon _{C})$ is both a left $H$-comodule coalgebra and a left $H$%
-module coalgebra.

\begin{lemma}
\label{lem: C1C2}Let $C_{1}$, $D_1$ be left $H$-comodule coalgebras and let $%
C_{2}$, $D_2$ be left $H$-module coalgebras. Then $C:=C_{1}\otimes C_{2}$ is
a coalgebra via%
\begin{equation*}
\Delta _{C}\left( z\otimes w\right) =z^{\left( 1\right) }\otimes
z_{\left\langle -1\right\rangle }^{\left( 2\right) }w^{\left( 1\right)
}\otimes z_{\left\langle 0\right\rangle }^{\left( 2\right) }\otimes
w^{\left( 2\right) }\qquad \text{and}\qquad \varepsilon _{C}\left( z\otimes
w\right) =\varepsilon _{C_{1}}\left( z\right) \varepsilon _{C_{2}}\left(
w\right);
\end{equation*}
$D:=D_1 \otimes D_2$ is a coalgebra in the same way. Let $f:C_{1}\rightarrow
D_{1}$ be a morphism of left $H$-comodule coalgebras and let $%
g:C_{2}\rightarrow D_{2}$ be a morphism of left $H$-module coalgebras. Then $%
f\otimes g:C\rightarrow D$ is a coalgebra map.
\end{lemma}

\begin{proof}
We must show that $(C,\Delta_C, \varepsilon_C)$ is a coalgebra. The
comultiplication map $\Delta_C$ is coassociative since for $y= z \otimes w$
with $z\in C_{1},w\in C_{2},$%
\begin{eqnarray*}
\left( C\otimes \Delta \right) \Delta (y) &=&z^{\left( 1\right) }\otimes
z_{\left\langle -1\right\rangle }^{\left( 2\right) }w^{\left( 1\right)
}\otimes \Delta \left( z_{\left\langle 0\right\rangle }^{\left( 2\right)
}\otimes w^{\left( 2\right) }\right) \\
&=& z^{\left( 1\right) }\otimes z_{\left\langle -1\right\rangle }^{\left(
2\right) }w^{\left( 1\right) }\otimes \left( z_{\left\langle 0\right\rangle
}^{\left( 2\right) }\right) ^{\left( 1\right) }\otimes \left(
z_{\left\langle 0\right\rangle }^{\left( 2\right) }\right) _{\left\langle
-1\right\rangle }^{\left( 2\right) }w^{(2)} \otimes \left( z_{\left\langle
0\right\rangle }^{\left( 2\right) }\right) _{\left\langle 0\right\rangle
}^{\left( 2\right) }\otimes w^{(3)} \\
&\overset{(\ref{form: Delta colin})}{=}& z^{\left( 1\right) }\otimes
z_{\left\langle -1\right\rangle }^{\left( 2\right) }z_{\left\langle
-1\right\rangle }^{\left( 3\right) }w^{\left( 1\right) }\otimes
z_{\left\langle 0\right\rangle }^{\left( 2\right) }\otimes \left(
z_{\left\langle 0\right\rangle }^{\left( 3\right) }\right) _{\left\langle
-1\right\rangle }w^{\left( 2\right) }\otimes \left( z_{\left\langle
0\right\rangle }^{\left( 3\right) }\right) _{\left\langle 0\right\rangle
}\otimes w^{\left( 3\right) } \\
&=& z^{\left( 1\right) }\otimes z_{\left\langle -1\right\rangle }^{\left(
2\right) }z_{\left\langle -2\right\rangle }^{\left( 3\right) }w^{\left(
1\right) }\otimes z_{\left\langle 0\right\rangle }^{\left( 2\right) }\otimes
z_{\left\langle -1\right\rangle }^{\left( 3\right) }w^{\left( 2\right)
}\otimes z_{\left\langle 0\right\rangle }^{\left( 3\right) }\otimes
w^{\left( 3\right) } \\
&\overset{(\ref{form: Delta lin})}{=}&z^{\left( 1\right) }\otimes
z_{\left\langle -1\right\rangle }^{\left( 2\right) }\left( z_{\left\langle
-1\right\rangle }^{\left( 3\right) }w^{\left( 1\right) }\right) ^{\left(
1\right) }\otimes z_{\left\langle 0\right\rangle }^{\left( 2\right) }\otimes
\left( z_{\left\langle -1\right\rangle }^{\left( 3\right) }w^{\left(
1\right) }\right) ^{\left( 2\right) }\otimes z_{\left\langle 0\right\rangle
}^{\left( 3\right) }\otimes w^{\left( 2\right) } \\
&=&\left( z^{\left( 1\right) }\right) ^{\left( 1\right) }\otimes \left(
z^{\left( 1\right) }\right) _{\left\langle -1\right\rangle }^{\left(
2\right) }\left( z_{\left\langle -1\right\rangle }^{\left( 2\right)
}w^{\left( 1\right) }\right) ^{\left( 1\right) }\otimes \left( z^{\left(
1\right) }\right) _{\left\langle 0\right\rangle }^{\left( 2\right) }\otimes
\left( z_{\left\langle -1\right\rangle }^{\left( 2\right) }w^{\left(
1\right) }\right) ^{\left( 2\right) }\otimes z_{\left\langle 0\right\rangle
}^{\left( 2\right) }\otimes w^{\left( 2\right) } \\
&=&\Delta \left( z^{\left( 1\right) }\otimes z_{\left\langle -1\right\rangle
}^{\left( 2\right) }w^{\left( 1\right) }\right) \otimes z_{\left\langle
0\right\rangle }^{\left( 2\right) }\otimes w^{\left( 2\right) }=\left(
\Delta \otimes C \right) \Delta (y).
\end{eqnarray*}%
Also, applying $C\otimes \varepsilon _{C}$ to $\Delta(y) ,$ we obtain
\begin{eqnarray*}
\left( C\otimes \varepsilon_C \right) \Delta ( y) &=&z^{\left(
1\right) }\otimes z_{\left\langle -1\right\rangle }^{\left( 2\right)
}w^{\left( 1\right) }\varepsilon_C \left( z_{\left\langle
0\right\rangle }^{\left( 2\right) }\otimes w^{\left( 2\right)
}\right) =z^{\left( 1\right) }\otimes z_{\left\langle
-1\right\rangle }^{\left( 2\right) }w^{\left( 1\right) }\varepsilon
_{C_{1}}\left( z_{\left\langle 0\right\rangle }^{\left( 2\right)
}\right) \varepsilon _{C_{2}}\left( w^{\left( 2\right) }\right)
\\
&=&z^{\left( 1\right) }\otimes z_{\left\langle -1\right\rangle }^{\left(
2\right) }\varepsilon _{C_{1}}\left( z_{\left\langle 0\right\rangle
}^{\left( 2\right) }\right) w\overset{(\ref{form: eps colin})}{=}y,
\end{eqnarray*}%
and similarly, using (\ref{form: eps lin}), we obtain $\left( \varepsilon_C
\otimes C\right) \Delta_C \left( y\right) =y$. The final statement is clear
since $f,g$ preserve the relevant coalgebra, $H$-module and $H$-comodule
structures.
\end{proof}

\begin{remark}
For $C_{1},C_{2},D_{1},D_{2},f,g$ as in Lemma \ref{lem: C1C2}, since $%
f\otimes g$ is a coalgebra map, then, for $A$ an algebra and maps $\alpha
,\beta :D \rightarrow A$, we have%
\begin{equation}
\left[ \alpha \circ \left( f\otimes g\right) \right] \ast \left[ \beta \circ
\left( f\otimes g\right) \right] =\left( \alpha \ast \beta \right) \circ
\left( f\otimes g\right) .  \label{form:
alfabeta}
\end{equation}
\end{remark}

\vspace{1mm}

\subsection{Some preliminary results}

We recall some key definitions and results from \cite{A.B.M.}.

\begin{definition}
\cite[Definition  2.1]{A.B.M.} \label{def: Psi} For $M$ a left $H$-comodule, define $%
\Psi :\mathrm{Hom}(M,K)\rightarrow \mathrm{Hom}^{H}(M,H)$ by
\begin{equation*}
\Psi \left( \alpha \right) =\left( H\otimes \alpha \right) \rho _{M}.
\end{equation*}
\end{definition}

\begin{remark}
\label{rem: hope1.8} (i) \cite[Remark  2.2]{A.B.M.} Let
$f:M\rightarrow L$ be a
morphism of left $H$-comodules and $\alpha \in \mathrm{Hom}(L,K)$. Then $%
\Psi \left( \alpha \right) \circ f=\Psi \left( \alpha \circ f\right) .$

(ii) \label{lem: Psi} \cite[Lemma 2.3]{A.B.M.} For $C$ a left
$H$-comodule
coalgebra, $\Psi :\mathrm{Hom}\left( C,K\right) \rightarrow \mathrm{Hom}%
^{H}\left( C,H\right) $ is an algebra isomorphism. The inverse $\Psi^{-1}$
is defined by $\Psi^{-1}(\alpha) = \varepsilon_H \alpha$.
\end{remark}

The next lemma notes that $\Psi, \Psi^{-1}$ in the preceding remark are isomorphisms
between the subalgebras $\mathrm{Hom}_H(C,K)$ of  $\mathrm{Hom}(C,K)$ where $K$ has the trivial $H$-action and
$\mathrm{Hom}_H(C,H)$ of $\mathrm{Hom}(C,H)$ where $H$ has the left adjoint $H$-action.

\begin{lemma}
\label{lem: PsiYD}For $C$ a coalgebra in ${_{H}^{H}\mathcal{YD}}$, $\Psi :
\mathrm{Hom}{{
_{H}}}\left( C,K\right) \rightarrow {_{H}^{H}\mathcal{YD}}\left(
C,H\right) $ is an algebra isomorphism.
\end{lemma}

\begin{proof}
Let $\alpha \in \mathrm{Hom}{{
_{H}}}\left( C,K\right) $ and we check that $%
\Psi(\alpha)$ is left $H$-linear.
\begin{eqnarray*}
\Psi \left( \alpha \right) \left( hc\right) &=&\left( H\otimes \alpha
\right) \rho _{M}\left( hc\right) =\left( hc\right) _{\left\langle
-1\right\rangle }\alpha( \left( hc\right) _{\left\langle 0\right\rangle }) \\
&=& h_{\left( 1\right) }c_{\left\langle -1\right\rangle } S_H\left(
h_{\left( 3\right) }\right) \alpha \left[ h_{\left( 2\right)
}c_{\left\langle 0\right\rangle }\right]
\\
&=& h_{\left( 1\right) }c_{\left\langle -1\right\rangle } S_H\left(
h_{\left( 3\right) }\right) \varepsilon _{H}\left( h_{\left(
2\right) }\right) \alpha \left( c_{\left\langle 0\right\rangle
}\right)
\\
&=&h_{\left( 1\right) }c_{\left\langle -1\right\rangle }\alpha
\left( c_{\left\langle 0\right\rangle }\right)  S_H\left( h_{\left(
2\right) }\right)
\\
&=&h_{\left( 1\right) }\left[ \Psi \left( \alpha \right) \left( c\right) %
\right]  S_H\left( h_{\left( 2\right) }\right) =\left[
h\rightharpoonup \Psi \left( \alpha \right) \left( c\right) \right]
.
\end{eqnarray*}%
Similarly if $\beta$ is an $H$-linear map from $C$ to $H$, then
\begin{equation*}
\Psi^{-1}(\beta)(hc) = \varepsilon_H \beta(hc) =
\varepsilon_H(h)\varepsilon_H(\beta(c)) = h \Phi^{-1}(\beta)(c).
\end{equation*}
\end{proof}

The previous result depends on the fact that the forgetful functor from ${_{H}^{H}%
\mathcal{YD}}$ to the category of left $H$-modules is left adjoint to the functor $H\otimes (-)$, see e.g. \cite[Claim 3.3]{Ardi-Separable}.

\begin{definition}
\label{de: dualSweedler 1 cocycle} For $C$ a coalgebra in $_{H}^{H}\mathcal{%
YD}$ and $\alpha \in \mathrm{Hom}(C,H)$, $\alpha $ is called a \emph{%
normalized dual Sweedler 1-cocycle} if $\Delta _{H}\alpha =(m_{H}\otimes
\alpha )(\alpha \otimes \rho _{C})\Delta _{C}$ and $\varepsilon _{H}\alpha
=\varepsilon _{C}$. Equivalently, for $x\in C$,
\begin{equation}
\alpha (x)_{(1)}\otimes \alpha (x)_{(2)}=\alpha (x^{(1)})x_{\langle
-1\rangle }^{(2)}\otimes \alpha (x_{\langle 0\rangle }^{(2)})\quad \text{and}%
\quad \varepsilon _{H}(\alpha (x))=\varepsilon _{C}(x).
\label{eq: Sweedler 1-cocycle}
\end{equation}
\end{definition}

\begin{remark}
\label{rem:ConvSw}Any $\alpha :C\rightarrow H$ satisfying (\ref{eq: Sweedler
1-cocycle}) above is convolution invertible by \cite[Proposition 2.6]{A.B.M.}%
.
\end{remark}

\begin{lemma}
\label{lem: Phi}\cite[Lemma 2.7]{A.B.M.} Let $C$ be a coalgebra and let $%
(M,\mu )$ be a left $H$-module. Define
\begin{equation*}
\Phi :\mathrm{Hom}\left( C,H\right) \rightarrow \mathrm{End}\left( C\otimes
M\right) \text{ by }\Phi \left( \alpha \right) :=\left( C\otimes \mu
_{M}\right) \circ \left[ \left( C\otimes \alpha \right) \Delta _{C}\otimes M%
\right] ,
\end{equation*}%
for $\alpha \in \mathrm{Hom}(C,H)$. The map $\Phi $ is an algebra
homomorphism.
\end{lemma}

The following observation is from \cite{A.B.M.}:

\begin{remark}
\label{rem: Takeuchi}\cite[Remark 2.4]{A.B.M.}

(i) For $C$ a left $H$-module coalgebra and ${v}\in \mathrm{Reg}(C,K)$,
 then ${v}$ is left $H$-linear if and only if ${v}%
^{-1}$ is.

(ii) A result of Takeuchi (see \cite[Lemma 5.2.10]{Mo}) shows that
if $(C,\Delta _{C},\varepsilon _{C},u_{C})$ is a coaugmented
connected coalgebra, then every ${v}:C\rightarrow K$ such that
${v}(1_{C})=1_{K}$ is convolution invertible and
${v}^{-1}(1_{C})=1_{K}$ also.
\end{remark}

\begin{lemma}
\label{lm: 2.12} Let $C$ be a left $H$-comodule coalgebra, $C^\prime $ a
left $H$-module coalgebra and let $C \otimes C^\prime$ have the coalgebra
structure from Lemma \ref{lem: C1C2}. Let $u,v \in \mathrm{Hom}(C,K)$ with $u
$ convolution invertible and let $u^\prime, v^\prime \in {\mathrm{Hom}_H}%
(C^\prime,K)$. %\textcolor[rgb]{0.50,0.00,0.51}{, the left $H$-linear maps from $C^\prime$ to $K$ -omit purple stuff}.
For vector spaces $D,D^\prime$, let $f:C\rightarrow D$ and $g: C^\prime \rightarrow
D^\prime$ and for a $K$-algebra $W$, let $\alpha: D \otimes D^\prime
\rightarrow W$. Then
\begin{equation}
\left( u\otimes u^{\prime }\right) \ast \left[ \alpha \circ (f\otimes g)%
\right] \ast \left( {v}\otimes {v}^{\prime }\right) =\alpha \circ
(f^{u}\otimes u^{\prime }\ast g\ast {v}^{\prime })\circ \Phi \left[ u\ast
\Psi \left( {v}\right) \right]  \label{form: aurea}
\end{equation}%
where $f^{u}:=u\ast f\ast u^{-1}$. In particular,  for $D':=C'$, one
has
\begin{equation}
\left( u\otimes \varepsilon _{C^{\prime }}\right) \ast \left[ \alpha \circ
(f\otimes C^{\prime })\right] =\alpha \circ \left( u\ast f\otimes C^{\prime
}\right)  \label{form: aurina}
\end{equation}%
If $v$ is also left $H$-colinear, we have%
\begin{equation}
\left[ \alpha \circ (f\otimes C^{\prime })\right] \ast \left( {v}\otimes
\varepsilon _{C^{\prime }}\right) =\alpha \circ (f\ast v\otimes C^{\prime })
\label{form: aurina2}
\end{equation}
\end{lemma}

{}

\begin{proof}
For every $z\in C$ and $t\in C^{\prime },$ we have
\begin{eqnarray*}
&&\alpha (f^{u}\otimes u^{\prime }\ast g\ast {v}^{\prime })\Phi \left[ u\ast
\Psi \left( {v}\right) \right] \left( z\otimes t\right) \\
&=&\alpha (f^{u}\otimes u^{\prime }\ast g\ast {v}^{\prime })\left( z^{\left(
1\right) }\otimes \left[ u\ast \Psi \left( {v}\right) \right] ( z^{\left(
2\right) } ) t\right) \\
&=&\alpha (\mathrm{Id}\otimes u^{\prime }\ast g\ast {v}^{\prime })\left(
f^{u} ( z^{\left( 1\right) } ) \otimes u ( z^{\left( 2\right) } )
z_{\left\langle -1\right\rangle }^{\left( 3\right) }{v} ( z_{\left\langle
0\right\rangle }^{\left( 3\right) } ) t\right) \\
&=&\alpha (\mathrm{Id}\otimes u^{\prime }\ast g\ast {v}^{\prime })\left(
f^{u} ( z^{\left( 1\right) } ) u ( z^{\left( 2\right) } ) \otimes
z_{\left\langle -1\right\rangle }^{\left( 3\right) }t\right) {v}\left(
z_{\left\langle 0\right\rangle }^{\left( 3\right) }\right) \\
&=& \alpha (\mathrm{Id}\otimes u^{\prime }\ast g\ast {v}^{\prime })\left( u
( z^{\left( 1\right) } ) f ( z^{\left( 2\right) } ) \otimes z_{\left\langle
-1\right\rangle }^{\left( 3\right) }t\right) {v}\left( z_{\left\langle
0\right\rangle }^{\left( 3\right) }\right) \\
&=& u ( z^{\left( 1\right) } ) \alpha (\mathrm{Id}\otimes u^{\prime }\ast
g\ast {v}^{\prime })\left( f ( z^{\left( 2\right) }) \otimes z_{\left\langle
-1\right\rangle }^{\left( 3\right) }t\right) {v}\left( z_{\left\langle
0\right\rangle }^{\left( 3\right) }\right) \\
&=& u ( z^{\left( 1\right) } ) \alpha \left( f ( z^{\left( 2\right) } )
\otimes \left( u^{\prime }\ast g\ast {v}^{\prime }\right) \left(z_{\left
\langle -1 \right \rangle }^{\left( 3\right) }t\right) \right) {\ v}\left(
z_{\left\langle 0\right\rangle }^{\left( 3\right) }\right) \\
&\overset{(\ref{form: Delta lin})}{=}&u ( z^{\left( 1\right) } ) u^{\prime
}\left( z_{\left\langle -3\right\rangle }^{\left( 3\right) }t^{\left(
1\right) }\right) \alpha \left( f ( z^{\left( 2\right) } ) \otimes g\left(
z_{\left\langle -2\right\rangle }^{\left( 3\right) }t^{\left( 2\right)
}\right) \right) {v}\left( z_{\left\langle 0\right\rangle }^{\left( 3\right)
}\right) {v}^{\prime }\left( z_{\left\langle -1\right\rangle }^{\left(
3\right) }t^{\left( 3\right) }\right) \\
&=&u ( z^{\left( 1\right) } ) u^{\prime }\left( t^{\left( 1\right) }\right)
\alpha \left( f ( z^{\left( 2\right) } ) \otimes g\left( z_{\left\langle
-1\right\rangle }^{\left( 3\right) }t^{\left( 2\right) }\right) \right) {v}%
\left( z_{\left\langle 0\right\rangle }^{\left( 3\right) }\right) {v}%
^{\prime }\left( t^{\left( 3\right) }\right) \text{by $u^{\prime },{v}%
^{\prime }$ left $H$-linear } \\
&=&u ( z^{\left( 1\right) } ) u^{\prime }\left( z_{\left\langle
-1\right\rangle }^{\left( 2\right) }z_{\left\langle -2\right\rangle
}^{\left( 3\right) }t^{\left( 1\right) }\right) \alpha (f\otimes g)\left(
z_{\left\langle 0\right\rangle }^{\left( 2\right) }\otimes z_{\left\langle
-1\right\rangle }^{\left( 3\right) }t^{\left( 2\right) }\right) {v}\left(
z_{\left\langle 0\right\rangle }^{\left( 3\right) }\right) {v}^{\prime
}\left( t^{\left( 3\right) }\right) \text{by $u^{\prime }$ left $H$-linear }
\\
&=& ( u\otimes u^{\prime } ) \left( z^{\left( 1\right) }\otimes
z_{\left\langle -1\right\rangle }^{\left( 2\right) }z_{\left\langle
-2\right\rangle }^{\left( 3\right) }t^{\left( 1\right) }\right) \alpha
(f\otimes g)\left( z_{\left\langle 0\right\rangle }^{\left( 2\right)
}\otimes z_{\left\langle -1\right\rangle }^{\left( 3\right) }t^{\left(
2\right) }\right) \left( {v}\otimes {v}^{\prime }\right) \left(
z_{\left\langle 0\right\rangle }^{\left( 3\right) }\otimes t^{\left(
3\right) }\right) \\
&=&\left[ \left( u\otimes u^{\prime }\right) \ast \alpha (f\otimes g)\ast
\left( {v}\otimes {v}^{\prime }\right) \right] \left( z\otimes t\right),
\end{eqnarray*}%
and (\ref{form: aurea}) is proved. It is then easy to see that  (\ref{form:
aurina}) holds by using (\ref{form: aurea}) with $u^\prime =
\varepsilon_{C^\prime}$, $g = C^\prime$, $v = \varepsilon_C$ and $v^\prime =
\varepsilon_{C^\prime}$.

If $v$ is also left $H$-colinear, we get
\begin{eqnarray*}
&&\left[ \alpha (f\otimes C^{\prime })\ast \left( {v}\otimes \varepsilon
_{C^{\prime }}\right) \right] \left( {z}\otimes {t}\right) \\
&\overset{ (\ref{form: aurea}) }{=}&\alpha (f\otimes C^{\prime
})\Phi \left[ \Psi \left( {v}\right) \right]
\left( {z}\otimes {t}\right) =\alpha (f\otimes C^{\prime })\left( {z}%
_{1}\otimes \left[ \Psi \left( {v}\right) \right] \left( {z}_{2}\right) {t}%
\right) \\
&=&\alpha (f\otimes C^{\prime })\left( {z}_{1}\otimes \left( {z}_{2}\right)
_{-1}{v}\left[ \left( {z}_{2}\right) _{0}\right] {t}\right) =\alpha \left(
f\left( {z}_{1}\right) \otimes \left( {z}_{2}\right) _{-1}{v}\left[ \left( {z%
}_{2}\right) _{0}\right] {t}\right) \\
&=&\alpha \left( f\left( {z}_{\left( 1\right) }\right) \otimes v\left(
z_{\left( 2\right) }\right) {t}\right) =\alpha \left( f\left( {z}_{\left(
1\right) }\right) v\left( z_{\left( 2\right) }\right) \otimes {t}\right) =%
\left[ \alpha (f\ast v\otimes C^{\prime })\right] \left( {z}\otimes {t}%
\right),
\end{eqnarray*}
so that (\ref{form: aurina2}) holds.
\end{proof}

\subsection{Pre-bialgebras with cocycle}

\label{sect: prebialgs} Following \cite[Definition 2.3, Definitions
3.1]{A.M.St.-Small}, we define:

\begin{definition}
A pre-bialgebra $R$ = $(R,m_{R},u_{R},\Delta _{R},\varepsilon _{R})$ in ${_{H}^{H}\mathcal{YD}}$ is a coaugmented coalgebra $\left( R,\Delta
_{R},\varepsilon _{R},u_{R}\right) $ in the category ${_{H}^{H}\mathcal{YD}}$
together with a left $H$-linear map $m_{R}:R\otimes R\rightarrow R$ such
that $m_{R}$ is a coalgebra homomorphism, i.e,
\begin{equation*}
\Delta _{R}m_{R}=(m_{R}\otimes m_{R})\Delta _{R\otimes R}\qquad \text{and}%
\qquad \varepsilon _{R}m_{R}=m_{K}(\varepsilon _{R}\otimes \varepsilon _{R}),
\end{equation*}%
and $u_{R}$ is a unit for $m_{R}$, i.e.,
\begin{equation}
m_{R}(R\otimes u_{R})=R=m_{R}(u_{R}\otimes R).  \label{eq:YD9'}
\end{equation}
\end{definition}

If clear from the context, the subscript $R$ on the structure maps may be
omitted.

Essentially a pre-bialgebra differs from a bialgebra in $_{H}^{H}\mathcal{YD}
$ in that \textbf{the multiplication need not be associative and need not be
a morphism of $H$-comodules}.

\begin{definition}
\label{de: xi} A cocycle for the pre-bialgebra $(R,m,u,\Delta ,\varepsilon )$
in ${_{H}^{H}\mathcal{YD}}$ is a map $\xi :R \otimes R \rightarrow H$ such
that:
\begin{eqnarray}
&& \xi \text{ is left $H$-linear with respect to the left $H$-adjoint action
on $H$};  \label{eq:YD3'} \\
&&\Delta _{H}\xi =(m_{H}\otimes \xi )(\xi \otimes \rho _{R\otimes R})\Delta
_{R\otimes R}\quad \text{and}\quad \varepsilon _{H}\xi =m_{K}(\varepsilon
\otimes \varepsilon );  \label{eq:YD5'} \\
&&c_{R,H}(m_{R}\otimes \xi )\Delta _{R\otimes R}=(m_{H}\otimes m_{R})(\xi
\otimes \rho _{R\otimes R})\Delta _{R\otimes R};  \label{eq:YD6'} \\
&&m_{R}(R\otimes m_{R})=m_{R}(m_{R}\otimes R)\Phi \left( \xi \right) ;
\label{eq:YD7'} \\
&&m_{H}\left( \xi \otimes H\right) \left[ R\otimes \left( m_{R}\otimes \xi
\right) \Delta _{R\otimes R}\right] =m_{H}\left( \xi \otimes H\right) \left(
R\otimes c_{H,R}\right) \left[ \left( m_{R}\otimes \xi \right) \Delta
_{R\otimes R}\otimes R\right] ;  \label{eq:YD8'} \\
&& \xi \mbox{ is unital, i.e., } \xi (R\otimes u)=\xi (u\otimes
R)=\varepsilon 1_{H}.  \label{eq:YD10'}
\end{eqnarray}
Then $(R,m,u,\Delta,\varepsilon, \xi)$, written as $(R,\xi)$ unless more detail is needed, is called a
pre-bialgebra with cocycle in ${_{H}^{H}\mathcal{YD}}$.
\end{definition}

Note that (\ref{eq:YD5'}) says that $\xi $ is a normalized dual Sweedler $1$%
-cocycle as in Definition \ref{de: dualSweedler 1 cocycle}. By Remark \ref%
{rem:ConvSw}, $\xi$ is convolution invertible. Condition (\ref{eq:YD6'}) has
a Yetter-Drinfeld-like form as follows.

\begin{lemma}
\label{lm: YD6' YD} Let $R$ be a pre-bialgebra in ${_{H}^{H}\mathcal{YD}}$,
and let $\xi: R \otimes R \rightarrow H$ be a convolution invertible map.
Write  $C:= R \otimes R$, and $m= m_R$. Then (\ref{eq:YD6'}) holds if and
only if, for all $z\in C$:
\begin{equation}
\rho_R(m(z)) = m( z)_{\langle -1\rangle }\otimes m ( z) _{\langle 0\rangle
}=\xi ( z^{( 1) }) z_{\langle -1\rangle }^{( 2) }\xi ^{-1}( z^{( 3) })
\otimes m ( z_{\langle 0\rangle }^{( 2) })  \label{form:
equivalent to YD6'}
\end{equation}
\end{lemma}

\begin{proof}
Applying (\ref{eq:YD6'}) to $z \in C$, yields $m( z^{( 1) }) _{\langle
-1\rangle }\xi ( z^{( 2) }) \otimes m( z^{( 1) }) _{\langle 0\rangle }=\xi (
z^{( 1) }) z_{\langle -1\rangle }^{( 2) }\otimes m( z_{\langle 0\rangle }^{(
2) }). $ Thus if (\ref{eq:YD6'}) holds, then $m( z^{( 1) }) _{\langle
-1\rangle }\xi ( z^{( 2) }) \otimes m( z^{( 1) }) _{\langle 0\rangle }
\otimes z^{(3)}=\xi ( z^{( 1) }) z_{\langle -1\rangle }^{( 2) }\otimes m(
z_{\langle 0\rangle }^{( 2) }) \otimes z^{(3)}, $ and applying $(m_H \otimes
R )\circ (H \otimes \tau) \circ (H \otimes R \otimes \xi^{-1})$ to both
sides of this equation where $\tau $ is the usual twist map, we obtain (\ref%
{form: equivalent to YD6'}). The argument that (\ref{form: equivalent to
YD6'}) implies (\ref{eq:YD6'}) is similar.
\end{proof}

Condition (\ref{eq:YD7'}) describes associativity of the multiplication $%
m_{R}$; it is shown in \cite[Remark 2.11]{A.B.M.} that if $\xi
(z)t=\varepsilon_{R\otimes R}  (z)t$ for all $z\in R\otimes R,t\in
R$ or, equivalently, if
$\Phi (\xi )=Id_{R^{\otimes 3}}$ then $m_{R}$ is associative. By \cite[%
Theorem 3.7]{A.B.M.} the converse holds if $R$ is connected.

\vspace{1mm}

\subsection{Splitting data}

\label{sec: splitting data} There is a correspondence between pre-bialgebras
with cocycle $(R,\xi)$ in $_H^H\mathcal{YD}$ and the $4$-tuples $%
(A,H,\pi,\sigma)$ known as splitting data.

\begin{definition}
A \textbf{splitting datum} $\left( A,H,\pi ,\sigma \right) $ consists of a
bialgebra $A,$ a Hopf algebra $H,$ a bialgebra homomorphism $\sigma
:H\rightarrow A$ and an $H$-bilinear coalgebra homomorphism $\pi
:A\rightarrow H$ such that $\pi \sigma =\mathrm{Id}_{H}$. Note that $H$%
-bilinear here means $\pi (\sigma (h)x\sigma (h^{\prime }))=h\pi
(x)h^{\prime }$ for all $h,h^{\prime }\in H$ and $x\in A$. We say that a
splitting datum is \textbf{trivial} whenever $\pi $ is a bialgebra
homomorphism.
\end{definition}

Given $(R,\xi)$, a splitting datum is constructed as follows. Let $A:= R
\#_\xi H$ have coalgebra structure equal to the smash coproduct $R\#H$ of $R$
by $H$, i.e., the coalgebra defined on $R\otimes H$ by setting, for every $%
r\in R$ and $h\in H$,%
\begin{equation}
\Delta _{R\#H}\left( r\#h\right) =r^{(1)}\#r_{\langle -1\rangle
}^{(2)}h_{(1)}\otimes r_{\langle 0\rangle }^{(2)}\#h_{\left( 2\right)
},\qquad \varepsilon _{R\#H}\left( r\#h\right) =\varepsilon _{R}\left(
r\right) \varepsilon _{H}\left( h\right).  \label{eq: DeltaA}
\end{equation}

The algebra structures are as follows. The unit $u_{A}(1):=1_{R}\#1_{H}$ and
multiplication is given by
\begin{equation*}
m_{A}=\left( R\otimes m_{H}\right) \left[ \left( m_{R}\otimes \xi \right)
\Delta _{R\otimes R}\otimes m_{H}\right] \left( R\otimes c_{H,R}\otimes
H\right)
\end{equation*}%
so that for $r,s\in R$, $h,l\in H$,
\begin{equation}
m_{A}(r\#h\otimes s\#l)=m_{R}\left( r^{(1)}\otimes r_{\langle
-1\rangle }^{(2)}h_{(1)}s^{(1)}\right) \#\xi \left( r_{\langle
0\rangle }^{(2)}\otimes h_{(2)}s^{(2)}\right) h_{(3)}l. \label{form:
multi smash xi}
\end{equation}

Unless $\xi (R\otimes R)=K$, the action of $\xi (R\otimes R)$ will
not be trivial. It is useful to note
that:
\begin{equation}
(R\otimes \varepsilon _{H})m_{A}(r\#h\otimes s\# l)=m_{R}(r\otimes
hs)\varepsilon_H  (l),\quad (\varepsilon _{R}\otimes
H)m_{A}(r\#h\otimes s\#l)=\xi (r\otimes h_{(1)}s)h_{(2)}l  .
\label{form: varepsilonHmA}
\end{equation}

Note that the canonical injection $\sigma :H\hookrightarrow R\#_{\xi }H$ is
a bialgebra homomorphism. \label{lem: sigma bialg}Furthermore
\begin{equation*}
\pi :R\#_{\xi }H\rightarrow H:r\#h\longmapsto \varepsilon \left( r\right) h
\end{equation*}%
is an $H$-bilinear coalgebra retraction of $\sigma $.\newline

Conversely, suppose that $(A,H,\pi ,\sigma )$ is a splitting datum and we
find $(R,m,u,\Delta ,\varepsilon ,\xi )$, the associated pre-bialgebra with
cocycle in $_{H}^{H}\mathcal{YD}$ as in \cite[2.2.3]{A.B.M.}. As when $\pi $
is a bialgebra morphism and $A$ is a Radford biproduct, set
\begin{equation*}
R=A^{co\pi }=\left\{ a\in A\mid a_{\left( 1\right) }\otimes \pi \left(
a_{\left( 2\right) }\right) =a\otimes 1_{H}\right\} ,
\end{equation*}%
and let
\begin{equation*}
\tau :A\rightarrow R,\hspace{1mm}\tau \left( a\right) =a_{\left(
1\right) }\sigma  S_H \pi \left( a_{\left( 2\right) }\right) .
\end{equation*}%
Define a left-left Yetter-Drinfeld structure on $R$ by
\begin{equation*}
h\cdot r=hr=\sigma \left( h_{\left( 1\right) }\right) r\sigma S_{H}\left(
h_{\left( 1\right) }\right) ,\text{\quad }\rho \left( r\right) =\pi \left(
r_{\left( 1\right) }\right) \otimes r_{\left( 2\right) },
\end{equation*}%
and define a coalgebra structure in $_{H}^{H}\mathcal{YD}$ on $R$ by
\begin{equation}
\Delta (r)=r^{(1)}\otimes r^{(2)}=r_{(1)}\sigma  S_H \pi
(r_{(2)})\otimes
r_{(3)}=\tau \left( r_{\left( 1\right) }\right) \otimes r_{\left( 2\right) },%
\text{\quad }\varepsilon =\varepsilon _{A\mid R}.  \label{form:
comulti R}
\end{equation}%
The map
\begin{equation*}
\omega :R\otimes H\rightarrow A , \quad \omega (r\otimes h)=r\sigma (h)
\end{equation*}%
is an isomorphism of $K$-vector spaces, the inverse being defined by%
\begin{equation*}
\omega ^{-1}:A\rightarrow R\otimes H, \quad \omega ^{-1}(a)=a_{\left(
1\right) }\sigma S_{H}\pi \left( a_{\left( 2\right) }\right) \otimes \pi
\left( a_{\left( 3\right) }\right) =\tau \left( a_{\left( 1\right) }\right)
\otimes \pi \left( a_{\left( 2\right) }\right) .
\end{equation*}%
Clearly $A$ defines, via $\omega $, a bialgebra structure on $R\otimes H$
that will depend on $\sigma $ and $\pi $. As shown in \cite[6.1]{Scha} and
\cite[Theorem 3.64]{A.M.S.}, $(R,m,u,\Delta ,\varepsilon )$ is a
pre-bialgebra in ${_{H}^{H}\mathcal{YD}}$ with cocycle $\xi $ where the maps
$u:K\rightarrow R$ and $m:R\otimes R\rightarrow R$, are defined by%
\begin{equation*}
u=u_{A}^{\mid R},\text{\qquad }m(r\otimes s)=r_{(1)}s_{(1)}\sigma
 S_H \pi (r_{(2)}s_{(2)})=\tau \left( r\cdot _{A}s\right)
\end{equation*}%
and the cocycle $\xi :R\otimes R\rightarrow H$ is the map defined by $\xi
(r\otimes s)=\pi (r\cdot _{A}s). $ Then $(R,\xi )$ is the pre-bialgebra with
cocycle in $_{H}^{H}\mathcal{YD}$ associated to $\left( A,H,\pi ,\sigma
\right) $. Moreover, $\omega :R\#_{\xi }H\rightarrow A$ is a bialgebra
isomorphism.

\subsection{Monoids of $H$-bilinear multibalanced and $H$-linear maps}

%{The isomorphism $\Omega _{H,C}^{n}:\mathrm{Reg}^b_{H,H}\left(
%A^{\otimes n},K\right) \rightarrow \mathrm{Reg}_{H}\left( {C}^{\otimes
%n},K\right) $}
\label{sec: isomorphism}

Let $C$ be a coalgebra in $^H_H\mathcal{YD}$ and set $A:=C\#H$, the smash
coproduct of $C$ by $H$. We point out (cf. \cite[Example 3.17]{A.M.S.}) that
$A$ becomes a coalgebra in the monoidal category $\left( {_{H}^{H}\mathcal{M}%
_{H}^{H}},\square _{H},H\right) $.
%with respect to the structures
%\begin{gather}
%\rho _{C\#H}^{l}\left( c\#h\right) =c_{\left\langle -1\right\rangle
%}h_{\left( 1\right) }\otimes (c_{\left\langle 0\right\rangle }\#h_{\left(
%2\right) }),\qquad \rho _{C\#H}^{r}\left( c\#h\right) =(c\#h_{\left(
%1\right) })\otimes h_{\left( 2\right) },  \label{form: coactionCsmashH} \\
%k\left( c\#h\right) =k_{\left( 1\right) }c\#k_{\left( 2\right) }h,\qquad
%\left( c\#h\right) k=c\#hk.  \label{form: actionCsmashH}
%\end{gather}%
%for all $c\in C,h,k\in H.$ In particular $\Delta _{A}$ and $\varepsilon _{A}$
%are $H$-bilinear i.e. for $a\in A$,
%\begin{equation}
%\Delta _{A}(k{a}l)=k_{\left( 1\right) }a_{\left( 1\right) }l_{(1)}\otimes
%k_{\left( 2\right) }a_{\left( 2\right) }l_{(2)},\quad \varepsilon _{A}\left(
%kal\right) =\varepsilon _{H}\left( k\right) \varepsilon _{A}\left( {a}%
%\right) \varepsilon _{H}\left( l\right) .  \label{form: DaltaHbilin}
%\end{equation}%

It is clear that we can regard $A^{\otimes n}$ as an $H$-bimodule via the
structures of the first left (resp. right) hand-side factor. Regard $%
C^{\otimes n}$ as a left $H$-module via the diagonal action. For $n>1$, a map $f:A^{\otimes n}\rightarrow K$ is called $H$-multibalanced
%\textcolor[rgb]{0.50,0.00,0.51}{[if it is $H$-bilinear - suggest we delete --
% maybe it is better to keep the two conditions of
%bilinearity and balanced separate -- we can change back if you prefer]}
if  %\textcolor[rgb]{0.50,0.00,0.50}{, for $n>1$,}
for all ${a}^{1},\ldots ,{a}^{n}\in A,h\in H$, one has $f\left( {a}^{1}\otimes \cdots \otimes a^{i}h\otimes
a^{i+1}\otimes \cdots \otimes {a}^{n}\right) =f\left( {a}^{1}\otimes \cdots
\otimes a^{i}\otimes ha^{i+1}\otimes \cdots \otimes {a}^{n}\right) $ for $1\leq i\leq
n-1 $. Sets of multibalanced maps will be denoted by a superscript $^b$. For example
\begin{displaymath}
\mathrm{Hom}_{H,H}^{b}\left( A^{\otimes n},K\right) =
 \{ f| f \in \mathrm{Hom}_{H,H}\left( A^{\otimes n},K\right)\text{ and,  if  } n>1, f \text{  is  }  H\text{-multibalanced} \}
\end{displaymath}
is the submonoid of $\mathrm{Hom}_{H,H}\left( A^{\otimes n},K\right)$ of $H$-multibalanced maps.

%\textcolor[rgb]{0.50,0.00,0.51}{I suggest that we delete:  We will denote the space of $H$-multibalanced   $H$-bilinear maps from $A^{\otimes n}$ to $K$ by
%$$\mathrm{Hom}_{H,H}^{b}\left( A^{\otimes n},K\right).$$
%%\begin{equation*}
%%\mathrm{Hom}_{H}\left( C^{\otimes n},K\right) =\left\{ f\in \mathrm{Hom}%
%%\left( {C}^{\otimes n},K\right) \mid f\text{ is left }H\text{-linear}%
%%\right\} .
%%\end{equation*}%
%For any coalgebra $\left( D,\Delta ,\varepsilon \right) $ and algebra $%
%\left( A,m,u\right) $ we set $\mathrm{Reg}\left( D,A\right) :=U\left(
%\mathrm{Hom}_{K}\left( D,A\right) \right) $, i.e. the group of all units in
%the monoid $\left( \mathrm{Hom}_{K}\left( D,A\right) ,\ast ,u\varepsilon
%\right) $. We also set
%\begin{eqnarray*}
%\mathrm{Reg}^b_{H,H}\left( A^{\otimes n},K\right) &:&=\mathrm{Hom}%
%_{H,H}^{b}\left( A^{\otimes n},K\right) \cap \mathrm{Reg}\left( A^{\otimes
%n},K\right) , \\
%\mathrm{Reg}_{H}\left( {C}^{\otimes n},K\right) &:&=\mathrm{Hom}_{H}\left(
%{C}^{\otimes n},K\right) \cap \mathrm{Reg}\left( {C}^{\otimes n},K\right)
%\text{.}
%\end{eqnarray*}}

\begin{lemma}
\label{lem:Omega} For $n\in \mathbb{N} $, there is an isomorphism of monoids

%\begin{itemize}
%\item[i)] $( \mathrm{Hom}_{H,H}^{b}\left( A^{\otimes n},K) ,\ast
%,\varepsilon _{A^{\otimes n}}\right) $ is a submonoid of $\left( \mathrm{Hom}%
%\left( A^{\otimes n},K\right) ,\ast ,\varepsilon _{A^{\otimes n}}\right) $,
%
%\item[ii)] $\left( \mathrm{Hom}_{H}\left( {C}^{\otimes n},K\right) ,\ast
%,\varepsilon _{{C}^{\otimes n}}\right) $ is a submonoid of $\left( \mathrm{%
%Hom}\left( {C}^{\otimes n},K\right) ,\ast ,\varepsilon _{{C}^{\otimes
%n}}\right) .$
%
%\item[iii)]

\begin{equation*}
\Omega ^{n}=\Omega _{H,C}^{n}:( \mathrm{Hom}_{H,H}^{b}\left( A^{\otimes
n},K\right) ,\ast ,\varepsilon _{A^{\otimes n}}) \rightarrow ( \mathrm{Hom}%
_{H}\left( {C}^{\otimes n},K\right) ,\ast ,\varepsilon _{{C}^{\otimes n}})
\end{equation*}
defined by $\gamma \mapsto \Omega ^{n}\left( \gamma \right) $ where
\begin{equation*}
\Omega ^{n}\left( \gamma \right) \left( {c}^{1}\otimes {c}^{2}\otimes \cdots
\otimes {c}^{n-1}\otimes {c}^{n}\right) :=\gamma \left( \left( {c}%
^{1}\#1_{H}\right) \otimes \left( {c}^{2}\#1_{H}\right) \otimes \cdots
\otimes \left( {c}^{n-1}\#1_{H}\right) \otimes \left( {c}^{n}\#1_{H}\right)
\right) ,
\end{equation*}%
with inverse $\mho ^{n}=\mho _{H,C}^{n}$ given by $v\mapsto \mho ^{n}\left(
v\right) $ where
\begin{eqnarray*}
&&\mho ^{n}\left( v\right) \left( \left( {c}^{1}\#h^{1}\right) \otimes
\left( {c}^{2}\#h^{2}\right) \otimes \cdots \otimes \left( {c}%
^{n-1}\#h^{n-1}\right) \otimes \left( {c}^{n}\#h^{n}\right) \right) \\
&&:=v\left( {c}^{1}\otimes h_{\left( 1\right) }^{1}{c}^{2}\otimes \cdots
\otimes h_{\left( n-2\right) }^{1}\cdots h_{\left( 2\right) }^{n-3}h_{\left(
1\right) }^{n-2}{c}^{n-1}\otimes h_{\left( n-1\right) }^{1}\cdots h_{\left(
3\right) }^{n-3}h_{\left( 2\right) }^{n-2}h^{n-1}{c}^{n}\right) \varepsilon
_{H}\left( h^{n}\right) .
\end{eqnarray*}
%\end{itemize}
\end{lemma}

\begin{proof}
The proof is straightforward, cf. \cite[Proposition 4.9]{A.B.M.}.
\end{proof}

Clearly the isomorphism $\Omega ^{n}$ induces an isomorphism%
\begin{equation*}
\Omega _{H,C}^{n}:\mathrm{Reg}^b_{H,H}\left( A^{\otimes n},K\right)
\rightarrow \mathrm{Reg}_{H}\left( {C}^{\otimes n},K\right)
\end{equation*}%
with inverse $\mho _{H,C}^{n}.$

Assume furthermore that $C$ is coaugmented and call a map $\phi$ in $\mathrm{Hom}(A^{\otimes n},K)$ or in $\mathrm{Hom}%
(C^{\otimes n},K)$ unital if $\phi = \varepsilon$ on elements of the form $%
x^1 \otimes \ldots \otimes x^n$ with at least one of the $x^i$ equal to $1$.
If $\phi$ is unital and convolution invertible, then $\phi^{-1}$ is also
unital. It is easy to see that $\Omega^n $ and $\mho^n$ preserve unitality.

\section{Cohomology of (pre-)bialgebras}\label{sect: cohomology}

\subsection{Cohomology of a $K$-bialgebra} \label{subs: coalgmultunit}

Recall (cf. \cite[dual to page 368]{Kassel}) that a coalgebra with
multiplication and unit is a datum $\left( E,m,u,\Delta ,\varepsilon \right)
$ where $\left( E,\Delta ,\varepsilon \right) $ is a $K$-coalgebra, $%
m:E\otimes E\rightarrow E$ is a coalgebra homomorphism called multiplication
(which may fail to be associative) and $u:K\rightarrow E$ is a coalgebra
homomorphism called unit. Let $\left( {E},m,u,\Delta ,\varepsilon \right) $
be a coalgebra with multiplication and unit. For $t\in \mathbb{N}$ and $%
0\leq i\leq t+1$ define the maps%
\begin{equation*}
m_{i}^{t+1}:{E}^{\otimes (t+1)}\rightarrow {E}^{\otimes t}
\end{equation*}%
as follows. If $t=0$ we set $m_{0}^{1}=\varepsilon =m_{1}^{1}$ while for $%
t>0 $ we set
\begin{eqnarray*}
& m_{i}^{t+1} (x_{1}\otimes \cdots \otimes x_{i}\otimes \cdots \otimes
x_{t+1})=&%
\begin{cases}
\varepsilon \left( x_{1}\right) x_{2}\otimes \cdots \otimes x_{i}\otimes
\cdots x_{t+1} & \text{for }i=0, \\
x_{1} \otimes \cdots \otimes x_{i}x_{i+1}\otimes \cdots \otimes x_{t+1} &
\text{for }1\leq i\leq t, \\
x_{1} \otimes \cdots \otimes x_{i} \otimes \cdots \otimes x_{t}\varepsilon
\left( x_{t+1}\right) & \text{for }i=t+1.%
\end{cases}%
\end{eqnarray*}
For ${w}\in \mathrm{Reg}\left( {E}^{\otimes t},K\right) $, note that $w
m_i^{t+1} \in \mathrm{Reg}\left( {E}^{\otimes (t+1)},K\right) $ since $%
m_i^{t+1} $ is a coalgebra map so that ${w}m_{i}^{t+1}\ast {w}%
^{-1}m_{i}^{t+1}=\left( {w}\ast {w}^{-1}\right) m_{i}^{t+1}=\varepsilon _{{E}%
^{\otimes t}}m_{i}^{t+1}=\varepsilon _{{E}^{\otimes (t+1)}}$ and thus%
\begin{equation*}
{w}^{-1}m_{i}^{t+1}=\left( {w}m_{i}^{t+1}\right) ^{-1}.
\end{equation*}

For $w \in \mathrm{Reg}(E^{\otimes t},K)$, define
(cf. \cite[pages 60, 53]{Majid}) the two elements $\partial^t_+ (w) $ and $%
\partial^t_-  (w) $ in $\mathrm{Reg}\left( {E}^{\otimes (t+1)},K\right) $ to be
the convolution products:
\begin{equation*}
\left( \partial _{E}^{t}\right) _{+}\left( {w}\right) =\partial
_{+}^{t}\left( {w}\right) =\prod_{\substack{ i=0,...,t+1  \\ i\text{ even}}}{%
w} m_{i}^{t+1},\qquad \left( \partial _{E}^{t}\right) _{-}\left( {w}%
\right) =\partial _{-}^{t}\left( {w}\right) =\prod_{\substack{ i=0,...,t+1
\\ i\text{ odd}}}{w} m_{i}^{t+1}.
\end{equation*}

In particular, for $t=0,1,2,3$:
\begin{eqnarray*}
t=0: &&\partial _{+}^{0}\left( {w}\right) ={w}m_{0}^{1}={w}\varepsilon ;
\hspace{2.5mm} \partial _{-}^{0}\left( {w}\right) ={w}m_{1}^{1}={w}%
\varepsilon, \\
 t=1: &&   \partial _{+}^{1}\left( {w}\right) ={w}m_{0}^{2}\ast {w}%
m_{2}^{2}=m_{K}\left( \varepsilon \otimes {w}\right) \ast m_{K}\left( {w}
\otimes \varepsilon \right)=m_{K}\left( {w}\otimes w \right),   \\
&&   \partial _{-}^{1}\left( {w}\right) ={w}m_{1}^{2}=  {w}m  ,
\\
t=2: && \partial _{+}^{2}\left( {w}\right) ={w}m_{0}^{3}\ast {w}%
m_{2}^{3}=m_{K}\left( \varepsilon \otimes {w}\right) \ast {w}\left( {E}%
\otimes m\right) , \\
&& \partial _{-}^{2}\left( {w}\right) ={w}m_{1}^{3}\ast {w}m_{3}^{3}={w}%
\left( m\otimes {E}\right) \ast m_{K}\left( {w}\otimes \varepsilon \right) ,
\\
t=3: && \partial _{+}^{3}\left( {w}\right) ={w}m_{0}^{4}\ast {w}%
m_{2}^{4}\ast {w}m_{4}^{4}=m_{K}\left( \varepsilon \otimes {w}\right) \ast {w%
}\left( {E}\otimes m\otimes {E}\right) \ast m_{K}\left( {w}\otimes
\varepsilon \right) , \\
&& \partial _{-}^{3}\left( {w}\right) ={w}m_{1}^{4}\ast {w}m_{3}^{4}={w}%
\left( m\otimes {E}\otimes {E}\right) \ast {w}\left( {E}\otimes {E}\otimes
m\right) .
\end{eqnarray*}

Now define the maps:
\begin{equation*}
\partial ^{t}=\partial _{{E}}^{t}:\mathrm{Reg}\left( {E}^{\otimes
t},K\right) \rightarrow \mathrm{Reg}\left( {E}^{\otimes (t+1)},K\right) :{w}%
\mapsto \partial _{+}^{t}\left( {w}\right) \ast \partial _{-}^{t}\left( {w}%
^{-1}\right) ,
\end{equation*}

and by the definition of $\partial ^{t}$:
\begin{eqnarray*}
&&  \partial ^{0}\left( {w}\right) =\varepsilon ,
\\
&&  \partial ^{1}\left( {w}\right) =m_{K}\left( {w}\otimes {w}\right) \ast {w}%
^{-1}m,
\\
&&  \partial ^{2}\left( {w}\right) =m_{K}\left( \varepsilon \otimes
{w}\right) \ast {w}\left( {E}\otimes m\right) \ast {w}^{-1}\left(
m\otimes {E}\right) \ast m_{K}\left( {w}^{-1}\otimes \varepsilon
\right) ,
\\
&&\partial ^{3}\left( {w}\right) =m_{K}\left( \varepsilon \otimes
{w}\right) \ast {w}\left( {E}\otimes m\otimes {E}\right) \ast
m_{K}\left( {w}\otimes
\varepsilon \right) \ast {w}^{-1}\left( m\otimes {E}\otimes {E}\right) \ast {%
w}^{-1}\left( {E}\otimes {E}\otimes m\right) .
\end{eqnarray*}

By definition, a $t$-cocycle is an element ${w}\in \mathrm{Reg}\left( {E}%
^{\otimes t},K\right) $ such that $\partial ^{t}\left( {w}\right) = \varepsilon_{E^{\otimes (t+1)}}
% = 1_{%
%\mathrm{Reg}\left( {E}^{\otimes ( t+1)},K\right) }
$ so that $Z^2(A,K)$ (\cite%
{Doi-braided},\cite[Section 4]{A.B.M.}) is just $\mathrm{Ker}(\partial^2)$.
A $t$-coboundary is an element in $\mathrm{Im}\left( \partial ^{t-1}\right) $%
.

In general it is not clear whether $\partial ^{t}$ is a group homomorphism
and any $t$-coboundary is a $t$-cocycle. This holds for some ${E}$; for example it is true if  $%
{E}$ is a cocommutative bialgebra.  In general, $\partial ^{1}\partial
^{0}\left( {w}\right) =\varepsilon _{E\otimes E}.$ Moroever, if $m$ is
associative, then $\partial ^{2}\partial ^{1}\left( {w}\right) =\varepsilon
_{E^{\otimes 3}}$, see \cite[Lemma 1]{BichonCarnovale}.

\subsection{Cohomology of a pre-bialgebra with cocycle}

\label{sec: coh prebi}

Let $(R,\xi )$ be a connected pre-bialgebra with
cocycle\emph{\ }in ${_{H}^{H}\mathcal{YD}}$ with associated splitting datum $%
(A:=R\#_{\xi }H,H,\pi ,\sigma )$ as outlined in the preliminaries. Since $A$
is a bialgebra, we can consider the maps%
\begin{equation*}
\left( \partial _{{A}}^{t}\right) _{+},\left( \partial _{{A}}^{t}\right)
_{-},\partial _{{A}}^{t}:\mathrm{Reg}\left( {A}^{\otimes t},K\right)
\rightarrow \mathrm{Reg}\left( {A}^{\otimes (t+1)},K\right) ,
\end{equation*}%
where $\partial _{{A}}^{t}\left( {w}\right) =\partial _{+}^{t}\left( {w}%
\right) \ast \partial _{-}^{t}\left( {w}^{-1}\right) $ as in the previous
section. Since the multiplication of $A$ is $H$-bilinear and $H$-balanced,
(cf. \cite[Theorem 3.62]{A.M.S.})  it is clear that these  %maps
induce maps%
\begin{equation*}
\left( \partial _{{A}}^{t}\right) _{+},\left( \partial _{{A}}^{t}\right)
_{-},\partial _{{A}}^{t}:\mathrm{Reg}^b_{H,H}\left( {A}^{\otimes t},K\right)
\rightarrow \mathrm{Reg}^b_{H,H}\left( {A}^{\otimes (t+1)},K\right) .
\end{equation*}

Since $(R,\Delta ,\varepsilon ,m,u)$ is not a coalgebra with
multiplication   and unit  over $K$   in the sense of Subsection
\ref{subs: coalgmultunit} (the coalgebra structure on $R\otimes R$
is different as here it depends on the braiding of
$^H_H\mathcal{YD}$), the definition of cohomology given in the
previous section does
not apply to $R$. In order to define a suitable cohomology for $\mathrm{Reg}%
(R^{\otimes t},K)$ we use the isomorphisms $\Omega$ and $\mho$ from Section %
\ref{sec: isomorphism}.

Note that this is the same approach taken in \cite[Section 4]{A.B.M.} where $%
Z^2_H(R,K)$ is defined to be the set of $H$-linear unital maps $\nu$ from $R
\otimes R$ to $K$ satisfying:
\begin{equation*}
(\varepsilon_R \otimes \nu) \ast \nu (R \otimes m_R) = (\nu \otimes
\varepsilon_R) \ast [\nu(m_R \otimes R) \Phi(\xi)],
\end{equation*}
and then $Z^2_H(R,K) = \Omega^2(Z^2_H(A,K))$; see \cite[Theorem 4.10]{A.B.M.}.

\begin{definition}
Cohomology for $\mathrm{Reg}(R^{\otimes t},K)$ is defined in terms of $
\partial _{{A}}^{t}   $ by:
%\begin{eqnarray*}
%&&
\begin{displaymath}{\left( \partial _{R}^{t}\right) _{+} := \Omega _{H,R}^{t+1}\circ \left(
\partial _{{A}}^{t}\right) _{+}\circ \mho _{H,R}^{t}, \hspace{3mm} \left(
\partial _{R}^{t}\right) _{-} :=\Omega _{H,R}^{t+1}\circ \left( \partial _{{A%
}}^{t}\right) _{-}\circ \mho _{H,R}^{t},
%&& \hspace{2mm}
\text{ and } \hspace{2mm} \partial _{R}^{t} := \Omega
_{H,R}^{t+1}\circ \partial _{{A}}^{t}\circ \mho _{H,R}^{t}.}
\end{displaymath}
%\end{eqnarray*}
\end{definition}

Thus the diagram below commutes as do similar diagrams with $\partial$
replaced by $\partial_+$ or $\partial_-$.
\begin{equation*}
\begin{array}{cccccccc}
\mathrm{Reg}_{H,H}^b\left( K,K\right) & \overset{\partial _{A}^{0}}{%
\longrightarrow } & \mathrm{Reg}_{H,H}\left( A,K\right) & \overset{\partial
_{A}^{1}}{\longrightarrow } & \mathrm{Reg}^b_{H,H}\left( A^{\otimes 2},K\right)
& \overset{\partial _{A}^{2}}{\longrightarrow } & \mathrm{Reg}^b_{H,H}\left(
A^{\otimes 3},K\right) & \cdots \\
\Omega ^{0}\downarrow \uparrow \mho ^{0} &  & \Omega ^{1}\downarrow \uparrow
\mho ^{1} &  & \Omega ^{2}\downarrow \uparrow \mho ^{2} &  & \Omega
^{3}\downarrow \uparrow \mho ^{3} &  \\
\mathrm{Reg}_{H}\left( K,K\right) & \overset{\partial _{R}^{0}}{%
\longrightarrow } & \mathrm{Reg}_{H}\left( R,K\right) & \overset{\partial
_{R}^{1}}{\longrightarrow } & \mathrm{Reg}_{H}\left( R^{\otimes 2},K\right)
& \overset{\partial _{R}^{2}}{\longrightarrow } & \mathrm{Reg}_{H}\left(
R^{\otimes 3},K\right) & \cdots%
\end{array}%
\end{equation*}
\bigbreak
\begin{lemma}
For any $w\in \mathrm{Reg}_{H}\left( {R}^{\otimes t},K\right) $, $\partial
_{R}^{t}\left( {w}\right) =\left( \partial _{R}^{t}\right) _{+}\left( {w}%
\right) \ast \left( \partial _{R}^{t}\right) _{-}\left( {w}^{-1}\right) .$
\end{lemma}

\begin{proof}
We have%
\begin{eqnarray*}
&&\left( \partial _{R}^{t}\right) _{+}\left( {w}\right) \ast \left( \partial
_{R}^{t}\right) _{-}\left( {w}^{-1}\right) =\left[ \Omega _{H,R}^{t+1}\left(
\partial _{{A}}^{t}\right) _{+}\mho _{H,R}^{t}\left( {w}\right) \right] \ast %
\left[ \Omega _{H,R}^{t+1}\left( \partial _{{A}}^{t}\right) _{-}\mho
_{H,R}^{t}\left( {w}^{-1}\right) \right] \\
&=&\Omega _{H,R}^{t+1}\left\{ \left[ \left( \partial _{{A}}^{t}\right)
_{+}\mho _{H,R}^{t}\left( {w}\right) \right] \ast \left[ \left( \partial _{{A%
}}^{t}\right) _{-} (\left( \mho _{H,R}^{t}\left( {w}\right) \right) ^{-1}) %
\right] \right\} =\Omega _{H,R}^{t+1}\partial _{{A}}^{t}\mho
_{H,R}^{t}\left( {w}\right) =\partial _{R}^{t}\left( {w}\right) .
\end{eqnarray*}
\end{proof}

In the following sections of this paper we will need formulas for $%
(\partial^2_R)_\pm$.

\begin{proposition}
\label{pro:DeltaR}Let $w\in \mathrm{Reg}_{H}\left( {R}^{\otimes t},K\right) $
with $t=0,1,2.$

\begin{enumerate}
\item[i)] For $t=0,$ we have
\begin{equation*}
\left( \partial _{R}^{{0}}\right) _{+}\left( w\right) =w\varepsilon
_{R},\qquad \left( \partial _{R}^{{0}}\right) _{-}\left( w\right)
=w\varepsilon _{R}\qquad \text{and}\qquad \partial _{R}^{{0}}\left( {w}%
\right) =\varepsilon _{R}.
\end{equation*}

\item[ii)] For $t=1,$
\begin{eqnarray*}
&&\left( \partial _{R}^{{1}}\right) _{+}\left( w\right) =m_{K}\left( {w}%
\otimes {w}\right) ,\qquad \left( \partial _{R}^{{1}}\right) _{-}\left(
w\right) =wm_{R}\qquad \text{and} \\
&&\partial _{R}^{{1}}\left( {w}\right) =\left[ m_{K}\left( {w}\otimes {w}%
\right) \right] \ast \left( w^{-1}m_{R}\right) .
\end{eqnarray*}

\item[iii)] For $t=2,$
\begin{eqnarray*}
\left( \partial _{R}^{{2}}\right) _{+}\left( w\right) &=&\left[ m_{K}\left(
\varepsilon _{R}\otimes w\right) \right] \ast \left[ w\left( R\otimes
m_{R}\right) \right] =w(R\otimes m^w)\ast   [m_K\left( \varepsilon _{R}\otimes {w}\right)]  , \\
\left( \partial _{R}^{{2}}\right) _{-}\left( w\right) &=&\left[ w\left(
m_{R}\otimes R\right) \Phi \left( \xi \right) \right] \ast \left[
m_{K}\left( w\otimes \varepsilon _{R}\right) \right] , \\
\partial _{R}^{{2}}\left( {w}\right) &=&  [ m_K\left( \varepsilon _{R}\otimes
w\right)]  \ast \left[ w\left( R\otimes m_{R}\right) \right] \ast
\left[ w^{-1}\left( m_{R}\otimes R\right) \Phi \left( \xi \right)
\right] \ast   [ m_K \left( w^{-1}\otimes \varepsilon _{R}\right) ]
,
\end{eqnarray*}%
where (as in Lemma \ref{lm: 2.12}) $m^{w}:= m_R^w = w\ast m_{R}\ast w^{-1}.$
\end{enumerate}
\end{proposition}

\begin{proof}
i) Since
\begin{equation*}
\left( \partial _{R}^{{0}}\right) _{\pm }\left( w\right) =\Omega
_{H,R}^{1}\left( \partial _{{A}}^{{0}}\right) _{\pm }\mho _{H,R}^{{0}}\left(
w\right) =\Omega _{H,R}^{1}\left( \partial _{{A}}^{{0}}\right) _{\pm }\left(
w\right) =\Omega _{H,R}^{1}\left( w\varepsilon _{A}\right) =w\varepsilon
_{R},
\end{equation*}%
we obtain $\partial _{R}^{{0}}\left( {w}\right) :=\left( \partial
_{R}^{0}\right) _{+}\left( {w}\right) \ast \left( \partial _{R}^{0}\right)
_{-}\left( {w}^{-1}\right) =w\varepsilon _{R}\ast w^{-1}\varepsilon
_{R}=\varepsilon _{R}.$

ii) We compute
\begin{eqnarray*}
\left( \partial _{R}^{{1}}\right) _{+}\left( w\right) &=&\Omega
_{H,R}^{2}\left( \partial _{{A}}^{{1}}\right) _{+}\mho _{H,R}^{{1}}\left(
w\right) =\Omega _{H,R}^{2}\left( \partial _{{A}}^{{1}}\right) _{+}\left(
w\otimes \varepsilon _{H}\right) \\
&=&\Omega _{H,R}^{2}\left[ m_{K}\left( w\otimes \varepsilon _{H}\otimes
w\otimes \varepsilon _{H}\right) \right] =m_{K}\left( {w}\otimes {w}\right) ,%
\text{ and } \\
\left( \partial _{R}^{{1}}\right) _{-}\left( w\right) &=&\Omega _{H,R}^{{2}%
}\left( \partial _{{A}}^{{1}}\right) _{-}\mho _{H,R}^{{1}}\left( w\right)
=\Omega _{H,R}^{{2}}\left( \partial _{{A}}^{{1}}\right) _{-}\left( w\otimes
\varepsilon _{H}\right) \\
&=&\Omega _{H,R}^{{2}}\left( \left( w\otimes \varepsilon _{H}\right)
m_{A}\right) =\left( w\otimes \varepsilon _{H}\right) m_{A}\left( R\otimes
u_{H}\otimes R\otimes u_{H}\right) \overset{(\ref{form: varepsilonHmA})}{=}%
wm_{R}
\end{eqnarray*}%
so that $\partial _{R}^{{1}}\left( {w}\right) =\left( \partial _{R}^{{1}%
}\right) _{+}\left( {w}\right) \ast \left( \partial _{R}^{{1}}\right)
_{-}\left( {w}^{-1}\right) =m_{K}\left( {w}\otimes {w}\right) \ast
w^{-1}m_{R}.$

iii) We first show that
\begin{equation}
m_{K}\circ \left( \varepsilon _{A}\otimes \mho _{H,R}^{{2}}\left( w\right)
\right) =\mho _{H,R}^{{3}}\left[ m_{K}\circ \left( \varepsilon _{R}\otimes
w\right) \right] ,  \label{form: om3a}
\end{equation}%
\begin{equation}
\mho _{H,R}^{{2}}\left( w\right) \circ \left( {A}\otimes m_{A}\right) =\mho
_{H,R}^{{3}}\left[ w\circ \left( R\otimes m_{R}\right) \right] .
\label{form: om3b}
\end{equation}

By the bijectivity of $\Omega _{H,R}^{{3}}$, showing that (\ref{form: om3a})
holds is equivalent to showing that $\Omega _{H,R}^{{3}} (\mathrm{lhs}(\ref%
{form: om3a}))= \Omega _{H,R}^{{3}} (\mathrm{rhs}(\ref{form: om3a}))$. Then
it is only necessary to check (\ref{form: om3a}) on elements of the form $%
\left( {r}\#1_{H}\right) \otimes \left( {s}\#1_{H}\right) \otimes \left( {t}%
\#1_{H}\right) $ with $r,s,t\in R$:
\begin{eqnarray*}
&&m_{K}\left( \varepsilon _{A}\otimes \mho _{H,R}^{{2}}\left( w\right)
\right) \left[ \left( {r}\#1_{H}\right) \otimes \left( {s}\#1_{H}\right)
\otimes \left( {t}\#1_{H}\right) \right] \\
&=&\varepsilon _{A}\left( {r}\#1_{H}\right) \mho _{H,R}^{{2}}\left( w\right) %
\left[ \left( {s}\#1_{H}\right) \otimes \left( {t}\#1_{H}\right) \right]
=\varepsilon _{R}\left( {r}\right) w\left( {s}\otimes {t}\right)
=m_{K}\left( \varepsilon _{R}\otimes w\right) \left( {r}\otimes {s}\otimes {t%
}\right) \\
&=&\mho _{H,R}^{{3}}\left[ m_{K}\left( \varepsilon _{R}\otimes w\right) %
\right] \left[ \left( {r}\#1_{H}\right) \otimes \left( {s}\#1_{H}\right)
\otimes \left( {t}\#1_{H}\right) \right] .
\end{eqnarray*}%
Similarly, to prove (\ref{form: om3b}), we compute:%
\begin{eqnarray*}
&&\mho _{H,R}^{{2}}\left( w\right) \left( {A}\otimes m_{A}\right) \left[
\left( {r}\#1_{H}\right) \otimes \left( {s}\#1_{H}\right) \otimes \left( {t}%
\#1_{H}\right) \right] \\
&=&\mho _{H,R}^{{2}}\left( w\right) \left[ \left( {r}\#1_{H}\right) \otimes
m_{A}\left[ \left( {s}\#1_{H}\right) \otimes \left( {t}\#1_{H}\right) \right]
\right] \\
&=&w\left[ {r}\otimes \left( R\otimes \varepsilon _{H}\right) m_{A}\left[
\left( {s}\#1_{H}\right) \otimes \left( {t}\#1_{H}\right) \right] \right]
\overset{(\ref{form: varepsilonHmA})}{=}w\left[ {r}\otimes m_{R}\left( {s}%
\otimes {t}\right) \right] \\
&=&w\left( R\otimes m_{R}\right) \left( {r}\otimes {s}\otimes {t}\right)
=\mho _{H,R}^{{3}}\left[ w\left( R\otimes m_{R}\right) \right] \left[ \left(
{r}\#1_{H}\right) \otimes \left( {s}\#1_{H}\right) \otimes \left( {t}%
\#1_{H}\right) \right].
\end{eqnarray*}%
Since $\mho _{H,R}^{{3}}$ is a convolution preserving isomorphism, using the
definitions of $\left( \partial _{R}^{{2}}\right) _{+}$ and $\left( \partial
_{{A}}^{{2}}\right) _{+},$
\begin{eqnarray*}
\left[ \mho _{H,R}^{{3}}\circ \left( \partial _{R}^{{2}}\right) _{+}\right]
\left( w\right) &\overset{\text{def}}{=}&\left[ \left( \partial _{{A}}^{{2}%
}\right) _{+}\circ \mho _{H,R}^{{2}}\right] \left( w\right) \\
&\overset{\text{def }}{=}&\left[ m_{K}\left( \varepsilon _{A}\otimes \mho
_{H,R}^{{2}}\left( w\right) \right) \ast \mho _{H,R}^{{2}}\left( w\right)
\left( {A}\otimes m_{A}\right) \right] \\
&\overset{(\ref{form: om3a}),(\ref{form: om3b})}{=}&\mho _{H,R}^{{3}}\left[
m_{K}\left( \varepsilon _{R}\otimes w\right) \right] \ast \mho _{H,R}^{{3}}%
\left[ w\left( R\otimes m_{R}\right) \right] \\
&=&\mho _{H,R}^{{3}}\left[ m_{K}\left( \varepsilon _{R}\otimes w\right) \ast
w\left( R\otimes m_{R}\right) \right]
\end{eqnarray*}%
so that $\left( \partial _{R}^{{2}}\right) _{+} \left( w\right) =\left[
m_{K}\left( \varepsilon _{R}\otimes w\right) \right] \ast \left[ w\left(
R\otimes m_{R}\right) \right] $ as claimed. It is straightforward to verify
the second equality for $\left( \partial _{R}^{{2}}\right) _{+} \left(
w\right)$, or one can use (\ref{form: aurea}).

Similarly, to find the formula for $\left( \partial _{R}^{{2}}\right) _{-}$,
we first prove
\begin{equation}
\left[ \mho _{H,R}^{{2}}\left( w\right) \right] \circ \left( m_{A}\otimes {A}%
\right) =\mho _{H,R}^{{3}}\left[ w\circ \left( m_{R}\otimes R\right) \circ
\Phi \left( \xi \right) \right] ,  \label{form: om3c}
\end{equation}%
\begin{equation}
m_{K}\circ \left( \mho _{H,R}^{{2}}\left( w\right) \otimes \varepsilon
_{A}\right) =\mho _{H,R}^{{3}}\left[ m_{K}\circ \left( w\otimes \varepsilon
_{R}\right) \right] .  \label{form: om3d}
\end{equation}

For $r,s,t\in R$
\begin{eqnarray*}
&&\mho _{H,R}^{{2}}\left( w\right) \left( m_{A}\otimes {A}\right) \left[
\left( {r}\#1_{H}\right) \otimes \left( {s}\#1_{H}\right) \otimes \left( {t}%
\#1_{H}\right) \right] \\
&=&\mho _{H,R}^{{2}}\left( w\right) \left[ m_{A}\left[ \left( {r}%
\#1_{H}\right) \otimes \left( {s}\#1_{H}\right) \right] \otimes \left( {t}%
\#1_{H}\right) \right] \\
&=&\mho _{H,R}^{{2}}\left( w\right) \left[ m_{R}\left( {r}^{(1)}\otimes {r}%
_{\langle -1\rangle }^{(2)}  s ^{(1)}\right) \#\xi \left( {r}%
_{\langle 0\rangle }^{(2)}\otimes {s}^{(2)}\right) \otimes \left( {t}%
\#1_{H}\right) \right] \\
&=&w\left[ m_{R}\left( {r}^{(1)}\otimes {r}_{\langle -1\rangle
}^{(2)}   s^{(1)}\right) \otimes \xi \left( \left( {r}\right)
_{\langle
0\rangle }^{(2)}\otimes    s ^{(2)}\right) {t}\right] \\
&=&w\left( m_{R}\otimes R\right) \left( R\otimes R\otimes \mu _{R}\right)
\left( R\otimes R\otimes \xi \otimes R\right) \left( \Delta _{R\otimes
R}\otimes R\right) \left[ {r}\otimes {s}\otimes {t}\right] \\
&=&w\left( m_{R}\otimes R\right) \Phi \left( \xi \right) \left[ {r}\otimes {s%
}\otimes {t}\right] \\
&=&\mho _{H,R}^{{3}}\left[ w\left( m_{R}\otimes R\right) \Phi \left( \xi
\right) \right] \left[ \left( {r}\#1_{H}\right) \otimes \left( {s}%
\#1_{H}\right) \otimes \left( {t}\#1_{H}\right) \right] ,
\end{eqnarray*}%
and%
\begin{eqnarray*}
&&m_{K}\left( \mho _{H,R}^{{2}}\left( w\right) \otimes \varepsilon
_{A}\right) \left[ \left( {r}\#1_{H}\right) \otimes \left( {s}\#1_{H}\right)
\otimes \left( {t}\#1_{H}\right) \right] \\
&=&\mho _{H,R}^{{2}}\left( w\right) \left[ \left( {r}\#1_{H}\right) \otimes
\left( {s}\#1_{H}\right) \right] \varepsilon _{A}\left( {t}\#1_{H}\right) \\
&=&w\left( {r}\otimes {s}\right) \varepsilon _{R}\left( {t}\right) \\
&=&m_{K}\left( w\otimes \varepsilon _{R}\right) \left( {r}\otimes {s}\otimes
{t}\right) \\
&=&\mho _{H,R}^{{3}}\left[ m_{K}\left( w\otimes \varepsilon _{R}\right) %
\right] \left[ \left( {r}\#1_{H}\right) \otimes \left( {s}\#1_{H}\right)
\otimes \left( {t}\#1_{H}\right) \right]
\end{eqnarray*}%
so that (\ref{form: om3c}) and (\ref{form: om3d}) hold. Since $\mho _{H,R}^{{%
3}}$ is a convolution preserving isomorphism, using the definitions of $%
\left( \partial _{R}^{{2}}\right) _{+}$ and $\left( \partial _{{A}}^{{2}%
}\right) _{+},$ we get%
\begin{eqnarray*}
&&\left[ \mho _{H,R}^{{3}}\circ \left( \partial _{R}^{{2}}\right) _{-}\right]
\left( w\right) \overset{\text{def}}{=}\left[ \left( \partial _{{A}}^{{2}%
}\right) _{-}\circ \mho _{H,R}^{{2}}\right] \left( w\right) \\
&\overset{\text{def}}{=}&\mho _{H,R}^{{2}}\left( w\right) \left(
m_{A}\otimes {A}\right) \ast m_{K}\left( \mho _{H,R}^{{2}}\left( w\right)
\otimes \varepsilon _{A}\right) \\
&\overset{(\ref{form: om3c}),(\ref{form: om3d})}{=}&\mho _{H,R}^{{3}}\left[
w\left( m_{R}\otimes R\right) \Phi \left( \xi \right) \right] \ast \mho
_{H,R}^{{3}}\left[ m_{K}\left( w\otimes \varepsilon _{R}\right) \right] \\
&=&\mho _{H,R}^{{3}}\left[ w\left( m_{R}\otimes R\right) \Phi \left( \xi
\right) \ast m_{K}\left( w\otimes \varepsilon _{R}\right) \right]
\end{eqnarray*}%
so that the formula for $\left( \partial _{R}^{{2}}\right) _{-}\left(
w\right)$ is proved and the formula for $\partial _{R}^{{2}}\left( {w}\right)
$ follows immediately.
\end{proof}

\begin{corollary}
\label{co:Delta2R} Suppose that $w \in \mathrm{Reg}_H(R \otimes R,K)$ is
such that $\xi = u_H w^{-1} \ast \Psi(w)$. Then for $m^w = w \ast m_R \ast
w^{-1}$ as in the proposition, $\left( \partial _{R}^{{2}}\right) _{-}\left(
w^{-1}\right) =\left( {w} ^{-1}\otimes \varepsilon _{R}\right) \ast {w}%
^{-1}(m^w\otimes R) $ and thus
\begin{equation*}
\left( \partial _{R}^{{2}}\right) \left( w\right) ={w}\left( R\otimes
m^w\right) \ast \left( \varepsilon _{R}\otimes {w}\right) \ast \left( {w}%
^{-1}\otimes \varepsilon _{R}\right) \ast {w}^{-1}(m^w\otimes R).
\end{equation*}
\end{corollary}

\begin{proof}
By Proposition \ref{pro:DeltaR}, since $\xi =u_{H}{w}^{-1}\ast \Psi \left( {w%
}\right) $, we have:
\begin{eqnarray*}
\left( \partial _{R}^{{2}}\right) _{-}\left( {w}^{-1}\right) &=&\left[ {w}%
^{-1}\left( m_{R}\otimes R\right) \Phi \left( \xi \right) \right] \ast
\left( {w}^{-1}\otimes \varepsilon _{R}\right) \\
&=& \left[ {w}^{-1}\left( \left( m^w\right) ^{w^{-1}}\otimes R\right) \Phi
\left( u_H{w}^{-1}\ast \Psi \left( {w}\right) \right) \right] \ast \left( {w}%
^{-1}\otimes \varepsilon _{R}\right) \\
&\overset{\left( \ref{form: aurea}\right) }{=}&\left[ \left( w^{-1}\otimes
\varepsilon _{R}\right) \ast {w}^{-1}(m^w\otimes R)\ast \left( {w}\otimes
\varepsilon _{R}\right) \right] \ast \left( {w}^{-1}\otimes \varepsilon
_{R}\right) \\
&=&\left( w^{-1}\otimes \varepsilon _{R}\right) \ast {w}^{-1}(m^w \otimes R).
\end{eqnarray*}

Since by Proposition \ref{pro:DeltaR}, $\left( \partial _{R}^{{2}}\right)
_{+}\left( {w}\right) = {w}\left( R\otimes m^w\right) \ast \left(
\varepsilon _{R}\otimes {w}\right)$, the result follows.
\end{proof}

\section{A new interpretation for the cocycle $\protect\xi $ for a
pre-bialgebra $R$}\label{sect: new interpretation}

Throughout this section $(R, m,u,\Delta, \varepsilon)$ will denote a
connected pre-bialgebra in $_{H}^{H}\mathcal{YD}$.  Define
\begin{equation*}
\Xi :=\{\xi |\xi \text{ is a cocycle for }R\}=\{\xi |\xi \text{ satisfies }(%
\ref{eq:YD3'})-(\ref{eq:YD10'})\text{ in Section \ref{sect: prebialgs}}\},
\end{equation*}%
in other words, the set of $\xi$ such that $(R,\xi)$ is a
pre-bialgebra with cocycle so that a bosonization $A:=R \#_\xi H$
can be built. By Remark \ref{rem:ConvSw}, and by (\ref{eq:YD3'}), $\Xi \subset {\rm Reg}_H(R \otimes R, H)$  where $H$ has the left adjoint action.  In this section, we
establish a bijective correspondence between $\Xi $ and a set $V
\subset \mathrm{\ Reg}_H(R\otimes R,K)$.

First for $C$ any connected coaugmented coalgebra in $_{H}^{H}\mathcal{YD}$,
we find a set in $\mathrm{Reg}_H(C,K)$ in bijective correspondence with the
set of normalized dual Sweedler 1-cocycles in $\mathrm{Reg}_H(C,H)$ which
map $1_{C}$ to $1_{H}$. Recall that $\lambda$ denotes the ad-invariant
integral for $H$, recall from Definition \ref{def: Psi} that $\Psi ({v}%
)=\left( H\otimes {v}\right) \circ \rho _{C},$ for ${v}$ in $\mathrm{Hom}%
(C,K)$ and define the following subsets of $\mathrm{Reg}_H(C,K)$ and $%
\mathrm{Reg}_H(C,H)$:
\begin{equation*}
\mathcal{G}=\left\{ v \in \mathrm{Hom}_H(C,K) \mid {v}(1_{C})=1_{K}\text{
and }\lambda \circ \left[ \Psi ({v})\right] =\lambda \circ \left( H\otimes {v%
}\right) \circ \rho _{C}=\varepsilon _{C}\right\} ,\text{ and }
\end{equation*}%
\begin{eqnarray*}
\mathcal{S} &=&\left\{ \xi\in \mathrm{Hom}_H(C,H) \mid \xi \text{ is a
normalized dual Sweedler $1$-cocycle, and }\xi (1_{C})=1_{H}\right\} \\
&=&\left\{ \xi\in \mathrm{Hom}_H(C,H)\mid \Delta _{H}\xi =(m_{H}\otimes \xi
)(\xi \otimes \rho _{C})\Delta _{C},\varepsilon _{H}\xi =\varepsilon _{C}%
\text{ and }\xi u_{C}:=u_{H}\right\} .
\end{eqnarray*}%
By Remarks \ref{rem:ConvSw} and \ref{rem: Takeuchi} the sets $\mathcal{G}$
and $\mathcal{S}$ consist of convolution invertible maps.

\begin{theorem}
\label{teo: bijection} For $C$, $\mathcal{G}$,$\mathcal{S}$ as above, the
maps $F :\mathcal{G}\rightarrow \mathcal{S}$ and $G :\mathcal{S}\rightarrow
\mathcal{G}$ defined by%
\begin{equation*}
F \left( {v}\right) =u_{H}{v}^{-1}\ast \Psi \left( {v}\right) \qquad \text{%
and}\qquad G \left( \xi \right) =(\lambda \xi )^{-1}
\end{equation*}%
are inverse bijections.
\end{theorem}

\begin{proof}
Let ${v}\in \mathcal{G}$, and $\xi :=u_{H}{v}^{-1}\ast \Psi \left( {v}%
\right) $. We check that $\xi \in \mathcal{S}.$ First we see that $\xi
u_{C}=u_{H}$ since
\begin{equation*}
\xi \left( 1_{C}\right) =\left[ u_{H}{v}^{-1}\ast \Psi \left( {v}\right) %
\right] \left( 1_{C}\right) =u_{H}{v}^{-1}\left( 1_{C}\right) \Psi \left( {v}%
\right) \left( 1_{C}\right) =1_{H}.
\end{equation*}%
Also
\begin{equation*}
\varepsilon _{H}\xi =\varepsilon _{H}\left[ u_{H}{v}^{-1}\ast \Psi \left( {v}%
\right) \right] =\varepsilon _{H}u_{H}{v}^{-1}\ast \varepsilon _{H}\Psi
\left( {v}\right) ={v}^{-1}\ast {v}=\varepsilon _{C}.
\end{equation*}%
Moreover, for $x\in C$%
\begin{eqnarray*}
&&(m_{H}\otimes \xi )(\xi \otimes \rho _{C})\Delta _{C}(x)=(m_{H}\otimes \xi
)(\xi (x^{(1)})\otimes x_{\langle -1\rangle }^{(2)}\otimes x_{\langle
0\rangle }^{(2)}) \\
&=&\xi (x^{(1)})x_{\langle -1\rangle }^{(2)}\otimes \xi (x_{\langle 0\rangle
}^{(2)})=[[u_{H}{v}^{-1}\ast \Psi ({v})](x^{(1)})]x_{\langle -1\rangle
}^{(2)}\otimes \lbrack u_{H}{v}^{-1}\ast \Psi ({v})](x_{\langle 0\rangle
}^{(2)}) \\
&=&u_{H}{v}^{-1}(x^{(1)})\Psi ({v})(x^{(2)})x_{\langle -1\rangle
}^{(3)}\otimes u_{H}{v}^{-1}[(x_{\langle 0\rangle }^{(3)})^{(1)}]\Psi ({v}%
)[(x_{\langle 0\rangle }^{(3)})^{(2)}] \\
&\overset{(\ref{form: Delta colin})}{=}&u_{H}{v}^{-1}(x^{(1)})\Psi ({v}%
)(x^{(2)})x_{\langle -1\rangle }^{(3)}x_{\langle -1\rangle }^{(4)}\otimes
u_{H}{v}^{-1}(x_{\langle 0\rangle }^{(3)})\Psi ({v})(x_{\langle 0\rangle
}^{(4)}) \\
&=&{v}^{-1}(x^{(1)})x_{\langle -1\rangle }^{(2)}{v}(x_{\langle 0\rangle
}^{(2)})x_{\langle -1\rangle }^{(3)}x_{\langle -1\rangle }^{(4)}\otimes {v}%
^{-1}(x_{\langle 0\rangle }^{(3)})\Psi ({v})(x_{\langle 0\rangle }^{(4)}) \\
&=&{v}^{-1}(x^{(1)})x_{\langle -1\rangle }^{(2)}x_{\langle -1\rangle
}^{(3)}x_{\langle -1\rangle }^{(4)}\otimes {v}(x_{\langle 0\rangle }^{(2)}){v%
}^{-1}(x_{\langle 0\rangle }^{(3)})\Psi ({v})(x_{\langle 0\rangle }^{(4)}) \\
&\overset{(\ref{form: Delta colin})}{=}&{v}^{-1}(x^{(1)})x_{\langle
-1\rangle }^{(2)}\otimes {v}((x_{\langle 0\rangle }^{(2)})^{(1)}){v}%
^{-1}((x_{\langle 0\rangle }^{(2)})^{(2)})\Psi ({v})((x_{\langle 0\rangle
}^{(2)})^{(3)}) \\
&=&{v}^{-1}(x^{(1)})x_{\langle -1\rangle }^{(2)}\otimes ({v}\ast {v}%
^{-1})((x_{\langle 0\rangle }^{(2)})^{(1)})\Psi ({v})((x_{\langle 0\rangle
}^{(2)})^{(2)}) \\
&=&{v}^{-1}(x^{(1)})x_{\langle -1\rangle }^{(2)}\otimes \Psi ({v}%
)(x_{\langle 0\rangle }^{(2)})={v}^{-1}(x^{(1)})x_{\langle -2\rangle
}^{(2)}\otimes x_{\langle -1\rangle }^{(2)}{v}(x_{\langle 0\rangle }^{(2)})
\\
&=&\Delta _{H}[{v}^{-1}(x^{(1)})x_{\langle -1\rangle }^{(2)}{v}(x_{\langle
0\rangle }^{(2)})]=\Delta _{H}\xi (x).
\end{eqnarray*}%
Thus $F $ maps $\mathcal{G}$ to $\mathcal{S}$. Now let $\xi \in \mathcal{S}.$
Then, $\lambda \xi u_{C}=\mathrm{Id}_{K}$ so that, by Remark \ref{rem:
Takeuchi}, $\lambda \xi \ $is convolution invertible and $\left( \lambda \xi
\right) ^{-1}u_{C}=\mathrm{Id}_{K}.$ Thus it makes sense to consider ${v}%
:=\left( \lambda \xi \right) ^{-1}$ and to see that ${v}\in \mathcal{G}$, it
remains to show that $\lambda \circ \Psi ({v})=\varepsilon _{C}$. Since $%
\lambda $ is a left integral, $u_{H}\lambda =(H\otimes \lambda )\Delta _{H}$%
. Thus,%
\begin{eqnarray*}
u_{H}{v}^{-1} &=&u_{H}\lambda \xi =\left( H\otimes \lambda \right) \Delta
_{H}\xi \overset{(\ref{eq: Sweedler 1-cocycle})}{=} (m_{H}\otimes \lambda
\xi )(\xi \otimes \rho _{C})\Delta _{C} \\
&=&(m_{H}\otimes {v}^{-1})(\xi \otimes \rho _{C})\Delta _{C}=\xi \ast \Psi
\left( {v}^{-1}\right)
\end{eqnarray*}%
so that $\xi =u_{H}{v}^{-1}\ast \Psi \left( {v}\right) .$ Hence we have $%
\lambda \Psi \left( {v}\right) =\lambda \left[ u_{H}{v}\ast \xi \right] ={v}%
\ast \lambda \xi =\varepsilon _{C}.$ Thus ${v}\in \mathcal{G}$ and $\xi =F ({%
v})$. This proves that $G$ maps $\mathcal{S}$ to $\mathcal{G}$ and $F \circ
G = Id_{\mathcal{S}} $.

Next we show that $G \circ F =Id_{\mathcal{G}}$. Let ${v}\in \mathcal{G}$.
Let $\xi :=u_{H}{v}^{-1}\ast \Psi \left( {v}\right) =F ({v})\in \mathcal{S}$%
. We check that ${v}=\left( \lambda \xi \right) ^{-1}$. We have%
\begin{equation*}
\lambda \xi =\lambda \circ \left[ u_{H}{v}^{-1}\ast \Psi \left( {v}\right) %
\right] ={v}^{-1}\ast \lambda \Psi \left( {v}\right) ={v}^{-1}\ast
\varepsilon _{C}={v}^{-1}.
\end{equation*}%
Thus ${v}=(\lambda \xi )^{-1}=\left( G \circ F \right) ({v})$.

Finally suppose that ${v}\in \mathcal{G}$ is left $H$-linear, i.e., ${v}%
(h\cdot c)=\varepsilon_H  (h){v}(c)$. Then ${v}^{-1}$ and $\Psi
({v})$ are also left $H$-linear by Remark \ref{rem: Takeuchi}, and
Lemma \ref{lem: PsiYD} so that their convolution product $F(v)$ is
left $H$-linear, as required.
Conversely, if $\xi$ is left $H$-linear, so is $\lambda \xi$ and thus so is $%
G(\xi)$.
\end{proof}

We note that if $C$ is coaugmented but not necessarily connected, then $F $
is still a map from $ \mathrm{\ Reg}_H(C,K)$ to the set $\mathcal{S}$.
\vspace{1.5mm}

Recall that if $(R,\xi )$ is a pre-bialgebra with cocycle,
%so that $%
%(A:=R\#_{\xi }H,H,\pi ,\sigma )$ is a splitting datum, then
the cohomology
of $(R,\xi )$ is defined using the cohomology of $R \#_\xi H$. By Corollary \ref%
{co:Delta2R}, if $v\in \mathrm{Reg}_{H}(R\otimes R,K)$ is such that
$\xi =u_{H}v^{-1}\ast \Psi (v)=F(v)$ then $\partial
_{R}^{2}(v)={v}\left( R\otimes m^v \right) \ast \left( \varepsilon
_{R}\otimes {v}\right) \ast \left( {v}^{-1}\otimes \varepsilon
_{R}\right) \ast {v}^{-1} (m^v\otimes R).$
However, this expression makes
sense for $v\in \mathrm{Reg}_{H}(R\otimes R,K)$ even if
$u_{H}v^{-1}\ast \Psi (v) \notin \Xi $.

\begin{definition}
\label{def: alpha} Let $(R,m,u,\Delta ,\varepsilon )$ be a connected
pre-bialgebra in ${_{H}^{H}\mathcal{YD}}$. For $w\in \mathrm{Reg}%
_{H}(R\otimes R,K)$ define $\alpha (w)=\alpha _{+}(w)\ast \alpha
_{-}(w^{-1}) $ where:
\begin{equation*}
\alpha _{+}(w):=w(R\otimes m^w) \ast (  \varepsilon \otimes
w),\mbox{ and }\alpha _{-}\left( w^{-1}\right) =\left(
{w^{-1}}\otimes   \varepsilon \right) \ast {w^{-1}}(m^w \otimes R).
\end{equation*}
\end{definition}

Theorem \ref{teo: bijection} establishes a bijection between $\mathcal{S}
\supset \Xi$ and a set $\mathcal{G} \subset \mathrm{Reg}_H(R \otimes R, K)$.
The rest of this section will be devoted to proving the following theorem:

\begin{theorem}
\label{teo: hope} Let $(R,m,u,\Delta ,\varepsilon )$ be a connected
pre-bialgebra in ${_{H}^{H}\mathcal{YD}}$, and let $\mathcal{S}$, $\mathcal{G%
}$,$F,G$ be as in Theorem \ref{teo: bijection}. Then $G(\Xi )=V$ where $%
V\subset \mathcal{G}$ is the set of $v$ in $\mathcal{G}$ satisfying the
following:
\begin{eqnarray}
&&m^v \text{is $H$-colinear, i.e. } \rho _{R}m^v =\left( H\otimes m^v\right)
\rho _{R\otimes R},  \label{form: mv colinear} \\
&&m^v\left( R\otimes m^v\right) \ast \alpha \left( v\right) =\alpha \left(
v\right) \ast m^v(m^v\otimes R),  \label{form: mv associative} \\
&&\Psi \left( \alpha (v)\right) =u_{H}\alpha \left( v\right) ,
\label{form: quasi
cocycle} \\
&& v \text{ is unital, i.e. for all } r \in R, {v}\left( r\otimes
1_{R}\right) ={v}\left( 1_{R}\otimes r\right) =\varepsilon \left( r\right).
\label{form: v strongly unitary}
\end{eqnarray}
\end{theorem}

Note that (\ref{form: mv associative}) says that $m^v$ is associative up to
inner action by the invertible element $\alpha (v)$.

\begin{remark}
\label{rem: alpha YD} We note that condition (\ref{form: quasi cocycle}%
) means that $\alpha \left( v\right) :R^{\otimes 3}\rightarrow K$ is left $H$%
-colinear. Since $\alpha \left( v\right) $ is  also  left
$H$-linear, $\alpha \left( v\right) $ is in
${_{H}^{H}\mathcal{YD}}\mathbf{.}$ If in addition $v$ is unital,
then $\alpha \left( v\right) $ is unital so that
one gets that $\alpha \left( v\right) $ is a retraction of the unit $%
u_{R^{\otimes 3}}:K\rightarrow R^{\otimes 3}$ of the coaugmented coalgebra $%
R^{\otimes 3}$ in the category ${_{H}^{H}\mathcal{YD}}\mathbf{.}$
\end{remark}

The proof of Theorem \ref{teo: hope} consists of a series of propositions
equating the conditions (\ref{eq:YD6'}) to (\ref{eq:YD10'}) from Definition %
\ref{de: xi} to the conditions (\ref{form: mv colinear}) to (\ref{form: v
strongly unitary}) listed above.

By Lemma \ref{lm: YD6' YD}, condition (\ref{eq:YD6'}) is equivalent to (\ref%
{form: equivalent to YD6'}), namely that for all $z\in C:=R \otimes R$,
\begin{equation*}
\rho _{R}(m(z))=m(z)_{\langle -1\rangle }\otimes m(z)_{\langle 0\rangle
}=\xi (z^{(1)})z_{\langle -1\rangle }^{(2)}\xi ^{-1}(z^{(3)})\otimes
m(z_{\langle 0\rangle }^{(2)}).
\end{equation*}
For $v\in \mathrm{Reg}_{H}(R\otimes R,K)$, $m^v$ is $H$-linear, being the
product of left $H$-linear maps. The next proposition shows that (\ref%
{eq:YD6'}) (equivalently (\ref{form: equivalent to YD6'})) holds for $\xi \in \mathcal{%
S}$ if and only if $m^v$ is left $H $-colinear where $v=G(\xi )$.

\begin{proposition}
\label{pro: YD6'} Let $v \in \mathcal{G}$ and let $\xi := F(v) = u_{H}{v}%
^{-1}\ast \Psi ({v} ) $. Then $\xi $ satisfies (\ref{eq:YD6'}) if and only
if $m^v$ is left $H$-colinear, i.e.,
\begin{equation*}
\rho _{R}m^v=(H\otimes m^v)\rho _{R\otimes R}.
\end{equation*}
\end{proposition}

\begin{proof}
Since $\xi =u_{H}{v}^{-1}\ast \Psi ( {v}) $ then $\xi ^{-1}=\Psi ( {v}^{-1})
\ast u_{H}{v},$ where $\Psi ( \gamma ) =( H\otimes \gamma ) \rho _{C}$ and $%
C = R \otimes R$. Hence $\xi ^{-1}( z) =z_{\langle -1\rangle }^{( 1) }{v}%
^{-1}( z_{\langle 0\rangle }^{( 1) }) {v}( z^{( 2) }) $. Then the right side
of (\ref{form: equivalent to YD6'}) applied to $z \in C$ is:
\begin{eqnarray*}
&&\xi ( z^{( 1) }) z_{\langle -1\rangle }^{( 2) }\xi ^{-1}( z^{( 3) })
\otimes m( z_{\langle 0\rangle }^{( 2) }) \\
&=&[ {v}^{-1}( z^{( 1) }) z_{\langle -1\rangle }^{( 2) }{v}( z_{\langle
0\rangle }^{( 2) }) ] z_{\langle -1\rangle }^{( 3) }[ z_{\langle -1\rangle
}^{( 4) }{v}^{-1}( z_{\langle 0\rangle }^{( 4) }) {v}( z^{( 5) }) ] \otimes
m( z_{\langle 0\rangle }^{( 3) }) \\
&=&{v}^{-1}( z^{( 1) }) z_{\langle -1\rangle }^{( 2) }z_{\langle -1\rangle
}^{( 3) }z_{\langle -1\rangle }^{( 4) }{v}( z^{( 5) }) \otimes {v}(
z_{\langle 0\rangle }^{( 2) }) m( z_{\langle 0\rangle }^{( 3) }) {v}^{-1}(
z_{\langle 0\rangle }^{( 4) }) \\
&\overset{(\ref{form: Delta colin})}{=}& {v}^{-1}( z^{( 1) }) z_{\langle
-1\rangle }^{( 2) }{v}( z^{( 3) }) \otimes {v}( ( z_{\langle 0\rangle }^{(
2) }) ^{( 1) }) m( ( z_{\langle 0\rangle }^{( 2) }) ^{( 2) }) {v}^{-1}( (
z_{\langle 0\rangle }^{( 2) }) ^{( 3) }) \\
&=& {v}^{-1}( z^{( 1) }) z_{\langle -1\rangle }^{( 2) }{v}( z^{( 3) })
\otimes ( {v}\ast m\ast {v}^{-1}) ( z_{\langle 0\rangle }^{( 2) }) \\
&=& z_{\langle -1\rangle }^{( 2) }\otimes m^v( {v}^{-1}( z^{( 1) })
z_{\langle 0\rangle }^{( 2) }{v}( z^{( 3) }) ) \\
&=&( H\otimes m^v) \rho_C ( {v}^{-1}( z^{( 1) }) z^{( 2) }{v
}( z^{( 3) }) ) \\
&=&( H\otimes m^v) \rho _{C}( u_{C}{v}^{-1}\ast \mathrm{Id}_{C}\ast u_{C}{v}%
) ( z) .
\end{eqnarray*}%
so that
\begin{equation}
\xi ( z^{( 1) }) z_{\langle -1\rangle }^{( 2) }\xi ^{-1}( z^{( 3) }) \otimes
m( z_{\langle 0\rangle }^{( 2) }) =( H\otimes m^v )\rho _{C}( u_{C}{v}%
^{-1}\ast \mathrm{Id}_{C}\ast u_{C}{v}) ( z) .  \label{form: YD6lem}
\end{equation}%
Note that $( u_{C}{v}^{-1}\ast \mathrm{Id}_{C}\ast u_{C}{v}) $ and $( u_{C}{v%
}\ast \mathrm{Id}_{C}\ast u_{C}{v}^{-1}) $ are composition inverses. Thus if
(\ref{form: equivalent to YD6'}) holds, by (\ref{form: YD6lem}) we have
\begin{equation*}
\rho _{R}  m =( H\otimes m^v) \rho _{C}( u_{C}{v}^{-1}\ast \mathrm{Id}%
_{C}\ast u_{C}{v})
\end{equation*}%
whence%
\begin{equation*}
\rho _{R}m^v = \rho _{R}  m ( u_{C}{v}\ast \mathrm{Id}_{C}\ast u_{C}{v}%
^{-1}) =( H\otimes m^v ) \rho _{C}
\end{equation*}%
as required. Conversely, if (\ref{form: mv colinear}) holds, then
\begin{equation*}
\xi ( z^{( 1) }) z_{\langle -1\rangle }^{( 2) }\xi ^{-1}( z^{( 3) })
\otimes m( z_{\langle 0\rangle }^{( 2) }) \overset{(\ref{form:
YD6lem})}{=}\rho _{R}m^v( u_{C}{v}^{-1}\ast \mathrm{Id}_{C}\ast
u_{C}{v}) ( z) =\rho _{R}  m ( z) .
\end{equation*}
\end{proof}

Now we can see the relationship between the associativity condition (\ref%
{eq:YD7'}) on $m_{R}$ with $\xi =F(v)$ and the associativity of the
multiplication $m^v$ on $R$.

\begin{proposition}
\label{pro: YD7'} Let $v \in \mathcal{G}$ and assume that $v$ satisfies (\ref%
{form: mv colinear}) so that $m^v$ is left $H$-colinear. Let $\xi := F(v) = {%
v} ^{-1}\ast \Psi \left( {v}\right) .$ Then $\left( \ref{eq:YD7'}\right) $
holds for $\xi $ if and only if (\ref{form: mv associative}) holds for $v$,
i.e.,
\begin{equation*}
m^v\left( R\otimes m^v \right) \ast \alpha(v) = \alpha(v) \ast m^v (m^v
\otimes R),
\end{equation*}
where $\alpha(v) $ was defined in Definition \ref{def: alpha}.
\end{proposition}

\begin{proof}
Since $m^v$ is left $H$-linear and colinear, we may apply $\left( \ref{form:
aurea}\right) $ and $\left( \ref{form: alfabeta}\right) $. We have
\begin{eqnarray*}
  m (R\otimes m ) &=&  m (R\otimes \left( {v}^{-1}\ast m^v \ast {v}%
\right) ) \\
&\overset{\left( \ref{form: aurea}\right) }{=}&\left(  \varepsilon
\otimes {v}^{-1}\right) \ast   m \left( R\otimes m^v \right) \ast
\left(   \varepsilon \otimes {v}\right) \\
&=&\left(   \varepsilon \otimes {v}^{-1}\right) \ast \left(
{v}^{-1}\ast
m^v \ast {v}\right) \left( R\otimes m^v \right) \ast \left(   \varepsilon \otimes {v}\right) \\
&\overset{\left( \ref{form: alfabeta}\right) }{=}&\left( \varepsilon
\otimes {v}^{-1}\right) \ast {v}^{-1}\left( R\otimes m^v \right)
\ast m^v \left( R\otimes m^v \right) \ast {v}\left( R\otimes m^v
\right) \ast \left(   \varepsilon \otimes {v}\right)
\end{eqnarray*}%
On the other hand,
\begin{eqnarray*}
  m ( m \otimes R)\Phi \left( \xi \right) &=&
m (  m \otimes R)\Phi
\left( {v}^{-1}\ast \Psi \left( {v}\right) \right) \\
&\overset{\left( \ref{form: aurea}\right) }{=}&\left(
{v}^{-1}\otimes   \varepsilon \right) \ast   m (m^v \otimes R)\ast
\left( {v}\otimes
  \varepsilon \right) \\
&=&\left( {v}^{-1}\otimes   \varepsilon \right) \ast \left(
{v}^{-1}\ast
m^v \ast {v}\right) (m^v \otimes R)\ast \left( {v}\otimes   \varepsilon \right) \\
&\overset{\left( \ref{form: alfabeta}\right) }{=}&\left(
{v}^{-1}\otimes   \varepsilon \right) \ast {v}^{-1}(m^v \otimes
R)\ast m^{{v} }(m^v \otimes R)\ast {v}(m^v \otimes R)\ast \left(
{v}\otimes   \varepsilon \right).
\end{eqnarray*}
Since for any coalgebra map $\varphi: R \otimes R \otimes R \rightarrow R
\otimes R$, $v \circ \varphi$ and $v^{-1} \circ \varphi$ are convolution
inverses, the statement follows.
\end{proof}

\begin{proposition}
\label{pro: YD8'} Let $v \in \mathcal{G}$ and assume that $v$ satisfies (\ref%
{form: mv colinear}) so that $m^v$ is left $H$-colinear. Set $\xi := F(v) = {%
v} ^{-1}\ast \Psi \left( {v}\right)$. Then $\left( \ref{eq:YD8'}\right) $
holds for $\xi $ if and only if (\ref{form: quasi cocycle}) holds for $v$,
i.e.,
\begin{equation*}
\Psi \left( \alpha\left( v\right) \right) =u_{H}\alpha \left( v\right) .
\end{equation*}
\end{proposition}

\begin{proof}
First consider the left side of $\left( \ref{eq:YD8'}\right)$. For $r\in R$
and $w\in R\otimes R,$ we have%
\begin{eqnarray*}
&& m_{H}\left( \xi \otimes H\right) \left[ R\otimes \left( m\otimes
\xi
\right) \Delta _{R\otimes R}\right] \left( r\otimes w\right) \\
&=& \xi \left[ r\otimes  m \left( w^{\left( 1\right) }\right)
\right] \xi
\left( w^{\left( 2\right) }\right) \\
&=& \xi \left[ r\otimes m\left( w^{\left( 1\right) }\right) \right] %
\left[ u_{H}{v}^{-1}\ast \Psi \left( {v}\right) \right] ( w^{\left(
2\right)} ) \\
&=& \xi \left[ r\otimes m( w^{\left( 1\right)} ) \right] {v}^{-1}(
w^{\left( 2\right)} ) \Psi \left( {v}\right) ( w^{\left( 3\right)} ) \\
&=& \xi \left[ r\otimes \left(m\ast {v}^{-1}\right) \left( w^{\left(
1\right) }\right) \right] \Psi \left( {v}\right) ( w^{\left( 2\right)} ) \\
&=& \xi \left[ r\otimes \left( {v}^{-1}\ast m^v \right) \left( w^{\left(
1\right) }\right) \right] \Psi \left( {v}\right) ( w^{\left( 2 \right ) } )
\\
&=&{\ v}^{-1}(w^{\left( 1\right) })\xi \left[ r\otimes m^v ( w^{\left(
2\right)} ) \right] \Psi \left( {v}\right) (w^{\left( 3\right) }) \\
&=& {v}^{-1}(w^{\left( 1\right) })\xi \left( R\otimes m^v \right) (r\otimes
w^{\left( 2\right) })\Psi \left( {v}\right) (w^{\left( 3\right) }) \\
&=& {v}^{-1}(w^{\left( 1\right) })\left[ u_{H}{v}^{-1}\ast \Psi \left( {v}%
\right) \right] \left( R\otimes m^v \right) (r\otimes w^{\left( 2\right)
})\Psi \left( {v}\right) (w^{\left( 3\right) }) \\
&\overset{\text{(\ref{form: alfabeta})}}{=}&{v}^{-1}(w^{\left( 1\right) })%
\left[ u_{H}{v}^{-1}\left( R\otimes m^v \right) \ast \Psi \left( {v}\right)
\left( R\otimes m^v \right) \right] (r\otimes w^{\left( 2\right) })\Psi
\left( {v}\right) (w^{\left( 3\right) }) \\
&=&( \varepsilon \otimes {v}^{-1})(r^{\left( 1\right) }\otimes
r_{\left\langle -1\right\rangle }^{\left( 2\right) }w^{\left( 1\right) })%
\left[ u_{H}{v}^{-1}(R\otimes m^v )\ast \Psi ({v)}(R\otimes m^v) \right]
(r_{\left\langle 0\right\rangle }^{\left( 2\right) }\otimes w^{\left(
2\right) })\Psi \left( {v}\right) (w^{\left( 3\right) }) \\
&=&\left[ u_{H}\left(  \varepsilon \otimes {v}^{-1}\right) \ast u_{H}{v}%
^{-1}\left( R\otimes m^v \right) \ast \Psi \left( {v}\right) \left( R\otimes
m^v \right) \right] \left( r\otimes w^{\left( 1\right) }\right) \Psi \left( {%
v}\right) \left( w^{\left( 2\right) }\right) \\
&=& \left[ u_{H}\left( \alpha\right) _{+}\left( v\right) ^{-1}\ast
\Psi \left( {v}\right) \left( R\otimes m^v \right) \right] \left(
r^{\left( 1\right) }\otimes r_{\left\langle -1\right\rangle
}^{\left( 2\right) }w^{\left( 1\right) }\right) \Psi \left(
\varepsilon \otimes {v}\right) \left( r_{\left\langle 0\right\rangle
}^{\left( 2\right) }\otimes w^{\left(
2\right) }\right) \\
&\overset{\left( \ref{rem: hope1.8}(i)\right) }{=}& \left[ u_{H}\left(
\alpha\right) _{+}\left( v\right) ^{-1}\ast \Psi \left( {v}\left( R\otimes
m^v \right) \right) \ast \Psi \left(   \varepsilon \otimes {v}\right) %
\right] \left( r\otimes w\right) \\
&\overset{\left( \ref{rem: hope1.8}(ii)\right) }{=}& \left[
u_{H}\left( \alpha\right) _{+}\left( v\right) ^{-1}\ast \Psi \left[
{v}\left( R\otimes m^v \right) \ast \left(   \varepsilon \otimes
{v}\right) \right] \right]
\left( r\otimes w\right) \\
&=& \left[ u_{H}\left( \alpha\right) _{+}\left( v\right) ^{-1}\ast \Psi
\left( \left( \alpha\right) _{+}\left( v\right) \right) \right] \left(
r\otimes w\right) .
\end{eqnarray*}%
so that the left hand side of $\left( \ref{eq:YD8'}\right) $ equals $%
u_{H}\left( \alpha\right) _{+}\left( v\right) ^{-1}\ast \Psi \left( \left(
\alpha\right) _{+}\left( v\right) \right) $.  Now consider the right side of
$\left( \ref{eq:YD8'}\right) .$ For $z\in R\otimes R$ and $t\in R$, we have%
\begin{eqnarray*}
&&m_{H}\left( \xi \otimes H\right) \left( R\otimes c_{H,R}\right)
\left[ \left( m\otimes \xi \right) \Delta _{R\otimes R}\otimes
R\right] \left(
z\otimes t\right) \\
&=&\xi \left[ m\left( z^{\left( 1\right) }\right) \otimes \xi \left(
z^{\left( 2\right) }\right) _{\left( 1\right) }t\right] \xi \left(
z^{\left(
2\right) }\right) _{\left( 2\right) } \\
&=&\xi \left[ m\left( z^{\left( 1\right) }\right) \otimes \left\{
\left[ u_{H}{v}^{-1}\ast \Psi \left( {v}\right) \right] \left(
z^{\left( 2\right)
}\right) \right\} _{\left( 1\right) }t\right] \left\{ \left[ u_{H}{v}%
^{-1}\ast \Psi \left( {v}\right) \right] \left( z^{\left( 2\right) }\right)
\right\} _{\left( 2\right) } \\
&=&\xi \left[ m\left( z^{\left( 1\right) }\right) \otimes
{v}^{-1}\left( z^{\left( 2\right) }\right) z_{\left\langle
-2\right\rangle }^{\left(
3\right) }t\right] z_{\left\langle -1\right\rangle }^{\left( 3\right) }{v}%
\left( z_{\left\langle 0\right\rangle }^{\left( 3\right) }\right) \\
&=&\xi \left[ \left( m_{R}\ast {v}^{-1}\right) \left( z^{\left( 1\right)
}\right) \otimes z_{\left\langle -1\right\rangle }^{\left( 2\right) }t\right]
\Psi \left( {v}\right) \left( z_{\left\langle 0\right\rangle }^{\left(
2\right) }\right) \\
&=&\xi \left[ \left( {v}^{-1}\ast m^v \right) \left( z^{\left( 1\right)
}\right) \otimes z_{\left\langle -1\right\rangle }^{\left( 2\right) }t\right]
\Psi \left( {v}\right) \left( z_{\left\langle 0\right\rangle }^{\left(
2\right) }\right) \\
&=&{v}^{-1}\left( z^{\left( 1\right) }\right) \xi \left[ m^v \left(
z^{\left( 2\right) }\right) \otimes z_{\left\langle -1\right\rangle
}^{\left( 3\right) }t\right] \Psi \left( {v}\right) \left( z_{\left\langle
0\right\rangle }^{\left( 3\right) }\right) \\
&=&{v}^{-1}\left( z^{\left( 1\right) }\right) \xi \left( m^v \otimes
R\right) \left( z^{\left( 2\right) }\otimes z_{\left\langle -1\right\rangle
}^{\left( 3\right) }t\right) \Psi \left( {v}\right) \left( z_{\left\langle
0\right\rangle }^{\left( 3\right) }\right) \\
&=&{v}^{-1}\left( z^{\left( 1\right) }\right) \left[ u_{H}{v}^{-1}\ast \Psi
\left( {v}\right) \right] \left( m^v \otimes R\right) \left( z^{\left(
2\right) }\otimes z_{\left\langle -1\right\rangle }^{\left( 3\right)
}t\right) \Psi \left( {v}\right) \left( z_{\left\langle 0\right\rangle
}^{\left( 3\right) }\right) \\
&\overset{\text{(\ref{form: alfabeta})}}{=}&{v}^{-1}\left( z^{\left(
1\right) }\right) \left[ u_{H}{v}^{-1}\left( m^v \otimes R\right) \ast \Psi
\left( {v}\right) \left( m^v \otimes R\right) \right] \left( z^{\left(
2\right) }\otimes z_{\left\langle -1\right\rangle }^{\left( 3\right)
}t\right) \Psi \left( {v}\right) \left( z_{\left\langle 0\right\rangle
}^{\left( 3\right) }\right) \\
&\overset{\left( \ref{rem: hope1.8}(i)\right) }{=}&{v}^{-1}\left(
z^{\left( 1\right) }\right) \left[ u_{H}{v}^{-1}\left( m^v \otimes
R\right) \ast \Psi \left( {v}\left( m^v \otimes R\right) \right)
\right] \left( z^{\left( 2\right) }\otimes z_{\left\langle
-1\right\rangle }^{\left( 3\right) }t^{\left( 1\right) }\right) \Psi
\left( {v}\otimes   \varepsilon \right) \left( z_{\left\langle
0\right\rangle }^{\left( 3\right) }\otimes t^{\left(
2\right) }\right) \\
&=& {v}^{-1}\left( z^{\left( 1\right) }\right) \left[
u_{H}{v}^{-1}\left( m^v \otimes R\right) \ast \Psi \left( {v}\left(
m^v \otimes R\right) \right) \ast \Psi \left( {v}\otimes
\varepsilon \right) \right] \left( z^{\left(
2\right) }\otimes t\right) \\
&\overset{(\ref{form: alfabeta})}{=}&{v}^{-1}\left( z^{\left( 1\right)
}\right) \left[ u_{H}{v}^{-1}\left( m^v \otimes R\right) \ast \Psi \left(
\left( \alpha\right) _{-}\left( v^{-1}\right) ^{-1}\right) \right] \left(
z^{\left( 2\right) }\otimes t\right) \\
&=&\left( {v}^{-1}\otimes   \varepsilon \right) \left( z^{\left(
1\right) }\otimes z_{\left\langle -1\right\rangle }^{\left( 2\right)
}t^{\left( 1\right) }\right) \left[ u_{H}{v}^{-1}\left( m^v \otimes
R\right) \ast \Psi \left( \left( \alpha\right) _{-}\left(
v^{-1}\right) ^{-1}\right) \right] \left( z_{\left\langle
0\right\rangle }^{\left( 2\right) }\otimes t^{\left(
2\right) }\right) \\
&=&\left[ u_{H}\left( {v}^{-1}\otimes  \varepsilon \right) \ast u_{H}{v}%
^{-1}\left( m^v \otimes R\right) \ast \Psi \left( \left( \alpha\right)
_{-}\left( v^{-1}\right) ^{-1}\right) \right] \left( z\otimes t\right) \\
&=&\left[ u_{H}\left( \alpha\right) _{-}\left( v^{-1}\right) \ast \Psi
\left( \left( \alpha\right) _{-}\left( v^{-1}\right) ^{-1}\right) \right]
\left( z\otimes t\right)
\end{eqnarray*}%
so that the right hand side of $\left( \ref{eq:YD8'}\right) $ is $%
u_{H}\left( \alpha\right) _{-}\left( v^{-1}\right) \ast \Psi \left( \left(
\alpha\right) _{-}\left( v^{-1}\right) ^{-1}\right) $. By Remark $\ref{rem:
hope1.8}(ii)$, $\Psi $ is an algebra map, and the statement follows.
\end{proof}

Finally we show that $\xi$ is unital if and only if $v= G(\xi)$ is unital.

\begin{proposition}
\label{pro: YD10'} Let ${v} \in \mathcal{G}$ and let $\xi : F(v) =u_{H}{v}%
^{-1}\ast \Psi \left( {v} \right) .$ Then $\xi $ satisfies $\left( \ref%
{eq:YD10'}\right) $ if and only if $v$ is unital, i.e., for all $r \in R$, $%
v(1_R \otimes r) = v(r \otimes 1_R) = \varepsilon(r)$.
\end{proposition}

\begin{proof}
Note that $v$ is unital if and only if $v^{-1}$ is. Since for all $r\in R$,
\begin{equation*}
\xi (r\otimes 1_{R})=\left[ u_{H}{v}^{-1}\ast \Psi \left( {v}\right) \right]
(r\otimes 1_{R})={v}^{-1}(r^{(1)}\otimes 1_{R})r_{\langle -1\rangle }^{(2)}{v%
}(r_{\langle 0\rangle }^{(2)}\otimes 1_{R}),
\end{equation*}%
then $v$ unital, i.e. (\ref{form: v strongly unitary}), implies $\xi
(r\otimes 1_{R})=\varepsilon (r)1_{H}$.

Conversely, if $\xi (r\otimes 1)=\varepsilon (r)1_{H}$, then applying $%
\lambda $ to $\varepsilon (r)1_{H}={v}^{-1}(r^{(1)}\otimes 1_{R})r_{\langle
-1\rangle }^{(2)}{v}(r_{\langle 0\rangle }^{(2)}\otimes 1_{R})$ and using
the fact that $\lambda \Psi ({v})=\varepsilon _{C}$, we obtain ${v}%
^{-1}(r\otimes 1)=\varepsilon (r)$ so that ${v}(r\otimes 1)=\varepsilon (r)$
also.

The argument for elements $1 \otimes r$ is the same.
\end{proof}

The propositions above now prove Theorem \ref{teo: hope}.

\begin{proof}
(of Theorem \ref{teo: hope}.) Theorem \ref{teo: bijection} and Propositions %
\ref{pro: YD6'}, \ref{pro: YD7'}, \ref{pro: YD8'}, \ref{pro: YD10'} show
that $V=G(\Xi )$ is the set of $v\in \mathcal{G}$ satisfying (\ref{form: mv
colinear}) through (\ref{form: v strongly unitary}).
\end{proof}

\begin{remark}
\label{rem: associative}In general the map ${v}= F(\xi)$ in Theorem \ref%
{teo: hope}, is not a cocycle
\cite[Example 5.11, Remarks 5.13(ii)]{A.B.M.} although in these cases $(R^{{v} },m^v,u)$ is an
associative algebra. See also
\cite{ABM proceedings} for a discussion of when $F(\xi)$ is a
cocycle. In general, it is unknown whether $(R^{{v}},m^v,u)$ is an
associative algebra or not.
\end{remark}

\begin{proposition}
\label{lem: A B} Let $\xi \in \Xi$ and $v:= G(\xi)$ as above. Then the
following are equivalent:
\begin{eqnarray*}
& \mathrm{(i)}& \quad {v}\in Z_{H}^{2}\left( R,K\right) ; \\
& \mathrm{(ii)} & \quad \partial _{R}^{{2}}\left( v\right) =\varepsilon _{D}
\text{ with } D:= R \otimes R \otimes R; \\
& \mathrm{(iii)}& \quad \lambda \circ \left[ \Psi \left( \partial _{R}^{{2}%
}\left( v\right) \right) \right] =\varepsilon _{D}; \\
& \mathrm{(iv)} &\quad \partial _{R}^{{2}}\left( v\right) =\varepsilon _{D}
\text{ on } {^{co( H) }(D)}.
\end{eqnarray*}
\end{proposition}

\begin{proof}
Set $A:=R\#_\xi H$. Since $\partial _{R}^{{2}}\left( v\right) =
\varepsilon_D $ if and only if $v\in \mathrm{Ker}(\partial _{R}^2)=\Omega^2(%
\mathrm{Ker}(\partial _{A}^2))=\Omega^2(Z_{H}^{2}\left( A,K\right))$, the
equivalence of (i) and (ii) follows from the definition of $Z_{H}^{2}\left(
R,K\right)$ from \cite{A.B.M.} (or see Section \ref{sec: coh prebi}). To see
the equivalence of (ii) and (iii), apply $\lambda$ to both sides of (\ref%
{form: quasi cocycle}) to obtain $\lambda \circ \Psi \left( \partial _{R}^{{2%
}}\left( v\right) \right) =\partial _{R}^{{2}}\left( v\right)$.

Since (ii) implies (iv) is trivial, it remains to show that (iv) implies
(iii).

%Note that for $z \in D$, $(\lambda \otimes D) \rho_D(z) = \lambda
%\Psi (z) = \lambda \left( z_{\left\langle -1\right\rangle }\right)
%z_{\left\langle 0\right\rangle }.$ Then
%\begin{equation*}
%\rho(\lambda \Psi (z) ) = \lambda \left( z_{\left\langle -2\right\rangle
%}\right) z_{\left\langle -1\right\rangle }\otimes z_{\left\langle
%0\right\rangle }=1_{H}\lambda \left( z_{\left\langle -1\right\rangle
%}\right) \otimes z_{\left\langle 0\right\rangle }=1_{H}\otimes \lambda \Psi
%(z).
%\end{equation*}%
%Thus $\mathrm{Im}\left( \lambda \Psi
%\right)\textcolor[rgb]{1.00,0.00,0.00}{\subseteq}{^{co\left(
%H\right) }\textcolor[rgb]{1.00,0.00,0.00}{D} }$ and the proof is
%complete.

  Note that for $z\in D$,%
\begin{equation*}
(\lambda \circ \left[ \Psi \left( \partial _{R}^{{2}}\left( v\right) \right) %
\right] )(z)=\lambda (z_{\left\langle -1\right\rangle }\partial _{R}^{{2}%
}\left( v\right) (z_{\left\langle 0\right\rangle }))=\lambda
(z_{\left\langle -1\right\rangle })\partial _{R}^{{2}}\left( v\right)
(z_{\left\langle 0\right\rangle })=\partial _{R}^{{2}}\left( v\right)
(\lambda (z_{\left\langle -1\right\rangle })z_{\left\langle 0\right\rangle
}).
\end{equation*}
To complete the proof it suffices to check that $\lambda (z_{\left\langle
-1\right\rangle })z_{\left\langle 0\right\rangle }\in {^{co(H)}D.}$ Indeed,
we have
\begin{equation*}
\rho (\lambda (z_{\left\langle -1\right\rangle })z_{\left\langle
0\right\rangle })=\lambda \left( z_{\left\langle -2\right\rangle }\right)
z_{\left\langle -1\right\rangle }\otimes z_{\left\langle 0\right\rangle
}=1_{H}\lambda \left( z_{\left\langle -1\right\rangle }\right) \otimes
z_{\left\langle 0\right\rangle }=1_{H}\otimes \lambda (z_{\left\langle
-1\right\rangle })z_{\left\langle 0\right\rangle }.
\end{equation*}%

\end{proof}

\section{The dual quasi-bialgebra $A^{{v}_{A}}$}\label{sect: dual quasi}

In the previous section, for $(R,\xi)$ a pre-bialgebra with cocycle in $_H^H%
\mathcal{YD}$, the defining properties of $\xi$ were translated to
properties of $v:=G(\xi) \in \mathrm{Reg}_H(R \otimes R,K)$. In other words, we showed
that the functor from $\mathcal{R}$ to $\mathcal{R}^\prime$ which takes an object
$(R,\xi)$ to $(R, v:=(\lambda \xi)^{-1})$ is an isomorphism. In this section
we will first show that $R^v$, the pre-bialgebra $R$ with its multiplication
twisted by $v$, is a dual quasi-bialgebra in $_H^H\mathcal{YD}$ with
reassociator $\partial^2_R(v)$, in other words that $T_1$ is a functor.  First we recall the definitions of dual
quasi-bialgebra in vector spaces and in $_H^H\mathcal{YD}$ and define a
process of bosonization taking dual quasi-bialgebras in $_H^H\mathcal{YD}$
to dual quasi-bialgebras in vector spaces over $K$.

\subsection{Dual quasi-bialgebras and bosonization.} \label{sec: dualquasiboso}

Recall from \cite[page 66]{Majid} that a \emph{dual quasi-bialgebra} $\left(
D,m,u,\Delta ,\varepsilon ,\alpha \right) $ is a coalgebra $\left( D,\Delta
,\varepsilon \right) $ with coalgebra homomorphisms $m:D\otimes D\rightarrow
D$ and $u:K\rightarrow D$ ($1_{D}:=u\left( 1_{K}\right) $) such that $%
1_{D}x=x=x1_{D}$ for all $x\in D$ and $\alpha \in \mathrm{Reg}\left( {D}%
^{\otimes 3},K\right) $ is such that $\partial _{D}^{3}\left( \alpha
\right) =\varepsilon _{D\otimes D\otimes D}$, $\alpha $ is unital and
   $m$ is $\alpha$-associative.
%This amounts to:
Equivalently,
\begin{eqnarray*}
& (i)&\quad \alpha \left( D\otimes D\otimes m\right) \ast \alpha \left(
m\otimes D\otimes D\right) = \left( \varepsilon \otimes \alpha \right) \ast
\alpha \left( D\otimes m\otimes D\right) \ast \left( \alpha \otimes
\varepsilon \right) \\
& (ii)&\quad \alpha \left( D\otimes 1_{D}\otimes D\right) = \alpha \left(
1_{D}\otimes D \otimes D\right) = \alpha \left( D\otimes D \otimes
1_D\right) = \varepsilon _{D\otimes D}, \\
& (iii)& \quad m\left( D\otimes m\right) \ast \alpha =\alpha \ast m\left(
m\otimes D\right) .
\end{eqnarray*}%
Note that in (ii) any of the three equalities such as $\alpha \left(
1_{D}\otimes D\otimes D\right) =\varepsilon _{D\otimes D} $ implies that $%
\alpha$ is unital. The map $\alpha$ is called the reassociator.

A unital map ${v}\in \mathrm{Reg}(D \otimes D,K)$ is called a \emph{gauge
transformation}. Then the twisted dual quasi-bialgebra $D^{{v}}=\left( D,m^v:=v\ast m \ast v^{-1},u,\Delta
,\varepsilon ,\alpha _{D^{{v}}}\right) $ is also a dual quasi-bialgebra
where the reassociator $\alpha_{D_v}$ is defined by:
\begin{equation}  \label{form: reassociator}
\alpha _{D^{{v}}}:=\left( \varepsilon \otimes {v}\right) \ast {v}\left(
D\otimes m\right) \ast \alpha \ast {v}^{-1}\left( m\otimes D\right) \ast
\left( {v}^{-1}\otimes \varepsilon \right) .
\end{equation}%
Note that, whenever $\alpha$ is trivial, one has $\alpha _{D^{{v}}}=\partial^2_D(v)$, in the sense of Subsection \ref{subs: coalgmultunit}.

A morphism of dual quasi-bialgebras $f:\left( D,m,u,\Delta ,\varepsilon
,\alpha \right) \rightarrow \left( D^{\prime },m^{\prime },u^{\prime
},\Delta ^{\prime },\varepsilon ^{\prime },\alpha ^{\prime }\right) $ is a
coalgebra homomorphism $f:\left( D,\Delta ,\varepsilon \right) \rightarrow
\left( D^{\prime },\Delta ^{\prime },\varepsilon ^{\prime }\right) $ such
that $m^{\prime }(f\otimes f)=fm,fu=u^{\prime }$ and $\alpha ^{\prime
}\left( f\otimes f\otimes f\right) =\alpha .$ It is an isomorphism of
quasi-bialgebras if, in addition, it is invertible.

\begin{definition}
Dual quasi-bialgebras $A$ and $B$ are called \emph{quasi-isomorphic} (or
equivalent) whenever $A \cong B^v$ as dual quasi-bialgebras for some gauge
transformation ${v}\in \left( B\otimes B\right) ^{\ast }$.
\end{definition}

We now give the definition of a dual quasi-bialgebra $Q$ in $_{H}^{H}%
\mathcal{YD}$ and show that a $K$-dual quasi-bialgebra can be constructed
from $Q$ by bosonization. Although our purpose in the end is to study Hopf
algebras whose coradicals are semisimple sub-Hopf algebras, this result is
interesting on its own and adds to the literature on constructions with dual
quasi-bialgebras.
(See, for example, \cite[Section
3]{BulacuNau} where a smash product $B\# H$ with a quasi-Hopf
structure is studied  for $H$ a quasi-Hopf algebra and $B$ a braided
Hopf algebra in $_H^H\mathcal{YD}$.)

\begin{definition}
\label{def: dual quasi braided} Let $H$ be a Hopf algebra. A \emph{dual quasi-bialgebra} $(Q,m,u,\Delta
,\varepsilon ,\alpha )$ in the braided monoidal category $^H_H\mathcal{YD}$ is a coalgebra $\left( Q,\Delta ,\varepsilon \right) $ in $^H_H\mathcal{YD}$ together with coalgebra homomorphisms $m:Q\otimes Q\rightarrow
Q $ and $u:K\rightarrow Q$ in $^H_H\mathcal{YD}$ and a convolution invertible
element $\alpha \in {^H_H\mathcal{YD}}\left( Q^{\otimes 3},K\right) $ (braided reassociator) such that%
\begin{eqnarray}
&&\alpha \left( Q\otimes Q\otimes m\right) \ast \alpha \left( m\otimes
Q\otimes Q\right) =\left( \varepsilon \otimes \alpha \right) \ast \alpha
\left( Q\otimes m\otimes Q\right) \ast \left( \alpha \otimes \varepsilon
\right) ,  \label{form: alpha 3-cocycle} \\
&&\alpha \left( Q\otimes u\otimes Q\right) =\alpha \left( u\otimes Q\otimes
Q\right) =\alpha \left( Q\otimes Q\otimes u\right) =\varepsilon _{Q\otimes
Q},  \label{form: alpha unital} \\
&&m\left( Q\otimes m\right) \ast \alpha =\alpha \ast m\left( m\otimes
Q\right) ,  \label{form: m quasi assoc} \\
&&m\left( u\otimes Q\right) =\mathrm{Id}_{Q}=m\left( Q\otimes u\right) .
\label{form: m unital}
\end{eqnarray}%
Note that in (\ref{form: alpha unital}) any of the three equalities such as $%
\alpha \left( u\otimes Q\otimes Q\right) =\varepsilon _{Q\otimes Q}$ implies
that $\alpha $ is unital.
\end{definition}

The next result proves there is a dual quasi-bialgebra associated to any
dual quasi-bialgebra $Q$ in the braided monoidal category ${_{H}^{H}%
\mathcal{YD}}$.

\begin{proposition}
\label{teo: smash quasi YD}Let $H$ be a Hopf algebra and let $\left(
Q,m,u,\Delta ,\varepsilon ,\alpha \right) $ be a dual quasi-bialgebra in ${%
_{H}^{H}\mathcal{YD}}\mathbf{.}$ Set $B:=Q\otimes H.$ Then
\begin{equation*}
Q\#H:=\left( B,m_{B},u_{B},\Delta _{B},\varepsilon _{B},\alpha _{B}\right)
\end{equation*}%
is an ordinary dual quasi-bialgebra where $Q\#H$ has the usual coalgebra
structure and multiplication and unit maps namely%
\begin{gather*}
m_{B}\left( r\#h\otimes s\#l\right) :=m\left( r\otimes h_{(1)}s\right)
\#h_{(2)}l,\qquad u_{B}\left( 1_{K}\right) :=1_{Q}\#1_{H}, \\
\Delta _{B}\left( r\#h\right) := ( r^{(1)}\#r_{\langle -1\rangle
}^{(2)}h_{(1)} ) \otimes ( r_{\langle 0\rangle }^{(2)}\#h_{\left( 2\right) }
) ,\qquad \varepsilon _{B}\left( r\#h\right) :=\varepsilon \left( r\right)
\varepsilon _{H}\left( h\right) ,
\end{gather*}%
and the reassociator is given by $\alpha _{B}=\mho _{H,Q}^{3}\left( \alpha
\right) $ namely%
\begin{equation*}
\alpha _{B}\left( r\#h\otimes s\#l\otimes t\#k\right) :=\alpha \left(
r\otimes h_{\left( 1\right) }s\otimes h_{\left( 2\right) }lt\right)
\varepsilon _{H}\left( k\right) .
\end{equation*}
\end{proposition}

\begin{proof}
Since $\left( Q,\Delta ,\varepsilon \right) $ is a coalgebra in ${_{H}^{H}%
\mathcal{YD}}$, $\left( B,\Delta _{B},\varepsilon _{B}\right) $ is the smash
coproduct of $Q$ by $H$. It is straightforward to verify that $%
u_{B}:K\rightarrow B$ is a coalgebra map and that $\varepsilon
_{B}m_{B}\left( r\#h\otimes s\#l\right) =\varepsilon \left( r\right)
\varepsilon _{H}\left( h\right) \varepsilon \left( s\right) \varepsilon
_{H}\left( l\right) $; these are left to the reader. Similarly it is clear
that $m_{B}\left( 1_{B}\otimes s\#l\right) =s\#l$ and also $m_{B}\left(
r\#h\otimes 1_{B}\right) =m\left( r\otimes h_{(1)}1_{Q}\right)
\#h_{(2)}=m\left( r\otimes 1_{Q}\right) \#h=r\#h.$ The map $m_B$ is clearly
right $H$-linear and is also left $H$-linear since
\begin{eqnarray*}
km_{B} ( r\#h\otimes s\#l ) &=&k_{\left( 1\right) }m ( r\otimes h_{(1)}s )
\#k_{\left( 2\right) }h_{(2)}l=m ( k_{\left( 1\right) }r\otimes k_{\left(
2\right) }h_{(1)}s ) \#k_{\left( 3\right) }h_{(2)}l \\
&=&m_{B} ( k_{\left( 1\right) }r\#k_{\left( 2\right) }h\otimes s\#l ) =m_{B}
( k(r\#h\otimes s\#l) ).
\end{eqnarray*}%
Furthermore, $m_B$ is $H$-balanced since
\begin{eqnarray*}
m_{B}\left( \left( r\#h\right) k\otimes s\#l\right) &=&m_{B}\left(
r\#hk\otimes s\#l\right) =m\left( r\otimes h_{(1)}k_{\left( 1\right)
}s\right) \#h_{(2)}k_{\left( 2\right) }l \\
&=&m_{B}\left( r\#h\otimes k_{\left( 1\right) }s\#k_{\left( 2\right)
}l\right) =m_{B}\left( r\#h\otimes k\left( s\#l\right) \right) .
\end{eqnarray*}

We check next that $m_{B}:B\otimes B\rightarrow B$ is a coalgebra
homomorphism. Since $\Delta _{B}$ is $H$-bilinear and $m_{B}$ is $H$%
-bilinear and $H$-balanced, it suffices to show that $m_{B}$ is a coalgebra
homomorphism on elements of the form $r\#1_{H}\otimes s\#1_{H} $. In the
third step we use the fact that $m:Q\otimes Q\rightarrow Q$ is a coalgebra
map and $\varepsilon _{B}m_{B}=\varepsilon _{B\otimes B}.$
\begin{eqnarray*}
\Delta _{B}m_{B}\left( r\#1_{H}\otimes s\#1_{H}\right) &=& \Delta _{B}\left[
m\left( r\otimes s\right) \#1_{H}\right]
\\
&=& \left[ m\left( r\otimes s\right) \right] ^{(1)}\#\left[ m\left( r\otimes
s\right) \right] _{\langle -1\rangle }^{(2)}\otimes \left[ m\left( r\otimes
s\right) \right] _{\langle 0\rangle }^{(2)}\#1_{H}
\\
&=& m(r^{\left( 1\right) }\otimes r_{\left\langle -1\right\rangle }^{\left(
2\right) }s^{\left( 1\right) })\#[m(r_{\left\langle 0\right\rangle }^{\left(
2\right) }\otimes s^{\left( 2\right) })]_{\langle -1\rangle }\otimes \lbrack
m(r_{\left\langle 0\right\rangle }^{\left( 2\right) }\otimes s^{\left(
2\right) })]_{\langle 0\rangle }\#1_{H}
 \\
&=& m ( r^{\left( 1\right) }\otimes r_{\left\langle -1\right\rangle
}^{\left( 2\right) }s^{\left( 1\right) } ) \# ( r_{\left\langle
0\right\rangle }^{\left( 2\right) }\otimes s^{\left( 2\right) } ) _{\langle
-1\rangle }\otimes m [ ( r_{\left\langle 0\right\rangle }^{\left( 2\right)
}\otimes s^{\left( 2\right) } ) _{\langle 0\rangle } ] \#1_{H}
\\
&=&m(r^{\left( 1\right) }\otimes r_{\left\langle -2\right\rangle }^{\left(
2\right) }s^{\left( 1\right) })\#r_{\left\langle -1\right\rangle }^{\left(
2\right) }s_{\left\langle -1\right\rangle }^{\left( 2\right) }\otimes
m(r_{\left\langle 0\right\rangle }^{\left( 2\right) }\otimes s_{\left\langle
0\right\rangle }^{\left( 2\right) })\#1_{H}
\\
&=&m\left( r^{\left( 1\right) }\otimes r_{\left\langle -2\right\rangle
}^{\left( 2\right) }s^{\left( 1\right) }\right) \#r_{\left\langle
-1\right\rangle }^{\left( 2\right) }s_{\left\langle -1\right\rangle
}^{\left( 2\right) }\otimes m_{B}\left[ r_{\left\langle 0\right\rangle
}^{\left( 2\right) }\#1_{H}\otimes s_{\left\langle 0\right\rangle }^{\left(
2\right) }\#1_{H}\right]
\\
&=&m_{B}\left[ r^{\left( 1\right) }\#r_{\left\langle -1\right\rangle
}^{\left( 2\right) }\otimes s^{\left( 1\right) }\#s_{\left\langle
-1\right\rangle }^{\left( 2\right) }\right] \otimes m_{B}\left[
r_{\left\langle 0\right\rangle }^{\left( 2\right) }\#1_{H}\otimes
s_{\left\langle 0\right\rangle }^{\left( 2\right) }\#1_{H}\right] \\
&=&m_{B}\left[ \left( r\#1_{H}\right) _{(1)}\otimes \left( s\#1_{H}\right)
_{(1)}\right] \otimes m_{B}\left[ \left( r\#1_{H}\right) _{(2)}\otimes
\left( s\#1_{H}\right) _{(2)}\right] .
\end{eqnarray*}%
It is clear that $\alpha _{B} = \mho^3_{H,R}(\alpha)$ is convolution
invertible with inverse $\mho _{H,R}^{3}\left( \alpha ^{-1}\right) .$

If $\phi$ is any $H$-multibalanced, $H$-bilinear map from $B \otimes B
\otimes B$ to $B$, then
\begin{equation*}
\phi \left( r\#h\otimes s\#l\otimes t\#k\right) =\phi \left( r\#1_{H}\otimes
h_{\left( 1\right) }s\#1_{H}\otimes h_{\left( 2\right) }l_{\left( 1\right)
}t\#1_{H}\right) h_{\left( 3\right) }l_{\left( 2\right) }k.
\end{equation*}
Since $m_{B}$ is $H$-balanced and $H$-bilinear, then $m_{B}^{l}:=m_{B}\left(
m_{B}\otimes B\right) $ and $m_{B}^{r}:=m_{B}\left( B\otimes m_{B}\right) $
are $H$-multibalanced and $H$-bilinear too and so are $m_{B}^{r}\ast \alpha
_{B}$ and $\alpha _{B}\ast m_{B}^{l}$. Thus it suffices to check that $%
m_{B}^{r}\ast \alpha _{B}= \alpha _{B}\ast m_{B}^{l}$, i.e., $m_{B}$
is associative up to multiplication with the reassociator, on
elements of the form $r\#1_{H}\otimes s\#1_{H}\otimes t\#1_{H}$ for
$r,s,t\in R $:
\begin{eqnarray*}
&&\left[ m_{B}^{r}\ast \alpha _{B}\right] \left( r\#1_{H}\otimes
s\#1_{H}\otimes t\#1_{H}\right)
 \\
&=&m_{B}^{r}(r^{(1)}\#r_{\langle -1\rangle }^{(2)}\otimes
s^{(1)}\#s_{\langle -1\rangle }^{(2)}\otimes t^{(1)}\#t_{\langle -1\rangle
}^{(2)})\alpha _{B}(r_{\langle 0\rangle }^{(2)}\#1_{H}\otimes s_{\langle
0\rangle }^{(2)}\#1_{H}\otimes t_{\langle 0\rangle }^{(2)}\#1_{H})
\\
 &=&m_{B}^{r}(r^{(1)}\#r_{\langle -1\rangle }^{(2)}\otimes
s^{(1)}\#s_{\langle -1\rangle }^{(2)}\otimes t^{(1)}\#t_{\langle
-1\rangle }^{(2)})\alpha (r_{\langle 0\rangle }^{(2)}\otimes
s_{\langle 0\rangle }^{(2)}\otimes t_{\langle 0\rangle }^{(2)})
 \\
&=&m_{B}^{r}(r^{(1)}\# 1 \otimes  r_{\langle -3\rangle
}^{(2)}s^{(1)} \# 1  \otimes r_{\langle -2\rangle }^{(2)} s_{\langle
-2\rangle }^{(2)} t^{(1)}\# 1 )r_{\langle -1\rangle
}^{(2)}s_{\langle -1\rangle }^{(2)}t_{\langle -1\rangle
}^{(2)}\alpha (r_{\langle 0\rangle }^{(2)}\otimes s_{\langle
0\rangle }^{(2)}\otimes t^{(2)}_{\langle 0\rangle })
 \\
&=&m_{B} (r^{(1)}\# 1 \otimes  m(r_{\langle -3\rangle }^{(2)}s^{(1)}
\otimes r_{\langle -2\rangle }^{(2)} s_{\langle -2\rangle }^{(2)}
t^{(1)})\# 1 )r_{\langle -1\rangle }^{(2)}s_{\langle -1\rangle
}^{(2)}t_{\langle -1\rangle }^{(2)}\alpha (r_{\langle 0\rangle
}^{(2)}\otimes s_{\langle 0\rangle }^{(2)}\otimes t_{\langle
0\rangle }^{(2)})
\\
&=& [m  (r^{(1)}  \otimes  m(r_{\langle -3\rangle }^{(2)}s^{(1)}
\otimes r_{\langle -2\rangle }^{(2)} s_{\langle -2\rangle }^{(2)}
t^{(1)})]\#  r_{\langle -1\rangle }^{(2)}s_{\langle -1\rangle
}^{(2)}t_{\langle -1\rangle }^{(2)}\alpha (r_{\langle 0\rangle
}^{(2)}\otimes s_{\langle 0\rangle }^{(2)}\otimes t^{(2)})
\\
& =& m \left( Q\otimes m \right) \left[ r^{(1)}\otimes r_{\langle
-2\rangle }^{(2)}s^{(1)}\otimes r_{\langle -1\rangle
}^{(2)}s_{\langle -1\rangle }^{(2)}t^{(1)}\right] \# \Psi(\alpha)
\left( r_{\langle 0\rangle }^{(2)}\otimes s_{\langle 0\rangle
}^{(2)}\otimes t^{(2)}\right)
\\
& \overset{ ( \ref{form: quasi cocycle} )}{=} & m \left( Q\otimes m
\right) \left[ r^{(1)}\otimes r_{\langle -2\rangle
}^{(2)}s^{(1)}\otimes r_{\langle -1\rangle }^{(2)}s_{\langle
-1\rangle }^{(2)}t^{(1)}\right] \#
 \alpha  \left( r_{\langle 0\rangle }^{(2)}\otimes s_{\langle
0\rangle }^{(2)}\otimes t^{(2)}\right)
\\
&=&\left[ m\left( Q\otimes m\right) \ast \alpha \right] \left(
r\otimes s\otimes t\right) \#1_{H}
\\
&\overset{(\ref{form: m quasi assoc})}{=}&\left[ \alpha \ast m\left(
m
\otimes Q\right) \right] \left( r\otimes s\otimes t\right) \#1_{H} \\
&=&\alpha \left( r^{(1)}\otimes r_{\langle -2\rangle
}^{(2)}s^{(1)}\otimes r_{\langle -1\rangle }^{(2)}s_{\langle
-1\rangle }^{(2)}t^{(1)}\right) m\left( m\otimes Q\right) \left(
r_{\langle 0\rangle }^{(2)}\otimes s_{\langle 0\rangle
}^{(2)}\otimes t^{(2)}\right) \#1_{H}
\\
&=&\alpha \left( r^{(1)}\otimes r_{\langle -2\rangle
}^{(2)}s^{(1)}\otimes r_{\langle -1\rangle }^{(2)}s_{\langle
-1\rangle }^{(2)}t^{(1)}\right) m \left[ m\left( r_{\langle 0\rangle
}^{(2)}\otimes s_{\langle 0\rangle }^{(2)}\right) \otimes
t^{(2)}\right] \#1_{H}
\\
&=&\alpha \left( r^{(1)}\otimes r_{\langle -2\rangle
}^{(2)}s^{(1)}\otimes r_{\langle -1\rangle }^{(2)}s_{\langle
-1\rangle }^{(2)}t^{(1)}\right) m_{B}\left( m \left( r_{\langle
0\rangle }^{(2)}\otimes s_{\langle 0\rangle }^{(2)}\right)
\#1_{H}\otimes t^{(2)}\#1_{H}\right)
\\
&=&\alpha \left( r^{(1)}\otimes r_{\langle -2\rangle
}^{(2)}s^{(1)}\otimes r_{\langle -1\rangle }^{(2)}s_{\langle
-1\rangle }^{(2)}t^{(1)}\right) m_{B}^{l}\left( r_{\langle 0\rangle
}^{(2)}\#1_{H}\otimes s_{\langle 0\rangle }^{(2)}\#1_{H}\otimes
t^{(2)}\#1_{H}\right)
\\
&=&\alpha _{B}\left( r^{(1)}\#r_{\langle -1\rangle }^{(2)}\otimes
s^{(1)}\#s_{\langle -1\rangle }^{(2)}\otimes t^{(1)}\#t_{\langle
-1\rangle }^{(2)}\right) m_{B}^{l}\left( r_{\langle 0\rangle
}^{(2)}\#1_{H}\otimes s_{\langle 0\rangle }^{(2)}\#1_{H}\otimes
t_{\langle 0\rangle }^{(2)}\#1_{H}\right)
\\
&=&\alpha _{B}\left( \left( r\#1_{H}\right) _{\left( 1\right)
}\otimes \left( s\#1_{H}\right) _{\left( 1\right) }\otimes \left(
t\#1_{H}\right) _{\left( 1\right) }\right) m_{B}^{l}\left( \left(
r\#1_{H}\right) _{\left( 2\right) }\otimes \left( s\#1_{H}\right)
_{\left( 2\right) }\otimes \left( t\#1_{H}\right) _{\left( 2\right)
}\right)
\\
&=&(\alpha _{B}\ast m_{B}^{l})\left( r\#1_{H}\otimes s\#1_{H}\otimes
t\#1_{H}\right).
\end{eqnarray*}

It remains to prove the cocycle condition for $\alpha _{B}$. First we note
that the following equalities can be easily checked by applying
$\Omega _{H,Q}^{4}$ to both sides.
\begin{eqnarray}
&& \alpha _{B}\left( B\otimes B\otimes m_{B}\right) = \mho _{H,Q}^{4}[\alpha
\left( Q\otimes Q\otimes m\right) ],  \label{form: 3cocycle1}
\\
&&\alpha _{B}\left( m_{B}\otimes B\otimes B\right) = \mho _{H,Q}^{4}\left[
\alpha \left( m\otimes Q\otimes Q\right) \right] ,  \label{form: 3cocycle2}
\\
&&\alpha _{B}\left( B\otimes m_{B}\otimes B\right) = \mho _{H,Q}^{4}\left(
\alpha \left( Q\otimes m\otimes Q\right) \right) ,  \label{form: 3cocycle4}
\\
&& \varepsilon _{B}\otimes \alpha _{B}= \mho _{H,Q}^{4}\left( \varepsilon
\otimes \alpha \right) ,  \label{form: 3cocycle3} \text{ and  }
\alpha _{B}\otimes \varepsilon _{B}= \mho _{H,Q}^{4}\left( \alpha \otimes
\varepsilon \right) .  \label{form: 3cocycle5}
\end{eqnarray}
Therefore we get%
\begin{eqnarray*}
&&\alpha _{B}\left( B\otimes B\otimes m_{B}\right) \ast \alpha _{B}\left(
m_{B}\otimes B\otimes B\right) \\
&\overset{(\ref{form: 3cocycle1}),(\ref{form: 3cocycle2})}{=}&\mho _{H,Q}^{4}%
\left[ \alpha \left( Q\otimes Q\otimes m\right) \right] \ast \mho _{H,Q}^{4}%
\left[ \alpha \left( m\otimes Q\otimes Q\right) \right] \\
&=&\mho _{H,Q}^{4}\left[ \alpha \left( Q\otimes Q\otimes m\right) \ast
\alpha \left( m\otimes Q\otimes Q\right) \right] \\
&\overset{(\ref{form: alpha 3-cocycle})}{=}&\mho _{H,Q}^{4}\left[ \left(
\varepsilon \otimes \alpha \right) \ast \alpha \left( Q\otimes m\otimes
Q\right) \ast \left( \alpha \otimes \varepsilon \right) \right] \\
&=&\mho _{H,Q}^{4}\left( \varepsilon \otimes \alpha \right) \ast \mho
_{H,Q}^{4}\left( \alpha \left( Q\otimes m\otimes Q\right) \right) \ast \mho
_{H,Q}^{4}\left( \alpha \otimes \varepsilon \right)
\\
&\overset{(\ref{form: 3cocycle4}),(\ref{form: 3cocycle3})}
{=}&\left( \varepsilon _{B}\otimes \alpha _{B}\right) \ast
\alpha _{B}\left( B\otimes m_{B}\otimes B\right) \ast \left( \alpha
_{B}\otimes \varepsilon _{B}\right)
\end{eqnarray*}

Hence we have proved the cocycle condition for $\alpha _{B}$ as desired.
\end{proof}

\begin{definition}
Let $H$ be a Hopf algebra and let $\left( Q,m,u,\Delta ,\varepsilon ,\alpha
\right) $ be a dual quasi-bialgebra in ${_{H}^{H}\mathcal{YD}}\mathbf{.}$
Set $B:=Q\otimes H.$ Then
\begin{equation*}
Q\#H:=\left( B,m_{B},u_{B},\Delta _{B},\varepsilon _{B},\alpha _{B}\right)
\end{equation*}%
will be called the \emph{bosonization }of the dual quasi-bialgebra $\left(
Q,m,u,\Delta ,\varepsilon ,\alpha \right) $ in ${_{H}^{H}\mathcal{YD}}$ by $%
H $.
\end{definition}

If $Q$ is connected, by the same proof as \cite[Theorem 3.9]{A.M.St.-Small}, the coradical
of $Q \# H $ is $K \otimes H$. Thus we have shown that the second bosonization functor $B_2$ mentioned in the Introduction
actually maps objects in
$\mathcal{Q}$ to objects in $\mathcal{B}$.  The fact that $B_2$ preserves morphisms is straightforward and so we omit it.
\subsection{The dual quasi-bialgebra from a pre-bialgebra with cocycle}

In this section we show that if $(R,\xi)$ is a pre-bialgebra with cocycle,
then for $v:= G(\xi)$ with $G$ the map from Theorem \ref{teo: hope},
twisting $R$ by $v$ gives $R^v$ the structure of a dual quasi-bialgebra in ${%
_{H}^{H}\mathcal{YD}}$.  In other words, we show that the first twisting functor $T_1$ from the introduction
maps objects in $\mathcal{R}$ to objects in $\mathcal{Q}$. Again, we leave the verification that $T_1$ also preserves morphisms to the reader.

\begin{proposition}
\label{pro: R dual quasi}Let $(R,\xi )$ be a connected pre-bialgebra with
cocycle in ${_{H}^{H}\mathcal{YD}}$. Let $v:=G(\xi )=\left( \lambda \xi
\right) ^{-1}$. Then $R^v:=\left( R,m^v,u,\Delta ,\varepsilon ,\alpha:=\partial
_{R}^{2}\left( v\right) \right) $ is a connected dual quasi-bialgebra in ${_{H}^{H}%
\mathcal{YD}}$.
\end{proposition}

\begin{proof}
By construction, $(R,\Delta, \varepsilon)$
is a connected coalgebra in ${_{H}^{H}\mathcal{YD}}$ and $m^v$ and $u$ are coalgebra
homomorphisms. From Corollary \ref{co:Delta2R},
\begin{equation*}
\alpha := \partial _{R}^{2}\left( v\right) ={v}\left( R\otimes m^v
\right) \ast \left(   \varepsilon  \otimes {v}\right) \ast \left(
{v}^{-1}\otimes   \varepsilon  \right) \ast {v}^{-1}(m^v\otimes R),
\end{equation*}
so that $\alpha $ is convolution invertible and by Remark \ref{rem: alpha YD}%
, $\alpha $ is in ${\ _{H}^{H}\mathcal{YD}}$. It remains to check that $%
\alpha$ satisfies (\ref{form: alpha 3-cocycle}) through (\ref{form: m unital}%
). It is straightforward to show that unitality of $v,v^{-1}$ implies (\ref%
{form: alpha unital}), the unital property for $\alpha $. Theorem \ref{teo:
hope} implies that $m^v\ $satisfies (\ref{form: m quasi assoc}). It is
straightforward to prove (\ref{form: m unital}). It remains to verify (\ref%
{form: alpha 3-cocycle}). To simplify notation in the following computation,
we set:
\begin{eqnarray*}
&& m_{ {v}}^{1} :=m^{{v}}\otimes R\otimes R,\quad m_{ {v}}^{2}:=R\otimes m^{
{v}}\otimes R,\quad m_{ {v}}^{3}:=R\otimes R\otimes m^{ {v}} ; \\
&& m_{ {v}}^{r} :=m^{ {v}}\left( R\otimes m^{ {v}}\right) ,\quad m_{ {v}%
}^{l}:=m^{ {v}}\left( m^{ {v}}\otimes R\right); \\
&& \alpha_+:=( \partial _{R}^{{2}}) _{+}(v) = {v}\left( R\otimes
m^{v}\right) \ast \left( \varepsilon \otimes {v}\right) ,\quad \alpha_-:= (
\partial _{R}^{{2}}) _{-}(v^{-1}) = \left( {v}^{-1}\otimes   \varepsilon
 \right) \ast {v}^{-1}(m^{ {v}}\otimes R).
\end{eqnarray*}%
Now note that by (\ref{form: alfabeta}),
\begin{eqnarray*}
\alpha m^3_v &=& v (R\otimes m^{ {v}})m_{ {v}}^{3}\ast (\varepsilon \otimes {%
v})m_{ {v}}^{3}\ast ({v}^{-1}\otimes \varepsilon )m_{ {v}}^{3}\ast {v}%
^{-1}(m^{ {v}}\otimes R)m_{ {v}}^{3} \\
&=& {v}(R\otimes m_{ {v}}^{r})\ast (\varepsilon \otimes {v}(R\otimes m^{ {v}%
}))\ast ({v}^{-1}\otimes \varepsilon_{R \otimes R})\ast {v}^{-1}(m^v \otimes
m^v)\hspace{1mm} \text{ and } \\
\alpha m^1_v &=& v(R\otimes m^v)m_{ v}^{1} \ast (\varepsilon \otimes v
)m_v^{1}\ast ({v}^{-1}\otimes \varepsilon )m_{ {v}}^{1}\ast {v }^{-1}(m^{ {v}%
}\otimes R)m_v^{1} \\
&=& {v}(m^v \otimes m^v) \ast (\varepsilon _{R\otimes R}\otimes {v})\ast ({v}%
^{-1}(m^{ {v}}\otimes R)\otimes \varepsilon _{R})\ast {v}^{-1}(m_{ {v}%
}^{l}\otimes R),
\end{eqnarray*}
so that $\alpha m^3_v \ast \alpha m^1_v$, the left hand side of (\ref{form:
alpha 3-cocycle}), is equal to
\begin{equation*}
{v}(R\otimes m_{ v }^{r})\ast (\varepsilon \otimes {v}(R\otimes m^v))\ast ({v%
}^{-1}\otimes v) \ast ({v}^{-1}(m^{ v}\otimes R)\otimes \varepsilon )\ast {v}
^{-1}(m_{ {v}}^{l}\otimes R).
\end{equation*}
Since
\begin{equation*}
(\varepsilon \otimes v(R \otimes m^v))\ast (\varepsilon \otimes \varepsilon
\otimes v) = \varepsilon \otimes [v(R \otimes m^v)\ast (\varepsilon \otimes
v)]= \varepsilon \otimes \alpha_+,
\end{equation*}
and
\begin{equation*}
(v^{-1} \otimes \varepsilon \otimes \varepsilon)\ast (v^{-1}(m^v \otimes R)
\otimes \varepsilon) = [ (v^{-1} \otimes \varepsilon) \ast (v^{-1}(m^v
\otimes R))] \otimes \varepsilon = \alpha_- \otimes \varepsilon ,
\end{equation*}
then $\alpha m^3_v \ast \alpha m^1_v$ equals:
\begin{eqnarray*}
&& {v}(R\otimes m_{ {v}}^{r})\ast (\varepsilon \otimes \alpha_{+} )\ast
(\alpha_{-} \otimes \varepsilon )\ast {v}^{-1}(m_{ {v}}^{l}\otimes R) \\
&\overset{(\ref{form: aurea})}{=}& {v}(R\otimes m_{v}^{r}\ast \alpha_{+}
)\ast (\alpha_{-} \otimes \varepsilon )\ast {v}^{-1}(m_{ {v}}^{l}\otimes R)
\\
&=& {v}(R\otimes m_{ {v}}^{r}\ast \alpha \ast (\alpha_{-} ) ^{-1})\ast
\alpha_{-} \otimes \varepsilon )\ast {v}^{-1}(m_{{v}}^{l}\otimes R) \\
&\overset{(\ref{form: m quasi assoc})}{=}& {v}(R\otimes \alpha \ast m_{{{v} }%
}^{l}\ast ( \alpha_{-} ) ^{-1})\ast (\alpha_{-}\otimes \varepsilon )\ast {v}
^{-1}(m^{l}_v \otimes R) \\
&\overset{(\ref{form: aurea})}{=}& (\varepsilon \otimes \alpha ) \ast {v}%
(R\otimes m_v^{l})\ast (\varepsilon \otimes ( \alpha_{-} ) ^{-1} \ast
(\alpha_{-} \otimes \varepsilon )\ast {v}^{-1}(m_v^{l}\otimes R) \\
&=& (\varepsilon \otimes \alpha )\ast {v}(R\otimes m^{ v})m_v^{2} \ast
(\varepsilon \otimes (\alpha_{-} ) ^{-1}) \ast (\alpha_{-}\otimes
\varepsilon )\ast {v}^{-1}(m_{ {v}}^{l}\otimes R) \\
&=& (\varepsilon _{R}\otimes \alpha )\ast \left[ \alpha_{+} \ast
(\varepsilon \otimes {v}^{-1})\right] m_{{v}}^{2}\ast (\varepsilon \otimes
(\alpha_{-} ) ^{-1}) \ast (\alpha_{-} \otimes \varepsilon )\ast {v}^{-1}(m_{{%
v}}^{l}\otimes R) \\
&\overset{(\ref{form: alfabeta})}{=}& (\varepsilon \otimes \alpha )\ast
\alpha_{+} m_{v}^{2}\ast (\varepsilon \otimes {v}^{-1})m_{v}^{2} \ast
(\varepsilon \otimes (\alpha_{-} )^{-1})\ast (\alpha_{-} \otimes \varepsilon
)\ast {v }^{-1}(m_{v}^{l}\otimes R)
\end{eqnarray*}

Moreover%
\begin{eqnarray*}
&&( \varepsilon \otimes {v}^{-1}) m_{ {v}}^{2}\ast ( \varepsilon \otimes (
\alpha_{-} ) ^{-1}) \ast ( \alpha_{-} \otimes \varepsilon ) \ast {v} ^{-1}(
m_{ {v}}^{l}\otimes R) \\
&=& ( \varepsilon \otimes {v}^{-1}( m^{v}\otimes R) ) \ast ( \varepsilon
\otimes ( \alpha _{-} ) ^{-1}) \ast ( \alpha_{-} \otimes \varepsilon ) \ast {%
v}^{-1}( m_{{v}}^{l}\otimes R) \\
&=& ( \varepsilon \otimes \lbrack {v}^{-1}( m_{{v}}\otimes R) \ast ( \alpha
_{-} ) ^{-1}]) \ast ( \alpha _{-} \otimes \varepsilon ) \ast {v}^{-1}(
m_{v}^{l}\otimes R) \\
&=&( \varepsilon \otimes {v}\otimes \varepsilon ) \ast ( \alpha _{-} \otimes
\varepsilon ) \ast {v}^{-1}( m_{v}^{l}\otimes R) \\
&\overset{(\ref{form: aurina})}{=}& ( \varepsilon \otimes {v}\otimes
\varepsilon ) \ast {v}^{-1}( [ \alpha_{-} \ast m_{{v}}^{l} ] \otimes R) \\
&\overset{(\ref{form: aurina})}{=}&{v}^{-1}(\left[ ( \varepsilon \otimes {v}%
) \ast \alpha _{-} \ast m_{v}^{l})\right] \otimes R) \\
&=&{v}^{-1}(\left[ ( \varepsilon \otimes {v}) \ast (\alpha _{+} )^{-1}\ast
\alpha_{+} \ast \alpha_{-} \ast m_{v}^{l} \right] \otimes R) \\
&=& {\ v}^{-1}(\left[ ( \varepsilon \otimes {v}) \ast ( \alpha _{+}
)^{-1}\ast \alpha \ast m_{{v}}^{l} \right ] \otimes R) \\
& \overset{(\ref{form: m quasi assoc})}{=}& {v}^{-1}([ ( \varepsilon \otimes
{v}) \ast ( \alpha _{+} )^{-1}\ast m_{v}^{r}\ast \alpha ] \otimes R) \\
&=&{v}^{-1}(\left[ {v}^{-1}( R \otimes m^{v}) \ast m_{v}^{r}\ast \alpha %
\right] \otimes R) \\
&\overset{(\ref{form: aurina})}{=}& ( {v}^{-1}( R\otimes m^{v}) \otimes
  \varepsilon  ) \ast {v}^{-1}( [ m_{v}^{r}\ast \alpha ] \otimes R) \\
&\overset{(\ref{form: aurina2})}{=}& ( {v}^{-1}( R\otimes m^{v}) \otimes
\varepsilon ) \ast {v}^{-1}(m_{ v}^{r}\otimes R)\ast ( \alpha \otimes
\varepsilon ) \\
&=& \left[ ( {v}^{-1}\otimes   \varepsilon  ) m_{v}^{2}\right] \ast \left[ {%
v}^{-1}( m^{v}\otimes R) m_{v}^{2}\right] \ast ( \alpha \otimes \varepsilon )
\\
&\overset{(\ref{form: alfabeta})}{=}& \left[ \alpha _{-} m_{v}^{2}\right]
\ast ( \alpha \otimes \varepsilon ).
\end{eqnarray*}

Thus%
\begin{equation*}
( \varepsilon \otimes {v}^{-1}) m_{v}^{2}\ast  ( \varepsilon \otimes
( \alpha_{-} ) ^{-1}) \ast ( \alpha_{-} \otimes \varepsilon )  \ast
{v} ^{-1}( m_{v}^{l}\otimes R) = ( \alpha _{-} m_{v}^{2} ) \ast (
\alpha \otimes \varepsilon ),
\end{equation*}

and so
\begin{eqnarray*}
&& \alpha m_{v}^{3}\ast \alpha m_{v}^{1}= (\varepsilon \otimes \alpha )\ast
\alpha_{+} m_{v}^{2}\ast (\varepsilon \otimes {v}^{-1})m_{v}^{2} \ast
  (\varepsilon \otimes (\alpha_{-} ) ^{-1})\ast (\alpha_{-} \otimes
\varepsilon ) \ast {v }^{-1}(m_v^{l}\otimes R) \\
&=& (\varepsilon \otimes \alpha )\ast   (\alpha_{+} m_{v }^{2})\ast
(\alpha_{-}m_{v}^{2}) \ast (\alpha \otimes \varepsilon )
\overset{(\ref{form: alfabeta})}{=} (\varepsilon \otimes \alpha
)\ast \alpha m_{v}^{2}\ast (\alpha \otimes \varepsilon ),
\end{eqnarray*}
and $\alpha$ satisfies the $3$-cocycle condition.
\end{proof}

\begin{remark}
\label{rem: unknown}It is unknown whether   $\left( R,m^v ,u,\Delta
,\varepsilon \right) $  is a
braided bialgebra in the braided monoidal category ${_{H}^{H}\mathcal{YD}%
}$, see Remark \ref{rem: associative}.
\end{remark}

\subsection{The bosonization of $R^{{v}} $ with $H$.}\label{subs: boso}

Let $(R,\xi)$ be a connected pre-bialgebra with cocycle in ${_{H}^{H}%
\mathcal{YD}}$, and let $A:= R \#_\xi H$. Since the coalgebra $A$ is the
smash coproduct of $R$ with $H$, comultiplication is given by (\ref{eq:
DeltaA}). It is useful to have:
\begin{equation}
\Delta _{A}^{2}\left( r\#h\right) = (r^{(1) } \# r_{\langle
-1\rangle}^{(2)}r_{\langle -2\rangle}^{(3)}h_{(1)} ) \otimes ( r_{\langle
0\rangle }^{(2)}\#r_{\langle -1\rangle }^{(3)}h_{(2)} ) \otimes ( r_{\langle
0\rangle }^{(3)}\#h_{\left( 3\right) } ).  \label{form: DeltaA2}
\end{equation}%
Let $v:=G(\xi)$ as in the preceding sections and let $v_A:= \mho^2(v) $ so that ${v}_{A}\left(
x\#h\otimes y\#h^{\prime }\right) =v\left( x\otimes
hy\right) \varepsilon _{H}\left( h^{\prime }\right) $. Since $v$ is unital, $%
v_A$ is also and thus is a gauge transformation.

In this subsection, we prove that the twisting of $A$ by $v_A$ is the
bosonization of the dual quasi-bialgebra $R^v$ and $H$ defined in Proposition %
\ref{teo: smash quasi YD}, i.e., that $A^{{v}_{A}}=R^{{v}}\#H$. Since in
general ${v}_{A}$ might not be a cocycle, we cannot say that $A^{{v}_{A}}$
is a bialgebra, but $A^{{v}_{A}}$ is always a dual quasi-bialgebra.

\begin{proposition}
\label{teo: smash bilin}For $(R,m,u,\Delta ,\varepsilon ,\xi )$ as
above, let $A:= R \#_\xi H$, let $v:= G(\xi)$ and let $v_A:=
 \mho^2(v)$. Then
\begin{equation*}
A^{{v}_{A}}=R^{{v}}\#H,
\end{equation*}
the bosonization of the dual quasi-bialgebra $\left( R,m^v,u,\Delta
,\varepsilon , \partial _{R}^{{2}} \left( v\right) \right) $ in ${_{H}^{H}%
\mathcal{YD}}$ by $H$.
\end{proposition}

\begin{proof}
Since $u_{A^{{v}_{A}}}(1_{K})=u_{A}(1_{K})=1_{R}\otimes 1_{H}=
1_{R^v}\otimes 1 =u_{R^{{v}}\#H}(1_{K})$, the unit maps for $A^{{v}_{A}}$
and $R^{{v}}\#H$ are the same. It remains to show that the multiplication
maps $m^{v_A} $ and $m_{R^v \#H}$ are the same on $R \#H$, and that the
reassociators coincide.

We begin by computing the product in $A^v$ of two elements from $R \# 1_H$.
As usual, when the context is clear, we omit the subscript $H$ from $1_H$,
and write $m$ instead of $m_R$.
\begin{eqnarray*}
&& m^{{v}_{A}}(r\#1 \otimes s\#1 ) \\
&=& ({v}_{A}\ast m_{A}\ast {v}_{A}^{-1})(r\#1\otimes s\#1) \\
&=& {v}_{A}\left[ (r\#1)_{(1)}\otimes (s\#1)_{(1)}\right] m_{A}\left[
(r\#1)_{(2)}\otimes (s\#1)_{(2)}\right] {v}_{A}^{-1}\left[
(r\#1)_{(3)}\otimes (s\#1)_{(3)}\right] \\
&\overset{(\ref{form: DeltaA2})}{=}& {v}_{A}\left[r^{(1)}\#r_{\langle
-1\rangle }^{(2)}r_{\langle -2\rangle }^{(3)}\otimes s^{(1)}\#s_{\langle
-1\rangle }^{(2)}s_{\langle -2\rangle }^{(3)}\right]m_{A}\left[r_{\langle
0\rangle }^{(2)}\#r_{\langle -1\rangle }^{(3)}\otimes s_{\langle 0\rangle
}^{(2)}\#s_{\langle -1\rangle }^{(3)}\right] {v}_{A}^{-1}\left[r_{\langle
0\rangle }^{(3)}\#1 \otimes s_{\langle 0\rangle }^{(3)}\#1\right] \\
&=&{v}\left[r^{(1)}\otimes r_{\langle -1\rangle }^{(2)}r_{\langle -2\rangle
}^{(3)}s^{(1)}\right ] m_{A}\left[ r_{\langle 0\rangle }^{(2)}\#r_{\langle
-1\rangle }^{(3)}\otimes s^{(2)}\#s_{\langle -1\rangle }^{(3)}\right ] {v}%
^{-1}\left [ r_{\langle 0\rangle }^{(3)}\otimes s_{\langle 0\rangle
}^{(3)}\right ]
\end{eqnarray*}
By {\eqref{form: multi smash xi}},
\begin{equation*}
m_{A}\left [ r_{\langle 0\rangle }^{(2)}\#r_{\langle -1\rangle
}^{(3)}\otimes s^{(2)}\#s_{\langle -1\rangle }^{(3)}\right ] = m \left [
(r_{\langle 0\rangle }^{(2)})^{(1)}\otimes (r_{\langle 0\rangle
}^{(2)})_{\langle -1\rangle }^{(2)}r_{\langle -3\rangle }^{(3)}s^{(2)}\right
]\#\xi ((r_{\langle 0\rangle }^{(2)})_{\langle 0\rangle }^{(2)}\otimes
r_{\langle -2\rangle }^{(3)}s^{(3)})r_{\langle -1\rangle }^{(3)}s_{\langle
-1\rangle, }^{(4)}
\end{equation*}
and so $m^{{v}_{A}}(r\#1 \otimes s\#1 ) $ equals:

\begin{eqnarray*}
&& {v}(r^{(1)}\otimes r_{\langle -1\rangle }^{(2)}r_{\langle -2\rangle
}^{(3)}x^{(1)})m (r_{\langle 0\rangle }^{(2)}\otimes r_{\langle -1\rangle
}^{(3)}x^{(2)})\#\xi (r_{\langle 0\rangle }^{(3)}\otimes x^{(3)})r_{\langle
-1\rangle }^{(4)}s_{\langle -1\rangle }^{(4)} {v}^{-1}(r_{\langle 0\rangle
}^{(4)}\otimes s_{\langle 0\rangle }^{(4)}) \\
&& \quad \text{ where } x:= r_{\langle -2\rangle }^{(4)}s^{(1)} \\
&=& {v}(r^{(1)}\otimes r_{\langle -1\rangle }^{(2)}(r_{\langle -1\rangle
}^{(3)}x^{(1)})^{(1)}) m (r_{\langle 0\rangle }^{(2)}\otimes (r_{\langle
-1\rangle }^{(3)}x^{(1)})^{(2)})\#\xi (r_{\langle 0\rangle }^{(3)}\otimes
x^{(2)})r_{\langle -1\rangle }^{(4)}s_{\langle -1\rangle }^{(3)} {v}^{-1} [
r_{\langle 0\rangle }^{(4)}\otimes s_{\langle 0\rangle }^{(3)} ] \\
&=& ({v}\ast m )(r^{(1)}\otimes r_{\langle -1\rangle }^{(2)}r_{\langle
-3\rangle }^{(3)}s^{(1)})\#\xi (r_{\langle 0\rangle }^{(2)}\otimes
r_{\langle -2\rangle }^{(3)}s^{(2)})r_{\langle -1\rangle }^{(3)}s_{\langle
-1\rangle }^{(3)}{v}^{-1}(r_{\langle 0\rangle }^{(3)}\otimes s_{\langle
0\rangle }^{(3)}) \\
&=& ({v}\ast m )(r^{(1)}\otimes r_{\langle -1\rangle }^{(2)}r_{\langle
-2\rangle }^{(3)}s^{(1)})\#\xi (r_{\langle 0\rangle }^{(2)}\otimes
r_{\langle -1\rangle }^{(3)}s^{(2)})\Psi ({v}^{-1})(r_{\left\langle
0\right\rangle }^{(3)}\otimes s^{(3)}) \\
&=& ({v}\ast m )(r^{(1)}\otimes r_{\langle -1\rangle }^{(2)}s^{(1)})\#\xi
((r_{\langle 0\rangle }^{(2)})^{(1)}\otimes (r_{\langle 0\rangle
}^{(2)})_{\left\langle -1\right\rangle }^{(2)}s^{(2)})\Psi ({v}%
^{-1})((r_{\langle 0\rangle }^{(2)})_{\left\langle 0\right\rangle
}^{(2)}\otimes s^{(3)}) \\
&=& ({v}\ast m )(r^{(1)}\otimes r_{\langle -1\rangle }^{(2)}s^{(1)})\#\left[
\xi \ast \Psi ({v}^{-1})\right] (r_{\langle 0\rangle }^{(2)}\otimes s^{(2)})
\\
&=&({v}\ast m )(r^{(1)}\otimes r_{\langle -1\rangle }^{(2)}s^{(1)})\#{v}%
^{-1}(r_{\langle 0\rangle }^{(2)}\otimes s^{(2)}) \\
&=&({v}\ast m \ast {v}^{-1})(r\otimes s)\#1=m^{v}(r\otimes s)\#1.
\end{eqnarray*}%
Since $\Delta _{A}$ is $H$-bilinear, and both $v_{A}$ and $m_{A}$ are $H $%
-bilinear $H$-balanced, $m^{{v}_{A}} $ is also $H$-bilinear $H$-balanced so
that:
\begin{equation*}
m^{{v}_{A}} (r\#h\otimes s\#l)=m^{{v}_{A}}(r\#1 \otimes h_{(1)}s\#1
)h_{(2)}l=m^v (r\otimes h_{(1)}s)\#h_{(2)}l=m_{R^{{v}}\#H}(r\#h\otimes s\#l),
\end{equation*}
and $A^{{v}_{A}}=R^{{v}}\#H$.

It remains to check that the reassociators are the same. Set $\alpha
:=\partial _{R}^{2}\left( v\right) .$ Since $A$ is a dual quasi-Hopf algebra
with trivial reassociator $\varepsilon$, by (\ref{form: reassociator}), the
reassociator $\gamma $ for $A^{{v}_{A}}$ is:
\begin{eqnarray*}
\gamma &:=&\left( \varepsilon _{A}\otimes {v}_{A}\right) \ast {v}_{A}\left(
A\otimes m_{A}\right) \ast {v}_{A}^{-1}\left( m_{A}\otimes A\right) \ast
\left( {v}_{A}^{-1}\otimes \varepsilon _{A}\right) \\
&=&\partial _{A}^{2}\left( v_{A}\right) =\partial _{A}^{2}\mho
_{H,R}^{2}\left( v\right) =\mho _{H,R}^{3}\partial _{R}^{2}\left( v\right)
=\mho _{H,R}^{3}\left( \alpha \right)
\end{eqnarray*}%
so that%
\begin{eqnarray*}
\gamma \left( r\#h\otimes s\#l\otimes t\#k\right) &=&\mho _{H,R}^{3}\left(
\alpha \right) \left( r\#h\otimes s\#l\otimes t\#k\right) \\
&=&\alpha \left( r\otimes h_{\left( 1\right) }s\otimes h_{\left( 2\right)
}lt\right) \varepsilon _{H}\left( k\right) =\alpha _{R^{{v}}\#H}\left(
r\#h\otimes s\#l\otimes t\#k\right) .
\end{eqnarray*}
\end{proof}

\section{The main theorem}\label{sect: main}

Recall from Section \ref{sec: splitting data}, that if $(A,H,\pi,\sigma)$ is a
splitting datum, then there is associated to this datum a
pre-bialgebra with cocycle $(R,\xi)$. Furthermore, given a
pre-bialgebra with cocycle $(R,\xi)$, then one can construct a
splitting datum $(R \#_\xi H, H, \pi, \sigma)$ and $R \#_\xi H \cong
A$. We can now prove the main result of this paper. Recall our assumption  that $H$ has an ad-invariant integral; for example $H$ could be a group algebra.

\begin{theorem}
\label{teo: splitting}Let $\left( A,H,\pi ,\sigma \right) $ be a
splitting datum with associated
pre-bialgebra with cocycle $(R,\xi)$ so that $A \cong B:= R \#_\xi
H$. Suppose $\sigma(H) =Corad\left( A\right) .$ Then $A$ is
quasi-isomorphic to the bosonization of a connected dual quasi-bialgebra in $%
{_{H}^{H}\mathcal{YD}}$ by $H$.
\end{theorem}

\begin{proof}
Let $(R,m,u,\Delta ,\varepsilon ,\xi )$ be the pre-bialgebra with cocycle in
$_{H}^{H}\mathcal{YD}$ associated to $(A,H,\pi ,\sigma ).$ Since $\sigma
\left( H\right) =Corad\left( A\right) $, by \cite[Corollary 5.3.5]{Mo}, $%
Corad\left( R\right) \subseteq \tau \left( Corad\left( A\right) \right)
\subseteq \tau \left( \sigma \left( H\right) \right) \subseteq K$ where $%
\tau(a) = a_{(1)} \sigma   S_H  \pi(a_{(2)})$ as described in
Section \ref{sec: splitting data}. Thus $\mathrm{Corad}\left(
R\right) =K$ whence $R$ is connected. Since $H$ has an
$ad$-invariant integral, by Theorem \ref{teo: hope} we have the
datum $(R,m^v,u,\Delta ,\varepsilon ,{v})$ and by Proposition
\ref{pro: R dual quasi}, $\left( R,m^v,u,\Delta ,\varepsilon ,\alpha
\right) $ is a dual quasi-bialgebra in the braided monoidal
category ${_{H}^{H}\mathcal{YD}}$. Set $B:=R\#_{\xi }H.$ By
Proposition \ref{teo: smash bilin}, there exists a gauge
transformation ${v}_{B}:B\otimes B\rightarrow K$ such that
\begin{equation*}
B^{{v}_{B}}=R^{{v}}\#H.
\end{equation*}%
where the latter is the bosonization of the dual quasi-bialgebra
$\left( R,m^v,u,\Delta ,\varepsilon ,\alpha \right) $ in
${_{H}^{H}\mathcal{YD}}$ by $H$.  In conclusion  $A$ is
quasi-isomorphic to the bosonization of the connected dual
quasi-bialgebra $\left( R, m^v,u,\Delta ,\varepsilon ,\alpha \right)
$ in ${_{H}^{H}\mathcal{YD}}$ by $H$.
\end{proof}

We can now give the proof of the main theorem.

\begin{proof}[Proof of Theorem \protect\ref{Theo1}]
By \cite[Theorem 2.35]{A.M.S.}, $A$ fits into a splitting datum
$\left( A,H,\pi ,\sigma \right) $ where $\sigma :H\rightarrow A$ is
the canonical inclusion.
  As mentioned in Section \ref{sec:preliminaries},
since $H$ is semisimple and cosemisimple, it has an $ad$-invariant
integral so that Theorem \ref{teo: splitting} applies.
\end{proof}

\begin{remark}Let $A$ be a bialgebra whose coradical $H$ is a subbialgebra of $A$ with antipode. Akira Masuoka
pointed out, see \cite{Masuoka-communic}, that, by Takeuchi's lemma \cite[%
Lemma 5.2.10]{Mo}, $A$ is necessarily a Hopf algebra. Thus $A$ is a
Hopf algebra with the  dual Chevalley property.

%Recall that a bialgebra $A$ has the dual Chevalley property whenever
%the coradical $H$ of $A$ is a subbialgebra of $A$ with antipode.
%Akira Masuoka
%pointed out, see \cite{Masuoka-communic}, that, by Takeuchi's lemma \cite[%
%Lemma 5.2.10]{Mo}, $A$ is necessarily a Hopf algebra.

\end{remark}

\smallbreak
\begin{center}Acknowledgement\end{center} We would like to thank the referee for an extremely helpful report, and especially for pointing out that our main result
could be explained more conceptually by the use of Diagram \eqref{diagram main}.

\end{document}